\documentclass[english]{amsart}

\usepackage{babel}
\usepackage{amstext}
\usepackage{amsmath}
\usepackage{amsfonts}
\usepackage{latexsym}
\usepackage{ifthen}
\usepackage{xypic}
\xyoption{all}
\newcommand{\I}{{\mathcal I}}
\newcommand{\J}{{\mathcal J}}
\newcommand{\E}{{\mathcal E}}

\newcommand{\PN}{{\mathbb P}}

\newcommand{\PF}{{\mathbb F}}

\newcommand{\codim}{{\rm codim}}
\newcommand{\Pic}{{\rm Pic}}

\newcommand{\lra}{\longrightarrow}
\newcommand{\KC}{{\mathbb C}}
\newcommand{\KZ}{{\mathbb Z}}
\newcommand{\KQ}{{\mathbb Q}}
\newcommand{\KN}{{\mathbb N}}

\newcommand{\Bs}{{\rm Bs}}
\newcommand{\Sl}{\rm Sl}

\newcommand\sE{{\mathcal E}}
\newcommand\sF{{\mathcal F}}
\newcommand\sG{{\mathcal G}}

\newcommand\sI{{\mathcal I}}

\newcommand\sO{{\mathcal O}}

\newcommand\bZ{{\mathbb Z}}

\newcommand\bQ{{\mathbb Q}}
\newcommand\bN{{\mathbb N}}

\newcommand\bP{{\mathbb P}}

\newcounter{lemma}

\newtheorem{lemma1}[lemma]{\setcounter{equation}{0}}

\newenvironment{lemma}{\begin{lemma1}{\bf Lemma.}}{\end{lemma1}}

\newenvironment{example}{\begin{lemma1}{\bf Example.}\rm}{\end{lemma1}}
\newenvironment{abs}{\begin{lemma1}\rm}{\end{lemma1}}
\newenvironment{theorem}{\begin{lemma1}{\bf Theorem.}}{\end{lemma1}}

\newenvironment{proposition}{\begin{lemma1}{\bf Proposition.}}{\end{lemma1}}
\newenvironment{proposition2}[1]{\begin{lemma1}{\bf Proposition [#1].}}{\end{lemma1}}
\newenvironment{corollary}{\begin{lemma1}{\bf Corollary.}}{\end{lemma1}}

\newenvironment{notation}{\begin{lemma1}{\bf Notation.}}{\end{lemma1}}

\newenvironment{setup}{\begin{lemma1}{\bf Setup.}}{\end{lemma1}}

\begin{document}

\title {Threefolds with big and nef \\ anticanonical bundles I} 
\author{Priska Jahnke, Thomas Peternell, Ivo Radloff}
\date{\today, {\em Address:}  Mathematisches Institut, Universit\"at Bayreuth, D--95440 Bayreuth/Germany, {\em e-mail:} priska.jahnke@uni-bayreuth.de, ivo.radloff@uni-bayreuth.de, thomas.peternell@uni-bayreuth.de}
\maketitle

%%%%%%%%%%%%%%%%%%%%%%%%%%%%

\section*{Introduction}

As one of the first applications of Mori theory, Mori and Mukai classified (smooth) Fano threefolds with Picard (or second Betti) number at least $2.$ 
In differential geometric terms, this is the same as classifying smooth threefolds with positive Ricci curvature. It is clearly interesting to consider
the situation when we ``degenerate'' the positivity condition, i.e. we consider threefolds whose anticanonical bundles are no longer ample but only 
big and nef. E.g. there exists a metric with semipositive Ricci curvature which is positive at some point. Recall that $-K_X$  nef is to say that 
$(-K_X) \cdot C \geq 0$ for all curves $C \subset X.$ This automatically implies $(-K_X)^3 \geq 0$ and bigness is just saying that $(-K_X)^3 > 0.$ 
We will call a threefold $X$ with $-K_X$ big and nef (but not ample) an {\it almost Fano threefold}.
\vskip .2cm \noindent 
With this paper we begin the classification of smooth almost Fano threefolds $X$. Due to the complexity of the problem we first study almost Fano 
threefolds
with Picard number two. The classification in the Fano case uses essentially the fact that there are two Mori contractions (= contractions of an
extremal ray) which are 
transversal in some sense. In our case we only have one Mori contraction at our disposal. The second Mori contraction is substituted by the morphism
associated with the base point free linear system $|{-}mK_X|, m \gg 0$. This morphism is clearly more
difficult to handle than a ``simple'' Mori contraction.  
The present article is dealing with the case that $|{-}mK_X|, m \gg 0$, is  {\it divisorial}, while the
second part will treat the case where this morphism is small, i.e. contracts just finitely many curves. 
\vskip .2cm \noindent
To be  a little more precise, the setup of the paper is as follows. Call the extremal ray contraction $\phi: X \lra Y$; on the other hand, the base point free theorem guarantees that $|{-}mK_X|$ is spanned for $m \gg 0$. After Stein factorization we get a second map $\psi: X \lra X'$ and a diagram
\[\xymatrix{X \ar[r]^{\psi} \ar[d]^{\phi} & X'\\
            Y. &}\]
Here $X'$ is a canonical Gorenstein Fano threefold, i.e., $-K_{X'}$ is ample, but $X'$ is singular with mild singularities. Using Mori's classification of what $\phi: X \lra Y$ can be, $\phi$ is a
either a del Pezzo fibration over $\bP_1$, a conic bundle over $\bP_2 $ or birational with a very precise structure. We shall treat all these
cases separately.  
\vskip .2cm \noindent The final result being a quite long list, we will not reproduce it here but just point out the places in the paper where
the explicit classification results occur: these are the theorems 2.2,
2.5, 2.6, 2.8, 2.9, 3.2, 3.4., 3.13, 4.9, 4.11, 5.2. 
\vskip .2cm \noindent It is interesting to note that not all almost Fano threefolds occur as degenerations of Fano threefolds, instead there
occur completely new cases.  
\vskip .2cm \noindent
The authors want to thank the DFG Schwerpunkt ``Global Methods in Complex Geometry'' for significant and indispensable support of our project. 
\tableofcontents

%%%%%%%%%%%%%%%%%%%%

\section{Preliminaries} \label{prel}
\setcounter{lemma}{0}

\begin{notation} {\rm Let $X$ be a smooth projective threefold with $-K_X$ big and nef, called {\em almost Fano}. We always assume that $X$ is {\it not Fano}. Moreover we assume $\rho (X) = 2$. Then $-mK_X$ will be spanned for suitable large $m.$ Throughout this paper $\psi: X \to X'$ will denote the morphism (with connected fibers) associated with $\vert -mK_X \vert$. We suppose that $\psi$ is not small; by our assumption this means that $\psi$ contracts an irreducible divisor which will be denoted by $D.$ Notice that $X'$ is Gorenstein with only canonical non-terminal singularities. In accordance with the literature, we define the {\em genus} $g$ of $X$ by  $g+2 = h^0(X, -K_X)$ and call a general member $S \in |{-}K_X|$ a {\em general elephant}.}
\end{notation} 

The existence of a $K3$-elephant was proved by Shokurov for the smooth case in \cite{Shokurov}, and by Reid in \cite{Reid} for Gorenstein Fano threefolds with canonical singularities.

\begin{proposition} 
 Let $X$ be a smooth almost Fano threefold. Then 
$$ g+2 = h^0(X, -K_X) = \chi(X,-K_X) = {{(-K_X)^3} \over {2}} + 3$$
and $(-K_X)^3$ is even (Riemann-Roch).
\end{proposition}

\begin{proposition} 
If $-K_X$ is spanned, the map $X \to W \subset \bP_{g+1}$ associated to the linear system $|{-}K_X|$ factorises as $\mu \circ \psi$ with a finite morphism 
 \[\mu: X' \lra W \subset \PN_{g+1}\] 
which has degree at most $2$ and actually it has degree $1$ unless
$({-}K_X)^3 \leq 8$ (the $X'$ where $\mu$ has degree $2$ are called hyperelliptic). If $\deg \mu = 1,$ then $\mu$ is an embedding.
\end{proposition} 

\begin{proof} The paper [IP99], Proposition~2.1.15 shows that $\deg \mu \leq 2.$ For the hyperelliptic case see Proposition~\ref{hyp} below. The last statement follows from Mumford's criterion \cite{Mumford}, considering a general elephant.
\end{proof}

As in Fano classification, Gorenstein Fano threefolds $X'$ with canonical singularities where $|{-}K_{X'}|$ is not base-point free or not very ample will occur as exceptional cases. We have the following classification from \cite{genprep}: 

\begin{theorem} \label{gen}
Let $X'$ be a Gorenstein Fano threefold with canonical singularities. Assume $|{-}K_{X'}|$ is not base point free. Then one of the following holds:
\begin{enumerate}
  \item $\dim(\Bs|{-}K_{X'}|) = 0$. In this case $X'$ is a complete intersection of a cone over a quadric in $\PN_3$ and a general sextic in the weighted projective space $\PN(1^4,2,3)$ and $\Bs|{-}K_{X'}| = X'_{sing}$ is a single, terminal singularity.
 \item $\dim(\Bs|{-}K_{X'}|) = 1$. Then $\Bs|{-}K_{X'}| \simeq \PN_1$ and either 
  \begin{enumerate}
    \item $X'$ is the blowup of a sextic in $\PN(1^3,2,3)$ along an irreducible curve of arithmetic genus one or
    \item $X' \simeq S \times \PN_1$, where $S$ is a del Pezzo surface of degree $1$ with at most Du Val singularities or
    \item $X' =X'_{2m-2}$ is an anticanonical model of the blowup of
      the variety $U_m$ (see below) along a smooth, rational complete intersection curve $\Gamma_0 \subset U_{m, \mathrm{reg}}$ for $3 \le m \le 12$.
 \end{enumerate}
\end{enumerate}
\end{theorem}

\noindent Here $U_m$ denotes a double cover of $\PN(\sO_{\PN_1}(m)
\oplus \sO_{\PN_1}(m-4) \oplus \sO_{\PN_1})$ with at worst canonical
singularities, such that ${-}K_{U_m}$ is the pullback of the
tautological line bundle $\sO(1)$. For $m \ge 4$, this is a hyperelliptic
Gorenstein almost Fano threefold of degree $4m-8$. The curve $\Gamma_0$ lies over the complete intersection of some general
element in $|\sO(1)|$ and the ``minimal surface'' $B \in |\sO(1) -mF|$,
where $|F|$ denotes the pencil (note that $\Gamma_0$ is always
contained in the ramification locus). If $m = 3$, then $\Gamma_0$ is
the only curve, on which ${-}K_{U_3}$ is not nef. For details of the construction see \cite{genprep}, section~5.

The threefolds in (i) and (ii) are the expected degenerations of the smooth case. New in a sense are (1) (see \cite{Mella}) and (iii). The intersection in (1) is a degeneration of $V_2$, a double cover of $\PN_3$ ramified along a smooth sextic, but here $|{-}K_{V_2}|$ is spanned. The examples in (iii) are not $\KQ$-factorial and have one isolated canonical singularity. We therefore get 

\begin{corollary} \label{gencor}
Let $X$ be a smooth almost Fano threefold with $\rho(X) = 2$, such that $X$ is not Fano and the anticanonical map $\psi$ is divisorial. Then $|{-}K_X|$ is base point free.
\end{corollary}

\begin{proof}
By assumption, the image $X'$ of the morphism $\psi$, defined by $|{-}mK_X|$ is a $\KQ$-factorial Gorenstein Fano threefold with canonical, non-terminal singularities and Picard number one.
\end{proof}

\vspace{0.2cm}

The classification of hyperelliptic Fano threefolds is due to Iskovskikh in the smooth case (\cite{Isk}) and to Cheltsov, Shramov and Przyjalkowski for Gorenstein Fano threefolds with canonical singularities (\cite{Cheltsov}). In our situation we have  

\begin{proposition} \label{hyp}
Let $X$ be a smooth almost Fano threefold with $\rho(X) = 2$, such that the anticanonical map is divisorial. If the anticanonical model $X' \stackrel{2:1}{\lra} W$ is hyperelliptic, then we are in one of the following cases 
 \begin{enumerate} 
   \item $(-K_X)^3 = 2$, $W = \PN_3$ and $X' \to W$ is ramified along a sextic; 
   \item $(-K_X)^3 = 4$, $W \subset \PN_4$ is a quadric and $X' \to W$ is ramified along a quartic;
   \item $(-K_X)^3 = 8$, $W$ is the cone in $\PN_6$ over the Veronese surface in $\PN_5$ and $X' = X_6 \subset \PN(1^3,2,3)$, i.e., $X' \to W$ is ramified along a cubic.
  \end{enumerate}
\end{proposition}
\noindent All of these threefolds are the expected degenerations of Iskovskikh's list. 

\begin{proof}
We have $X \stackrel{\psi}{\lra} X' \stackrel{\mu}{\lra} W \subset
\PN_{g+1}$ with $\psi$ divisorial and $\mu$ a double cover of a
threefold $W$ of minimal degree, i.e. $\deg(W) = \codim(W)
+1$. Varieties with this property are classified by Bertini
(\cite{Bertini}), whereafter $W$ is one of $\PN_3$, a quadric in
$\PN_4$, the cone over the Veronese surface or a (cone over a)
rational scroll. By a cone over a (rational) scroll we mean the image of 
 \[\PF(d_1, d_2, d_3) = \PN(\sO_{\PN_1}(d_1) \oplus \sO_{\PN_1}(d_2) \oplus \sO_{\PN_1}(d_3)), \quad d_1 \ge d_2 \ge d_3 \ge 0\]
in $\PN_{d_1 + d_2 + d_3 + 2}$ under the map associated to the
tautological system $|\sO(1)|$. Note that this is free. For $d_3 \ge
1$, the cone and $\PF(d_1, d_2, d_3)$ are
isomorphic. The pencil of $\PF(d_1, d_2, d_3)$ is denoted by $|F|$. If
in our situation $W$ is a cone over a scroll, then $d_3 = 0$, since $\rho(X') = 1$. 

1.) If $d_2 = 0$, then $W$ is a double cone over a rational normal curve of degree $d_1$. The double cover $X'$ will have canonical singularities along a curve. There are dissident singular points, which are not cDV, hence a resolution of singularities will not have Picard number two. 

2.) Assume $d_2 > 0$, i.e. $W$ is a cone over a Hirzebruch surface. In this case, the map $\sigma\colon \PF(d_1, d_2, 0) \to W$ is a small resolution and may be viewed as blowup of $W$ along the Weil divisor $\sigma(F)$. Since $X$ is smooth, we obtain an induced birational map $X \to \PF(d_1, d_2, 0)$, mapping the exceptional divisor $D$ of $\psi$ to the curve in $\PF(d_1, d_2, 0)$ contracted by $\sigma$. This contradicts $\rho(X) = 2$.
\end{proof}

\vspace{0.2cm}

For small genus we find in our situation (see \cite{AG5}, Proposition~4.1.12. for the smooth case):

\begin{proposition} \label{ci}
Let $X$ be a smooth almost Fano threefold with $\rho(X) = 2$, such that the anticanonical map $\psi: X \to X'$ is divisorial. Assume $X'$ is not hyperelliptic.
 \begin{enumerate}
  \item If $g = 3$, then $X'_4 \subset \PN_4$ is a quartic.
  \item If $g = 4$, then $X'_{2,3} \subset \PN_5$ is a complete intersection of a quadric and a cubic.
  \item If $g = 5$, then $X'_{2,2,2} \subset \PN_6$ is a complete intersection of three quadrics.
 \end{enumerate}
\end{proposition}

\begin{proof}
Since the canonical curve section $C \subset X'$ is a smooth canonical curve of genus $g$, (1) and (2) are easily obtained. 

Assume $g = 5$. We have two possible cases: either $X'$ is cut out by quadrics or it is trigonal. Since $X'$ is already a complete intersection in the first case, assume the latter one. Then by \cite{Cheltsov}, $X'$ is the anticanonical model of an almost Fano threefold $V$ with canonical singularities, where $V$ is a divisor in $|\sO(3) + \pi^*\sO_{\PN_1}(-1)|$ on the projective bundle $\pi: \PN(\oplus_{i=1}^4\sO_{\PN_1}(d_i)) \to \PN_1$, where either $d_1 = d_2 = d_3 = 1$, $d_4 = 0$ or $d_1 = 2$, $d_2 = 1$, $d_3 = d_4 = 0$. Neither case is possible in our situation: in the first case the map $V \to X'$ is small and in the second case the Picard number of $X'$ is greater than one.
\end{proof}

By definition, the anticanonical map $\psi$ is a birational contraction. The assumption $\rho(X) = 2$ guarantees that $\psi$ is primitive, i.e. it does not factor. The structure of the exceptional locus of such contractions is studied by Wilson, Paoletti and Minagawa. We have in our situation

\begin{proposition2}{Wilson, Paoletti, Minagawa}\label{Wilson}
Let $X$ be a smooth almost Fano threefold with $\rho(X) = 2$, such that the anticanonical map $\psi: X \to X'$ is divisorial, contracting the divisor $D$ to a curve $B$. Let $l_{\psi}$ be a general exceptional fiber. Then 
\begin{enumerate}
 \item $B$ is a smooth curve of cDV singularities, $X = Bl_B(X')$;
 \item $D$ is a conic bundle over $B$; each fiber is either isomorphic to a smooth conic or to a line pair. In particular, $D.l_{\psi} = -2$.
\end{enumerate}
\end{proposition2} 

\begin{proof}
\cite{Wilson}, \cite{Wilson2}, \cite{Paoletti} and \cite{Minagawa}.
\end{proof}

\

%%%%%%%%%%%%%%%%%%%%%%%%%

\section{Del Pezzo fibrations} \label{delPezzo}
\setcounter{lemma}{0}

In this section we consider threefolds with $-K_X$ big and nef admitting a del Pezzo fibration.

\begin{setup} {\rm We fix for this section the following setup. $X$ is
    almost Fano, i.e. a smooth projective threefold with $-K_X$ big
    and nef, but not ample. As usual $\rho (X) = 2$ and we suppose
    that $\psi$ is divisorial. Suppose that $\phi: X \to \bP_1 $ is a del Pezzo fibration, which is the contraction of an extremal ray, since $\rho (X) = 2.$ Let $F$ denote a general fiber of $\phi.$
Notice that $K_F^2 \ne 7$ and that $F$ is normal [Mo82]. We put $F' = \psi(F) $ and $F'' = \mu(F').$ }
\end{setup}

\begin{theorem} Suppose that $\phi$ is a $\bP_2-$bundle and write $X = \bP(\sE)$ with a rank 3-bundle $\sE.$ Normalize $\sE$ in the 
following way: $\sE = \sO(a_1) \oplus \sO(a_2) \oplus \sO$ with $a_1 \geq a_2 \geq 0.$ Then either
\begin{enumerate} 
\item $a_1 = a_2 = 1$, and $\psi $ contracts the unique section $C_0$ with normal bundle $N_{C_0} = \sO(-1) \oplus \sO(-1)$ (but here 
$\psi$ is small), or
\item $a_1 = 2, a_2 = 0$, and $\psi $ contracts the divisor $D = \bP(\sO \oplus \sO) \subset X$ (Case A.2,no.1).
\end{enumerate}
In both cases $-K_X$ is indeed nef. 
\end{theorem}

\begin{proof} The condition that $-K_X$ is nef, but not ample, can be rewritten as follows: 
$$ \sE \otimes {{\det \sE^*} \over {3}} \otimes \sO({{2} \over {3}}) = \sE \otimes \sO({{-a_1-a_2+2} \over {3}})$$
is nef but not ample. Hence $a_1+a_2 = 2$ and we are in one of the both cases stated in the theorem. 
The bigness of $-K_X$ is translated into 
$$ c_1(\sE \otimes {{\det \sE^*} \over {3}} \otimes \sO({{2} \over {3}})) > 0 $$
and is therefore automatically fulfilled. 
The rest of the claim is clear.
\end{proof}

\noindent
{\it From now on we shall assume that $F \ne \bP_2$ (for some or - equivalently - all fibers)}. \\

\begin{setup} {\rm  Consider the unique irreducible prime divisor $D$ contracted by $\psi$. 
Since $-K_X$ is $\phi-$ample, $K_X \vert D \not \equiv 0, $ hence the divisor 
$D$ cannot be contracted to a point
by $\psi,$ and we denote $B = \psi (D).$ 
Let $l_{\psi}$ be the general fiber of $\psi.$ Then by (\ref{Wilson})
$$ D \cdot l_{\psi} = -2  \eqno (2.3.1)$$
and $l_{\psi}$ is a conic in $\bP_2.$ Since $\rho(X/X') = 1,$ this must be true for all fibers of $\psi.$ 
Now write
$$ D = -\alpha K_X - \phi^*(\sO(\beta)) = -\alpha K_X - \beta F, \eqno (2.3.2)$$
with rational numbers $\alpha$ and $\beta$ and $F$ a fiber of $\phi.$ 
Thus $D_F \in \vert -\alpha K_F \vert$, hence $\alpha > 0$ and
actually $\alpha$ is an integer, unless possibly $F$ is a
quadric. This comes from the existence of $(-1)-$curves in del Pezzo surfaces $F$ with $K_F^2 \leq 7.$ If $F$ is a quadric, then we can at least say that $\alpha \in {{1} \over {2}} \bZ.$
This special case that $\alpha$ is not an integer is excluded from now on and will be treated in (2.9) separately. \\
Since now $\alpha$ is an integer, also $\beta $ is an integer and from (2.3.1) and (2.3.2) we defer
$$ \beta = 1, F \cdot l_{\psi} = 2 \ {\rm or} \ \beta = 2, F \cdot l_{\psi} = 1.  \eqno (2.3.3) $$
We finally notice that $K_X^2 \cdot D = 0$ gives
$$ \alpha (-K_X)^3 =  \beta K_F^2.  \eqno (2.3.4) $$
}
\end{setup}

\begin{proposition} Suppose in (2.3) that $\beta = 2.$ Then either $\alpha = 1$ or $(-K_X)^3 = 2, \alpha = 3$ or $(-K_X)^3 = 4, \alpha = 2.$ 
\end{proposition}

\begin{proof} We will compute $g = g(B)$ in two different ways. Comparing both formulas, we will then arrive at our claim. 
We notice that by (2.3.4): $ \alpha (-K_X)^3 = 2K_F^2,$ in particular $(-K_X)^3 \leq 16.$ \\
(1) We will use the exact sequence
$$ 0 \to H^0(X',\sI_B \otimes -\alpha K_{X'}) \to H^0(X',-\alpha K_{X'}) \to H^0(B,-\alpha K_{X'} \vert B) \to $$
$$ \to H^1(X',\sI_B \otimes -\alpha K_{X'}) \to 0.$$
Notice
$$h^0(X',\sI_B \otimes -\alpha K_{X'}) = h^0(X,-\alpha K_X - D) = h^0(X,2F) = 3.$$
Next we observe that 
$$H^q(X',\sI_B \otimes -\alpha K_{X'}) = 0  \eqno (*)$$
for $q \geq 1.$ In fact, $$H^q(X,-\alpha K_X - D) = H^q(X,2F) = 0$$ for $q \geq 1$ 
and thus the apparent vanishing
$$R^q\psi_*(-\alpha K_X - D) = R^q \psi_*(-D) \otimes (-\alpha K_{X'}) = 0$$ 
for $q \geq 1$ yields the vanishing (*). 
Next we compute
$$ h^0(B,-\alpha K_{X'} \vert B).$$
First of all, $-K_{X'} \cdot B = -K_X \cdot D \cdot F$ so that
$$ -K_{X'} \cdot B = \alpha K_F^2.  \eqno (A)$$
Thus $$ \chi(B,-\alpha K_{X'} \vert B) = 1 - g + \alpha^2K_F^2 = {{\alpha^3} \over {2}} (-K_X)^3 + 1-g.$$
Now by Kodaira vanishing and (*)
$$ H^1(B,-\alpha K_{X'} \vert B) = H^2(X',\sI_B \otimes -\alpha K_{X'}) = 0 $$
so that 
$$ h^0(B,-\alpha K_{X'} \vert B) = {{\alpha^3} \over {2}} (-K_X)^3+1-g.$$
We arrive finally at
$$ h^0(-\alpha K_{X'}) = 4 - g + {{\alpha^3} \over {2}} (-K_X)^3.$$
Putting in Riemann-Roch and Kodaira vanishing to compute $h^0(-\alpha K_X)$ we obtain
$$ g = ({{\alpha^3} \over {3}} - {{\alpha^2} \over {4}} - {{\alpha} \over {12}}) (-K_X)^3  - 2\alpha + 3. \eqno (**)$$ 
(2) On the other hand, $B \simeq D \cdot F$ since $\psi \vert F$ is an
isomorphism $(F\cdot l_{\psi} = 1)$. 
Thus adjunction and $D_F \in \vert -\alpha K_F \vert $ gives
$$ 2g - 2 = \alpha(\alpha-1)K_F^2 $$
so that
$$ g = {{\alpha^2(\alpha-1)} \over {4}} (-K_X)^3 + 1. \eqno (***)$$
Putting (**) and (***) together, we obtain
$$ ({{\alpha^3} \over {12}} - {{\alpha} \over {12}}) (-K_X)^3 - 2\alpha + 2 = 0. \eqno (+) $$
Of course (+) is always fulfilled if $\alpha = 1.$ If however $\alpha \geq 2$ we have only the two solutions 
stated above; here of course we also use $\alpha(-K_X)^3 = 2K_F^2.$  
\end{proof}

%2.5 
\begin{theorem}
In (2.4) the exceptional cases that $\alpha = 2$ and $(-K_X)^3 = 4$ resp. $\alpha = 3, (-K_X)^3 = 2$ really occur and
are described as follows: 
\begin{enumerate} 
\item let $B \subset  \bP_4$ be a smooth complete intersection of three quadrics and let $\pi: \hat P \to \bP_4$ be the 
blowup of $\bP_4$ along $B$ with exceptional divisor $E.$ Then $X$ is a smooth member of $\pi^*(\sO(4)) - 2E;$ here $(-K_X)^3 = 4,$ (A.2,no.2);
\item $X$ is a degree 2 covering over the 3-dimensional quadric blown up in the intersection of two quadrics; the ramification being in two fibers
of the natural del Pezzo fibration of the blown up quadric; again $(-K_X)^3 = 4,$ (A.2,no.2);
\item $X$ is a degree 2 covering over $\bP_3$ blown up in the smooth intersection of two cubics; the ramification being in two fibers
of the natural del Pezzo fibration of the blown up projective space; here $(-K_X)^3 = 2,$ (A.2,no.3).
\end{enumerate} 
\end{theorem}

\begin{proof}  We consider the two cases separately: 
\begin{enumerate}
\item $\alpha = 2, (-K_X)^3 = 4,$ or
\item $\alpha = 3, (-K_X)^3 = 2.$ 
\end{enumerate}
(1) Here we obtain $g(B) = 5$ and $\deg B = 8$ and $B$ is a complete intersection of three quadrics. \\
First suppose that $\mu: X' \to \bP_4$ is an embedding.
The converse construction is obvious: 
we take $B$ to be a complete intersection of three smooth quadrics in $\bP_4$.
Then $B$ has genus 5 and we let $\hat P$ be the blowup of $\bP_4$ along $B.$ Since $h^0(\bP_4,\sI_B(2)) = 3,$ the Fano 4-fold $\hat P$
carries a contraction $f: \hat P \to \bP_2$. If $Q$ is a smooth quadric containing $B$, then its strict transform $\hat Q \subset  \hat P$
is a Fano threefold and $f \vert Q: Q \to \bP_1$ is the second contraction whose fiber $F$ has $K_F^2 = 4.$ Now let $X$ be a smooth member
of $\vert f^*(\sO(2)) \vert.$ We have $D = -2K_X - 2F$ so that $\alpha = 2$ and $\beta = 2.$ \\
The other case is that $\mu$ is a degree 2 covering over a quadric $Q$. In that case the image $B' = \mu(B) \subset  Q \subset \bP_4$
is a smooth curve of genus 5 and degree 8 (since $B'$ cannot be in the singular locus of $Q,$ we must have $\deg \mu \vert B = 1$).
From $D = -2K_X - 2F$ and the numerical data, we easily obtain
$$ \phi_*(-K_X) = \sO^5 $$
(see the proof of (2.6) for some detailed computations). Since $H^0(-K_X-F) = 0,$ we also have $H^0(-K_X) = H^0(-K_F).$ 
Thus $X$ embeds in $\bP(\phi_*(-K_X)) = \bP_1 \times \bP_4.$ 
The covering part of $X \to Q_3$ is completely induced by $\bP_1;$ in other words, there exists a del Pezzo fibration $Z \to C = \bP_1$
and a degree 2 covering $\bP_1 \to C$ inducing $\phi$ by base change. Moreover we have a birational map $Z \to Q_3.$ \\
Conversely, let $\pi: \hat Q \to Q_3$ be the blowup of $B'$ in $Q_3.$ Then $\hat Q$ is Fano and carries a del Pezzo fibration $p: \hat Q \to \bP_1$ 
such that $K_{F'}^2 = 4$ for the fibers $F'$ of $p.$ Observe that $p$ is defined by $\pi^*(\sO(2) - E)$, where $E$ is the exceptional divisor.
Let
$$R \in \vert p^*(\sO(2))  \vert = \vert \pi^*(\sO(-2)) \otimes -2K_{\hat Q} \vert $$
be general and let $\tau: X \to \hat Q$ be the degree 2 covering branched along $R.$ Then $-K_X = \tau^* \pi^*(\sO_{Q_3}(1))$ and $X$ fulfils
all requirements.
\vskip .2cm \noindent
(2) Here we have a holomorphic map $\mu: X' \to \bP_3$ defined by $\vert -K_{X'} \vert$ of degree $2.$
Let $R \subset  \bP_3$ be the ramification locus. Since $K_{X'} = \mu^*(K_{\bP_3} + R/2),$ we must have 
$\deg R \leq 6.$ If $\deg R  \ne 6, $ then $-K_X$ would be divisible, contradicting $K_F^2 = 3.$ 
Let $R' = \mu^{-1}(R),$ of course $R' \simeq R.$ 
Then $B \subset R', $ since $X'$ is singular along $B.$ Then the formula (***) from the proof of (2.4)
yields $g(B) = 10$ and formula (A) gives $\deg B = 9$ in $\bP_3.$ \\
Conversely, let $B \subset \bP_3$ be the smooth complete intersection of two cubics and $\hat P$ the blowup of $B.$ Then $\hat P$
is Fano with a del Pezzo fibration $p: \hat P \to \bP_1$ defined by $\pi^*(\sO(3)) - E$. Now let $X$ be the degree 2 covering ramified
along a smooth element of $\vert p^*(\sO(2)) \vert.$ 
\end{proof}

%2.6
\begin{theorem} Suppose in (2.3) that $\beta = 2$ and $\alpha = 1.$ Then $(-K_X)^3 \in \{4,6,8,10,12,16\}.$ If $2d := (-K_X)^3, $
then $X \subset \bP(\sO_{\bP_1}(2) \oplus \sO^d)$ if $d  \geq 3$ resp. $X$ is a two-sheeted cover of $ \bP(\sO_{\bP_1}(2) \oplus \sO^d)$
if $d = 2$ (A.2,no.4-9).
\end{theorem}

The cases $d = 2,3,4,8$ really exist; see Example 2.7 below. However it seems to be a possibly difficult problem to determine whether $d = 5,6$
are possible. 

\begin{proof} Since $\beta = 2, $ (2.3.3) gives $F \cdot l_{\psi} = 1.$ Hence all $l_{\psi}$ are irreducible and thus $\psi \vert D$ is a
$\bP_1-$bundle. Moreover $\psi \vert F$ is an isomorphism, so $F' = \psi (F) \simeq F.$ From (2.3.1) we obtain
$$ -K_{X'} - 2F' \equiv 0.$$ On the other hand, $2F+D \cdot l_{\psi} = 0,$
hence $2F'$ is Cartier in $X'.$ Since $X'$ is $\bQ-$Fano, this implies 
$$ 2F' \in \vert -K_{X'} \vert. \eqno (2.6.1) $$
By (2.3.4) we have $(-K_X)^3 = 2K_F^2 \leq 16,$ and $(-K_X)^3 \ne 14. $
First we consider the cases that $d = 1$ and $d = 2.$ 

\vskip .2cm \noindent
{\it  Case $K_F^2 = 1.$} Then $(-K_X)^3  = 2$ and $h^0(-K_X) = 4.$ Moreover $\mu: X' \to \bP_3$ 
has $\deg \mu = 2;$ on the other hand, $-K_{X'} \cdot B = 1$ 
by (A) in (2.4), so that $\deg (\mu \vert B) = 1$ and $\mu(B)$ is a line, contradicting the fact that  
$g(B) = 1$ by (***) in (2.4).
So $K_F^2 = 1 $ is ruled out. 

\vskip .2cm \noindent
{\it Case $K_F^2 = 2,$ i.e. $(-K_X)^3 = 4.$ } Since $h^0(-K_F) = 3$ and $h^0(-K_X) = 5$, the map 
$H^0(-K_X) \to H^0(-K_F)$ is surjective. Therefore
$\mu: X' \to \bP_4$ is a 2:1-covering over a singular quadric $Q \subset \bP_4$ and $F$ is mapped to a plane in $\bP_4$ 
by $\mu \circ \psi.$ We compute easily from the given data and from $D = -K_X -2F$ that 
$$ \phi_*(-K_X) = \sO(2) \oplus \sO^2.$$ 
The canonical map $X \to \bP(\sO(2) \oplus \sO^2)$ is a degree 2 covering. By contracting the exceptional divisor 
in $\bP(\sO(2) \oplus \sO^2)$, we obtain a birational map  to the singular quadric $Q$ from above and all maps commute. \\
Of course this construction can be reversed, but $X$ can also be constructed as follows. Consider the Fano manifold $V_2$, the
degree 2 covering of $\bP_3$ of index 2. Let $B \subset V_2$ be a complete intersection of two general elements in the half-anticanonical
system. Then $B$ is an elliptic curve. Let $\pi: \hat V \to V_2$ be the blowup of $B.$ Then $\hat V$ has a natural del Pezzo fibration 
$p$ provided by $\pi^*(\sO(1)) - E.$ Now we let $\tau: X \to \hat V$ be the degree 2 covering ramified along $p^*(\sO(2)).$  
Notice that one section of $-K_X$ is not the pullback under $\tau$!

\vskip .2cm \noindent 
{\it From now on suppose $d \geq 3.$ }
Write
$$ \phi_*(-K_X) = \bigoplus_{i=1}^r \sO(a_i) $$
with $a_1 \geq a_2 \ldots \geq a_r.$
Then 
$$ \phi_*(\sO_X(D)) = \bigoplus_{i} \sO(a_i-2).$$
Since $h^0(\sO_X(D)) = 1,$ we have $a_1 = 2$ and $a_i \leq 1$ for $i \geq 2.$ 
Since $H^1(-K_X) = 0,$ the only negative summands can be of type $\sO(-1), $ hence
$$ \phi_*(-K_X) = \sO(2) \oplus \sO(1)^u \oplus \sO^v \oplus \sO(-1)^w.$$
Therefore  $h^0(-K_X) = 3 + 2 u + v$ and since $h^0(-K_X) = d+3, $ we obtain
$$ d = 2u + v.$$
By $h^0(-K_F) = d+1$, 
the bundle $\phi_*(-K_X) $ has rank $r = d+1$. Now we claim that
$h^0(-K_X-F) = 2.$ In fact, every section of $-K_X-F$ has to vanish on $D$, simply because $-K_X-F$ is negative on the fibers
of $D \to B.$ Hence 
$$h^0(-K_X-F) = h^0(-K_X-F-D) = h^0(F) = 2.$$
Therefore $u = 0$ and $v = d.$ Now the rank condition implies $w = 0$ 
so that 
$$ \phi_*(-K_X) = \sO(2) \oplus \sO^d.$$
Since $-K_X$ is $\phi-$very ample, we obtain an embedding $X \subset \bP(\phi_*(-K_X)).$ 
\end{proof}

\begin{example} {\rm (1) Let $\sE = \sO(2) \oplus \sO^d$ and $\eta = \sO_{\bP (\sE)}(1)$. For $d = 2$ see the proof of (2.6). In the case $d = 3, 4$ simply take $X$ as a general member of $|3 \eta|$ or a complete intersection of two members of $|2\eta|$, respectively. Note that $\eta$ is spanned. 
\vskip .2cm \noindent
(2) Suppose now that $d = 8,$ i.e. $\phi: X \to \bP_1$ is a quadric bundle. Here $(-K_X)^3 = 16.$ If $l_{\phi}$ denotes a line in a fiber of
$\phi,$ then $-K_X \cdot l_{\phi} = 2.$ Since moreover $-K_X \cdot l_{\psi} = 0$ and since every curve in $X$ is numerically a linear
combination of $l_{\phi}$ and $l_{\psi} $ with integer coefficients, $-K_X$ is  divisible by $2:$ there is a line bundle $L$ such that
$$ -K_X = 2L.$$ Observe $L^3 = 2$ and $h^0(L) = 4.$ Moreover $L = \psi^*(L') $ and $-K_{X'} = 2L'.$
By Fujita's results \cite{Fu90}, \cite{AG5}, $L'$ is spanned and provides a double cover $X' \to \bP_3$, ramified along a quartic.  
Conversely, we argue as follows. Two general quadrics $q_1, q_2$ in $\PN_3$ meet in a smooth elliptic curve $B$. Let $X'$ be the 2:1-covering of $\PN_3$, ramified along a general quartic in $H^0(\bP_3,\sO(4) \otimes \sI_B^2)$. The blowup of $X'$ in the reduced preimage of $B$ is a smooth threefold $X$. The anticanonical divisor is the pullback of $\sO_{\PN_3}(2)$, and $X$ admits a pencil defined by $q_1, q_2$. Here $D = -K_X - 2F$.

}
\end{example}  

%2.8
\begin{theorem} In the setting of (2.3) suppose that $\beta = 1$. Then $X$ is one of the following and all these cases really exist.

\vspace{0.2cm}

\begin{center}
\begin{tabular}{c||c||c||p{5.5cm}|c||c|}
  $No.$ & $-K_X^3$ & $K_F^2$ & $X'$ & $d/g_B$ & $(\alpha, \beta)$\\\hline
  $1$ & $2$ & $2$ & 2:1-cov. of $\PN_3$, ram. sextic & $1/0$ & $(1,1)$\\
  $2$ & $4$ & $4$ & $X'_4 \subset \PN_4$ or \newline 2:1-cov. of $Q_3$, ram. quartic & $2/0$ & $(1,1)$\\ 
  $3$ & $8$ & $8$ & 2:1-cov. of Ver. cone & $4/0$ & $(1,1)$\\ 
  $4$ & $2$ & $4$ & 2:1-cov. of $\PN_3$, ram. sextic & $4/1$ & $(2,1)$\\\hline
\end{tabular}
\end{center}
These are the cases (A.2),no.10-13.  
\end{theorem} 

\begin{proof} Suppose that $\beta = 1.$ Then $D = -\alpha K_X - F$ and $\alpha(-K_X)^3 = K_F^2.$ 
This leads to the following cases.
\begin{itemize}
\item $\alpha = 1 $ and $(-K_X)^3 = K_F^2 = 2,4,6,8;$
\item $\alpha = 2$ and $(-K_X)^3 = 2,4; K_F^2 = 4,8,$ respectively;
\item $\alpha = 3 $ and $(-K_X)^3 = 2, K_F^2 = 6;$
\item $\alpha = 4 $ and $(-K_X)^3 = 2, K_F^2 = 8.$
\end{itemize}
Since $F \cdot l_{\psi} = 2,$ we have
$$ -K_{X'} \cdot B = - {{1} \over {2}} K_X \cdot D \cdot F = {{\alpha} \over {2}} K_F^2 = \alpha {{\alpha} \over {2}} (-K_X)^3. \eqno (*)$$
Now we proceed as in the proof of (2.4) and obtain the following formula for the genus $g$ of the curve $B.$
$$ g = (-K_X)^3({{\alpha^3} \over {3}} - {{\alpha^2} \over {4}} - {{\alpha} \over {12}}) -2\alpha+2. \eqno (**)$$ 
The second method to compute $g$ unfortunately fails, since $D_F \to B$ is now a degree 2 covering. 
\vskip .2cm \noindent 
(1) Suppose $\alpha = 1$. Then $g = g(B) = 0$ and $D = -K_X-F.$ \\
(1.a) If $(-K_X)^3 = K_F^2 = 2,$ then $\phi_*(-K_X) = \sO(1) \oplus \sO^2 $. This follows easily from the intersection numbers
and $h^0(-K_X-F) = 1.$ The map $X' \to \bP_3$ is a two-sheeted cover and so does $X \to \bP(\sO(1) \oplus \sO^2)$ which is just 
$\bP_3$ blown up in a line. Conversely, let $R \in \vert 4 \zeta + \pi^*(\sO(2))\vert$ be smooth and $g: X \to \bP(\sO(1) \oplus \sO^2)$ be
the two-sheeted cover branched along $R.$ Let $D = g^*(\bP(\sO^2)).$ Then there is a contraction $X \to X'$ along a projection of $D$,
the variety $X'$ admits a covering $X' \to \bP_3$ and this commutes with the blowdown $\bP(\sO(1) \oplus \sO^2) \to \bP_3.$ 
\vskip .2cm \noindent 
(1.b) Now suppose $(-K_X)^3 = K_F^2 = 4.$ Here $X' \subset \bP_4$ of degree 4 or $X' \to Q$ is a degree 2 covering over a quadric.  
The bundle $\phi_*(-K_X) $ is or rank 5 with 5 sections and therefore is of the form $\sO(1) \oplus \sO^3 \oplus \sO(-1).$ Moreover
the canonical map $X \to \bP(\sO(1) \oplus \sO^3 \oplus \sO(-1))$ is either an embedding or a two-sheeted cover over its image. \\
In case of an embedding, we can describe $X$  also as follows. Notice that by $(*), $ the smooth rational curve $B$ is a conic, contained in some $\PN_2 \subset \PN_4$, defined by two hyperplanes $H_1$ and $H_2$. For $X'$ choose a general quartic in $|\sO_{\PN_4}(4) \otimes \I^2_{B/\PN_4}|$. After blowing up $B$ we obtain $X$ with $-K_X \simeq \psi^*\sO_{\PN_4}(1)$ big and nef and $(-K_X)^3 = 4$. The pencil $\langle H_1, H_2\rangle$ defines a del Pezzo fibration $X \lra \PN_1$ with general $F$ having $K_F^2 = 4$. The case $\mu: X' \to Q_3$ is analogous. 
\vskip .2cm \noindent 
(1.c) If $(-K_X)^3 = K_F^2 = 6,$ then $X' \subset \bP_5$ is a complete intersection of a cubic $Z$ and a quadric $Q$ and $B$ is a twisted cubic. In particular $\mu $ is an embedding.  
Let $\pi: \hat P \to \bP_5$ denote the blowup of $\bP_5$ along $B.$ Then $X \subset \hat P$ is the complete intersection of the strict transforms 
$\hat Z$ and $\hat Q.$ Using adjunction and $-K_X = \pi^*(\sO(1)),$ we have
$$\hat Z \in \vert \pi^*(\sO(3)) - 2E \vert$$ 
and 
$$ \hat Q \in \vert \pi^*(\sO(2)) - E \vert.$$ 
   
Now $B$ is contained in a 3-dimensional linear subspace $Y \subset \bP_5$.
We claim that $H^0(Y,\sI_B^2(3)) = 0.$ It is classical that we find a smooth quadric $Q_2 \subset  Y$ containing
$B.$ Now $B$ is a curve of type (say) $(2,1)$ (possibly after permuting the factors), therefore
$$H^0(Q_2,\sI_B^2(3)) = H^0(Q_2,\sO(-1,1)) = 0 $$
and thus the vanishing is achieved. \\
We conclude that $Y \subset Z.$ Any component of the singularity set of $Z$ is either contained in $B$ or disjoint from $Q \cap Z$. 
Consider the exact sequence of conormal sheaves
$$ 0 \to N^*_{Z/\bP_5} \vert Y \to N^*_{Y/\bP_5} \to N^*_{Y/Z} \to 0 $$
which reads
$$ 0 \to \sO(-3) \to \sO(-1) \oplus \sO(-1) \to N^*_{Y/Z} \to 0.$$
The map $\gamma: N^*_{Z/\bP_5} \vert Y \to N^*_{Y/\bP_5}$ is given by sections $s_1,s_2 \in H^0(\sO_Y(2)).$ Since we already
know that the zero locus of $\gamma $ is contained in $B$ and since $B$ has degree 3, we obtain a contradiction and 
rule out this case. 

\vskip .2cm \noindent
(1.d) If $(-K_X)^3 = K_F^2 = 8,$ then either $X' \subset \bP_6$ has degree 8 and $B$ has degree 4 or 
$X'$ is a degree 2 covering over the cone over the Veronese surface. We will exclude the first case  and show that the latter leads to an example.
In the first case $X'$ is the complete intersection of three quadrics and $B$ is a smooth rational curve of degree 4 in $X'.$ Blowing up $\bP_6$ along $B$ and computing as in (1.c), two of the quadrics are in $\vert \sO(2)-E \vert$ while the third one, 
say $Q$, is in $\vert \sO(2)-2E \vert.$ The quadric $Q$ is singular along $B$. Since $Q_{sing}$ is a linear subspace of $\PN_6$, while $B$ neither is a line nor a plane curve, we conclude $Q_{sing} \simeq \PN_3$. Hence $B$ is the complete intersection of two quadrics in $\PN_3$. But this is not a rational curve.\\
It remains to treat the case that $X'$ is a degree two covering of the Veronese cone in $\PN_6$. Choose two general hyperplanes $H_1$ and $H_2$. The intersection with the Veronese cone gives $C_4$. Choose a general cubic in $\sO_{\PN_6}(3) \otimes \I_{B/\PN_6}^2$ and let $X'$ be the double covering of the Veronese cone, ramified along this cubic. After blowing up $B$ we get $X$ and a pencil $X \lra \PN_1$ with general fiber $\PN_1 \times \PN_1$, $2:1$ over $\PN_2$ ramified along a conic.

\vskip .2cm \noindent 
(2.a) If $\alpha = 2$ and $(-K_X)^3 = 2,$ we obtain $g(B) = 1$ and $-K_{X'} \cdot B = 4.$ 
Since $B$ is the singular locus of $X',$ the curve $B$ must be inside the ramification of $\mu$ and thus $\deg \mu \vert B = 1.$ Thus $B' = \mu(B)$
is an elliptic curve in $\bP_3$ of degree 4. Therefore $B$ is the complete intersection of two quadrics. Also observe that $F$ is mapped
onto a quadric $F' \subset \bP_3.$ Therefore the image $Z$ of $X \to \bP(\phi_*(-K_X))$ is a quadric bundle over $\bP_1$ which
admits a birational map to $\bP_3.$ This must be the blowup of $B'.$  \\
In order to construct $X$, we start with $B' \subset \bP_3$ and let $\pi: \hat P \to \bP_3 $ be the blowup of $B' \subset \bP_3.$ Then $\hat P$
is a quadric bundle over $\bP_1.$ 
We want to construct
$X$ as 2:1 covering over $\hat P.$ Denote by $p: \hat P \to \bP_1$ the quadric projection and observe that $-K_{\hat P} - \pi^*(\sO(2))$
is trivial on the fibers $\hat F$ of $p$ so that  $-K_{\hat P} - \pi^*(\sO(2)) = p^*(\sO(1))$. In particular, the bundle
$-2K_{\hat P} - \pi^*(\sO(2))$ is spanned and we let $R$ be a general member of the associated linear system. Since $R$ is divisible 
by 2, we can take the two-sheeted cover $g: X \to \hat P$ ramified along $R.$ 
Then
$$ K_X = g^*(K_{\hat P} - R/2) = g^*\pi^*(\sO(-1)) $$
so that $-K_X$ is big and nef but not ample. The induced map $X \to \bP_1$ has general fiber $F$ with $K_F^2 = 4; $ moreover 
$(-K_X)^3 = 2.$ It is also easy to check that $\alpha = 2$ and $\beta = 1.$

\vskip .2cm \noindent 
(2.b) If $\alpha = 2$ and $(-K_X)^3 = 4,$ then $g(B) = 4$ and $-K_{X'} \cdot B = 8.$ The map $\mu: X' \to \bP_4$ is either an embedding or
a 2:1 covering over a quadric $Q_3 \subset  \bP_4.$ \\
First we consider the case where $\mu$ is an embedding.
In that case $B \subset  \bP_4$ is a smooth curve of degree $8$. From
$F = -2K_X-D$ we infer that $B$ is contained in a pencil of quadrics
in $\PN_4$. From $F\cdot l_{\psi} = 2$ we infer that all the quadrics must be singular along $B$. But $B$ neither is a line nor a plane curve. Hence $\mu$ is not an embedding. \\
If $\deg \mu = 2,$ then $W \subset  \bP_4$ is a quadric. Let $B' = \mu(B).$ We claim that $\mu \vert B$ is $1:1$ to $B'.$ This is clear if 
$W$ has at most an isolated singularity. 
If $W$ is singular along a curve, necessarily along a line, and if $\mu \vert B$ is not $1:1,$ then $B'$ must be the singular locus of $W$,
i.e. $B'$ is a line and $B \to B'$ is $2:1.$ By Riemann-Hurwitz, $K_B = \mu^*(\sO(-2)) + R'$ so that $\deg R' = 10.$ On the other hand,
$R' = B \cdot R$ with $R \subset  X'$ the branch locus of $\mu: X' \to W.$ Since $R = \mu^*(\sO(1)),$ we obtain a contradiction to
$-K_{X'} \cdot B = 8.$ \\
So $\mu \vert B$ is $1:1$ and $B \subset R.$ Since $D = -2K_X - F$,
the images of $F$ in $W$ are quadric surfaces in $\bP_4.$ Hence $B'$ is contained in the
intersection of three quadrics in $\bP_4$ and since $\deg B' = 8,$ we have equality. But a smooth complete intersection of three quadrics
has genus 5. \\

\vskip .2cm \noindent
If $\alpha = 3,$ then $(-K_X)^3 = 2,$ moreover $K_F^2 = 6$ and $g(B) = 9$ and $-K_{X'} \cdot B = 9.$ Let $B' = \mu(B);$ notice that $\mu \vert B$
is an isomorphism. Then
$B'$ is a smooth curve in $\bP_3$ of degree 9 and genus 9. Since $F'' = \mu(F') $ has degree 3, $B'$ is in the intersection of
two hypersurfaces of degree 3, hence equals this complete intersection. But then we obtain $g = 10,$ contradiction.

\vskip .2cm \noindent
If $\alpha = 4,$  the argument is completely the same: here $g(B) = 28 $ and has degree 16. Thus $B' \simeq B$ is the complete intersection
of two quartics. But then the genus must be 31, again a contradiction. 
\end{proof}

%2.9
\begin{theorem} Suppose in (2.3) that $K_F^2 = 8.$ If $\alpha$ is not an integer then either
$X$ is a general divisor of the form $2 \zeta $ in $\bP(\sO(2) \oplus \sO^3)$ over $\bP_1$ (A.2,no.14) \\
 or $X$ is a general divisor of the form $2 \zeta + 2F$ in $\bP(\sO(1) \oplus \sO^2 \oplus \sO(-1))$ over $\bP_1$ (A.2,no.15).
\end{theorem}

\begin{proof} We already know that $\alpha \in {{1} \over {2}} \bN.$ Suppose that $\alpha$ is not an integer. Then we find an odd integer $\tilde{\alpha}$
such that $\alpha = {{\tilde{\alpha}} \over {2}}. $ 
\vskip .2cm \noindent
(1) First suppose that $\beta \in \bZ.$ Then by (2.3.3) $\beta = 1, 2.$
\vskip .2cm \noindent
(1.a) $\beta = 2.$ Then ${{-K_X} \over {2}}$ is Cartier and using (2.3.4): $ \alpha(-K_X)^3 = 2K_F^2 = 16,$ which implies $\tilde{\alpha} = 1, \alpha = {{1} \over {2}} $ and
$(-K_X)^3 = 32.$ Let $L = {{-K_X} \over {2}} = D + 2F$ so that $X \subset \bP(\phi_*(L)).$ In order to determine $\phi_*(L)$ 
we compute by Riemann-Roch and Kawamata-Viehweg
$$ h^0(L) = 6.$$
Moreover $h^0(L_F) = 4$ so that $\phi_*(L)$ has rank 4. Using $h^1(L) = 0$ and $h^0(L-2F) = 1,$ we find
$$ \phi_*(L) = \sO(2) \oplus \sO^3 \eqno (A) $$
or
$$ \phi_*(L) = \sO(2) \oplus \sO(1) \oplus \sO \oplus \sO(-1). \eqno (B) $$
Write $X = r \zeta + s F$. A priori we know $r = 2$ and $L = \zeta \vert X.$\\
In case (A)  we get  $s = 0$
in order to have $-K_X = 2L.$ Let 
$$ \sigma: \bP(\sO(2) \oplus \sO^3) \to W $$
be the blowdown of $\bP(\sO^3)$ to a surface $S \simeq \bP_2$ or equivalently, the morphism defined by $\vert \zeta \vert.$ 
Then  $X = \sigma^*(X')$ and since $X$ is irreducible, $X'$ does not contain $S$ so that $X'$ meets $S$ in a curve $B$, so that $X \to X'$ contracts some divisor in $X.$ Thus $-K_X$ is nef (and big) but not ample. \\
In case (B), $-K_{\bP} = 4 \zeta $ is not nef exactly on $C_0 = \bP(\sO(-1)).$ Hence $-K_X = 2 \zeta \vert X$ is nef
exactly if $C_0$ is not contained in $X.$ However the base locus of each system $\vert m \zeta \vert$ contains $C_0$, so that
always $C_0 \subset X.$ Hence case (B) is ruled out.
\vskip .2cm \noindent
(1.b) $\beta = 1.$ Hence $\alpha(-K_X)^3 = K_F^2 = 8$ so that again $\tilde{\alpha} = 1$ and moreover $(-K_X)^3 = 16.$ Defining $L$ as before
we this time have $L = D+F.$ Since $h^0(L) = 4,$ it follows
$$ \phi_*(L) = \sO(1) \oplus \sO^2 \oplus \sO(-1).$$
Here the adjunction formula and $-K_X = 2L$ force $X = 2\zeta + 2F$. In order $-K_X$ is nef, we need $\bP(\sO(-1))$ not
to be contained in $X$ which of course is true for general choice of $X$. 
\vskip .2cm \noindent
(2) $\beta \not \in \bZ$. Of course $2\beta \in \bZ$   
and  $\beta = {{\tilde{\beta}}\over {2}}$ with $\tilde{\beta}$ odd. 
Since $D \cdot l_{\psi} = -2,$ we obtain $\tilde{\beta} (F \cdot l_{\psi}) = 4,$ hence $\tilde{\beta} = 1$.
By (2.3.4) (which is purely numerical and hence also true here) we have
$$ {{\tilde{\alpha}} \over {2}} (-K_X)^3  =  \alpha(-K_X)^3 = \beta K_F^2 = 8 {{\tilde{\beta}} \over {2}} = 4 ,$$
so $\tilde{\alpha} = 1$ and $(-K_X)^3 = 8.$ \\
From Riemann-Roch we compute
$$\chi(X,\sO_X(D)) = 1.$$
Now consider the rank 4-bundle $\phi_*(\sO_X(D)).$ Since $h^0(\sO_X(D)) = 1$ and, by the Leray spectral sequence, $H^q(X, \sO_X(D)) = 0$ for $q=2,3$, we deduce $h^1(\sO_X(D)) = 0.$ Hence
$$ \phi_*(\sO_X(D)) = \sO \oplus \sO(-1)^{\oplus 3}.$$
Writing again $X = 2 \zeta + s F$, we obtain $s = 4$ and $-K_X = 2 \zeta \vert X + F.$ 
Now $\bP(\sO \oplus \sO(-1)^3) $ is nothing than $\bP_3$ blown up in a line $l$. Let $W = \pi(X).$  Since
$$X = 2 \zeta + 4F = \pi^*(\sO(4)) - 2E $$
with $E$ the exceptional divisor of $\pi$, we have $l \subset W.$ Therefore $X$ contains curves contracted by $\pi$ and
thus $-K_X = 2 \zeta + F$ cannot be nef. The fact that $X$ contains curves contracted by $\pi$ is also clear from the
fact that $X$ is an ample divisor. 
\end{proof}
 
\

%%%%%%%%%%%%%%%%%%%%%

\section{Conic bundles} \label{conic}
\setcounter{lemma}{0}

%4.1
\begin{setup} {\rm In this section $\phi: X \to Y = \bP_2$ denotes a conic bundle with $\rho(X) = 2$. 
As always we assume $-K_X$ big and nef but not ample. We do {\em
  not} assume that the anticanonical morphism is divisorial until
after Proposition~\ref{Enotgen}. 

The discriminant locus is denoted by $\Delta,$ if
$\Delta = \emptyset, $ then $X$ is a $\bP_1-$bundle. Set 
$$d = \deg \Delta.$$
We introduce the rank 3-bundle
$$ \sE = \phi_*(-K_X)$$ 
which is relatively spanned, i.e., $\phi^*\sE \lra \sE$ is
surjective. We obtain an embedding 
$$ X \subset  \bP(\sE)$$
such that $-K_X = \zeta \vert X.$ 
The divisor $X \subset \bP(\sE)$ is of the form
  \begin{equation}\label{XinS2}
     [X] = 2 \zeta + \pi^*(\sO(\lambda)) 
  \end{equation}
with some integer $\lambda.$ Then the adjunction formula yields
$$ \lambda = 3 - c_1.$$  
Here we use the shorthand $c_i = c_i(\sE).$ 
Since
$$ H^q(\bP(\sE),- \zeta - \pi^*(\sO(\lambda))) = 0$$
for $q = 0,1,$ every section in $H^0(-K_X)$ uniquely lifts to a
section of $\zeta.$ We first consider the case where $\sE$ is not spanned:}
\end{setup}

\begin{proposition} \label{Enotgen}
  Assume $\sE = \phi_*(-K_X)$ is not globally generated. Then $X'$ is
  hyperelliptic and $X$ is a double covering of a Fano $\PN_1$--bundle over $\PN_2$
    \[\sigma: X \lra \PN(\sF),\]
  such that $-K_X = \sigma^*\sO_{\PN(\sF)}(1)$. The Fano bundle $\PN(\sF)$ is one of the following:
  \begin{enumerate}
     \item $\PN_3$ blown up in a point, i.e., $\sF =
       \sO_{\PN_2} \oplus \sO_{\PN_2}(1)$. Here $\sE = \sO_{\PN_2}(-2)
       \oplus \sF$ and $(-K_X)^3 = 2$ (A.3, no.9).
     \item $Q_3$ blown up in a line, i.e., $\sF$ sits in the exact sequence
        \[0 \lra \sO_{\PN_2}(1) \lra \sF \lra {\mathcal I}_p(1) \lra 0,\] 
       where ${\mathcal I}_p$ denotes the ideal sheaf of a point $p \in
       \PN_2$. Here $\sE = \sO_{\PN_2}(-1) \oplus \sF$ and $(-K_X)^3 =
       4$ (A.3, no.10).
     \item $\PN_3$ blown up in a twisted cubic, i.e., $\sF$ sits in
          \[0 \lra \sO_{\PN_2}(-1)^{\oplus 2} \lra \sO_{\PN_2}^{\oplus
            4} \lra \sF \lra 0.\] 
        Here $\sE = \sO_{\PN_2}(-1) \oplus \sF$ and $(-K_X)^3 = 2$ (A.3,
        no.11).
     \item The Veronese cone blown up in its vertex, i.e., $\sF =
       \sO_{\PN_2} \oplus \sO_{\PN_2}(2)$. Here $\sE =
       \sO_{\PN_2}(-1) \oplus \sF$ and $(-K_X)^3 = 8$ (A.3, no.12).
   \end{enumerate}
  In all these cases $\psi$ is divisorial. 
\end{proposition}

We have the following diagram in Proposition~\ref{Enotgen}:
 \[\xymatrix{X \ar[d]^{\sigma} \ar[rrr]^{\psi} & & & X' \ar[d]^{\mu} \\
       \PN(\sF) \ar[rrr]^{|\sO_{\PN(\sF)}(1)|} & & & W }\]
   
\vspace{0.2cm}

\begin{proof} Think of $C_p$, the scheme--theoretic fiber of $X$ over a point $p \in \PN_2$, as a conic in
  $\PN(\sE(p)) \simeq \PN_2$. The vertical maps in 
    \[\xymatrix{H^0(\PN_2, \sE) \ar@{=}[d]\ar[r] & \sE(p) \ar@{=}[d] \\
                H^0(X, -K_X) \ar[r] & H^0(C_p, -K_X|_{C_p})}\] 
 are isomorphisms by which we identify these spaces. In this way
 we may think of $|H^0(C_p,
  -K_X|_{C_p})|$ as lines in $\PN(\sE(p))$. Since
  $|-K_X|$ is spanned, the image of $H^0(X, -K_X) \lra H^0(C_p,
-K_X|_{C_p})$ is at least two dimensional, describing all lines, or the family of lines through a given point outside $C_p$.

Now assume that $\sE$ is not generated by global sections at the point
$b \in \PN_2$. Then 
 \[H^0(X, -K_X) \lra H^0(C_b, -K_X|_{C_b}) \eqno (*)\]
is not surjective and the image is two
dimensional. On $\PN(\sE(b))$, the map $\psi$ is projection
from a point. Except for the case $C_b$ a double line, $X'$ is clearly
hyperelliptic, and if $R$ denotes the
strict transform in $X$ of the ramification divisor of $X' \lra
W$, then $R\cdot C_p = 2$.

\vspace{0.2cm}

We now treat the case that $C_b$ is not a double line for some base
point $b$ of $\sE$. Denote by $C_p$ a {\em general} conic and by $C'_p$ its strict transform under the
hyperelliptic involution. At $b$ we have $C_b = C'_b$. Hence, if $H =
\phi^*\sO_{\PN_2}(1)$, then $0 = H\cdot C_p = H\cdot C_b = H\cdot C'_b =
H\cdot C'_p$, implying that $C'_p$ will be a fiber of $\phi$. But
$R\cdot C_p > 0$, hence $C_p \cap C'_p \not= \emptyset$ so that $C_p
= C'_p$. Then $\psi$ is 2:1 on the {\em generic}
conic and therefore given by projection in $\PN(\sE(p))$ from some point. 

We have proved the following in terms of vector bundles. If the canonical map of
global sections of
\[\sO_{\PN_2} \otimes H^0(\PN_2, \sE) \lra \sE\]
is not surjective at some point $b \in \PN_2$ (corresponding to a
reduced conic $C_b$), it has rank $2$ at the general point, hence its
rank is $2$ everywhere by (*). Thus we have an exact
sequence of vector bundles,
  \[0 \lra \sF \lra \E \lra L \lra 0,\]
with $L$ some line bundle. The rational projection
  \[\xymatrix{\PN(\E) \ar@{-->}[r] & \PN(\sF)}\]
is defined outside the section $S$ of $\PN(\sE)$ corresponding to $\sE
\lra L$. Since $S \cap X = \emptyset$, the projection is holomorphic on $X$, defining a 2:1--covering
  \[\sigma: X \lra \PN(\sF), \quad \sigma^*\sO_{\PN(\sF)}(1) = -K_X,\]
factorizing $\psi$. By the Leray spectral sequence, $H^1(\PN(\sF), \sO_{\PN(\sF)}(1)) = H^1(\PN_2, \sF) =
0$. Then $L$ cannot be spanned, i.e.,
 \[L = \sO_{\PN_2}(k), \quad k < 0.\]
Since $S \cap X = \emptyset$, we can now compute $L$ in terms of
$\sF$. Using (\ref{XinS2}) and $\zeta|_S = L$, we find
 \[L^{\otimes 2} \simeq \det \sE \otimes K_{\PN_2},\]
or, using $\det \sE \simeq \det \sF \otimes L$,
  \[-K_{\PN_2} \simeq \det \sF \otimes L^*.\]
As a first consequence, we note that $S^2\sF \otimes \det \sF^* \otimes (-K_{\PN_2}) \simeq S^2\sF \otimes
L^*$ is ample, which implies that $\PN(\sF)$ is Fano. Moreover we see that there are
only two possibilities for $(\det \sF, L)$:
 \[(\sO_{\PN_2}(2), \sO_{\PN_2}(-1)), \mbox{ or } (\sO_{\PN_2}(1),
 \sO_{\PN_2}(-2)).\]
The case $\det \sF \simeq \sO_{\PN_2}$ is impossible since $\sF$ cannot be trivial.

Using the list of ``Fano bundles'' over $\PN_2$ (\cite{SW90}) we
obtain $X$ as stated in the
Proposition. By Griffiths' vanishing theorem,
$H^1(\PN_2, \sF \otimes L^*) = 0$, which implies $\sE \simeq \sF \oplus L$. 

\vspace{0.2cm}

It remains to consider the case that for any $b$, where $\sE$ is not
spanned, the fiber $C_b$ is a double line. In the hyperelliptic case we may assume $R\cdot C_p = 0$,
for otherwise we can conclude exactly as above. In both cases, whether
$X'$ is hyperelliptic or not, $\psi$ is not an isomorphism at any
point of $C_b$, hence $\psi$ is divisorial, sending $D$ to the line
$\mu(\psi((C_b)_{red})) \subset \PN_{g+1}$. Note that $R\cdot C_p = 0$ in the
hyperelliptic case, which implies that all ramification comes from $\PN_2$.

A line is cut out by hyperplanes. Hence $-K_X - \alpha D$ will be generically spanned for some $\alpha
\in \KN$, which implies $(-K_X - \alpha D)\cdot C_p \ge 0$. From $-K_X\cdot
C_p = 2$ and the fact that $D\cdot C_p$ must be
even, since $\Delta \not= \emptyset$, we get $-K_X\cdot C_p = D\cdot C_p$ and
  \[-K_X - D = \phi^*\sO_{\PN_2}(k)\]
for some $k \in \KN$. Let $l_{\psi}$ denote the general positive
  dimensional fiber of $\psi$. Intersecting with $l_{\psi}$ using
  $D\cdot l_{\psi} = -2$ shows $k = 1$ or $2$.

We prove $H^1(X, \sO_X(D)) = 0$. By Serre duality this is equivalent
to $h^2(\sO_X(K_X-D)) = 0$. From
  \[0 \lra \sO_X(K_X-D)  \lra \sO_X(K_X) \lra \sO_D(K_X) \lra 0\]
and $h^2(\sO_X(K_X)) = 0$ we see that it suffices to prove
$h^1(\sO_D(K_X)) = 0$. Since $\psi(D)$ is a line, $h^1(\sO_D(K_X)) =
h^1(\PN_1, \sO_{\PN_1}(-1)) = 0$.

But $H^1(X, \sO_X(D)) = H^1(X, -K_X - \phi^*\sO_{\PN_2}(k)) = 0$, for
$k = 1, 2$. This says that the restriction map
 \[H^0(\PN_2, \sE) \lra H^0(C, \sE|_C)\]
is surjective for any line or conic $C$ in $\PN_2$. Since $\sE$ is globally generated away from
points, $\sE|_C$ is globally generated. This contradicts
$Bs(\sE) \not= \emptyset$. The Proposition is proved.
\end{proof}

\

From now on we assume that $\sE$ is globally generated. And we assume
again that $\psi$ is divisorial. We write for the
irreducible exceptional divisor $D:$
$$ D = \alpha (-K_X) + \phi^*(\sO(\beta))$$
with rational numbers $\alpha$ and $\beta$. Actually $\alpha $ and $\beta$ are integers unless $\Delta = \emptyset$
(intersect with an irreducible component of a reducible fiber).  
Thus $\vert \zeta \vert $ defines via Stein factorisation a  map $\tilde \psi: \bP(\sE) \to \bP'  $ extending $\psi$  and in total
a map $\tau \circ \tilde \psi: \bP(\sE) \to \tilde W \subset \bP_{g+1}.$ 

%4.2
\begin{theorem} $B = \psi(D)$ is a smooth curve unless $X = \bP(\sO_{\bP_2} \oplus \sO_{\bP_2}(3))$.
\end{theorem}

\begin{proof} We must first show that $\dim \psi(D) = 1.$ So suppose $\dim \psi(D) = 0.$ 
Then $K_X \vert D \equiv 0,$ hence $D$ must be a multi-section of $\phi.$ 
Write $l = l_{\phi}$ for a smooth fiber of $\phi$ over $p \in \bP_2.$ Since $\sE$ is spanned, 
$$ H^0(\sE) \to  H^0(\sE \vert \{p\})$$ 
is onto, hence so does
$$ H^0(-K_X) \to H^0(-K_X \vert l).$$
Hence $\mu \circ \psi \vert l_{\phi}$ is an embedding, $D$ is a section of $\phi$ and $\phi$ is a $\bP_1-$bundle. Write $X = \bP(\sF).$ 
Since $\phi$ has a section, $\sF$ splits and after normalising we can write
$$ \sF = \sO \oplus \sO(a).$$
Since $-K_X$ is nef but not ample, $S^2(\sF) \otimes \det \sF^* \otimes -K_{\bP_2}$ is nef but not ample. 
This gives $a = 3.$ \\
Finally the smoothness of $B$ is (\ref{Wilson}). 
\end{proof}

First we classify the $\bP_1-$bundles, i.e. the case $\Delta = \emptyset.$  
%4.3
\begin{theorem} If $\Delta = \emptyset,$ then $X = \bP(\sF) $ is one of the
following.
\begin{enumerate}
\item $\sF = \sO_{\bP_2} \oplus \sO_{\bP_2}(3),$ (A.3,no.1);
\item $\sF$ is given by an extension 
$$ 0 \to \sO \to \sF \to \sI_p(-1) \to 0 $$
with a point $p \in \bP_2,$ (A.3,no.2);
\item $\sF$ is given by an extension 
$$ 0 \to \sO(-1) \to \sF \to \I_Z \to 0$$
where $Z$ has length $4,$ (A.3,no.3); 
\item $\sF$ is given by an extension
$$ 0 \to \sO(-2) \to \sF \to \sI_Z(1) \to 0,$$
where $Z$ has length $6$, (A.3,no.4).   
\end{enumerate}
\end{theorem} 

\begin{proof} Write
$$ X = \bP(\sF) $$
with a rank 2-bundle $\sF$ on $\bP_2$ normalised such that $c_1(\sF) = 0,-1.$ Set $\eta = \sO_{\bP(\sF)}(1) $ 
and $D = a \eta + \phi^*(\sO(b))$ with integers $a,b.$ By (3.3) we may assume that $\dim \psi(D) = 1.$
\\
First suppose $c_1(\sF) = 0.$ The nefness of $-K_X$ translates into the nefness of $\sF({{3}\over {2}}).$ However this bundle is
not ample. 
Now consider a curve $l_{\psi}$
contracted by $\psi,$ i.e. $K_X \cdot l_{\psi} = 0,$ i.e. 
$$(\eta+\phi^*(\sO({{3} \over {2}}))) \cdot l_{\psi} = 0.$$ 
Since $l_{\psi} $ must be a section of $\phi$ over its image, $\phi(l_{\psi}) $ has degree at most 2 and so
we have $\phi^*(\sO(1)) \cdot l_{\psi} = 2.$ So $\sF({{3} \over {2}}) $ is ample on all lines but not ample on a 1-dimensional family of conics.
This is impossible: if $V$ is a vector bundle over $\bP_2$ which is ample on all lines, then $V$ cannot be non-ample on the general member of a 1-dimensional family
of conics. 
Indeed, if $V$ is non-ample on such a family of conics, then this family must contain a splitting member (since the splitting conics form an
ample divisor in the $\bP_5$ of all conics). But $V$ is ample on this splitting conic, hence ample on the general member of the conic, since ampleness
is an open condition. \\
Thus we cannot have $c_1(\sF) = 0.$ 
\vskip .2cm \noindent
 Now suppose $c_1(\sF) = -1.$ Here $\sF(2)$ is nef but not ample. 
The equation $K_X^2 \cdot D = 0$ leads to
$$ c_2(\sF) = 3{{b} \over {a}} + 1.$$ 
Since $(\eta + \phi^*(\sO(2))) \cdot l_{\psi} = 0,$ 
we obtain $b = 2a-1$ or $b = 2a-2.$ Since $c_2(\sF) $ is an integer, this leads to the cases $(a,b) = (1,1), (3,5)$ resp. $(1,0), (2,2), (3,4), (6,10)$ and $c_2 = 4,6$ resp. $c_2 = 1,4,5,6.$ \\
Riemann-Roch gives
$$ h^0(\sF(2)) = \chi(\sF(2)) = 9 - c_2(\sF) \geq 3.$$
Now choose the minimal integer $\lambda$ such that 
$$ H^0(\sF(\lambda)) \ne 0.$$ 
The existence of $D$ with $a,b > 0$ yields $\lambda \geq 0;$ actually $\lambda \geq {{b} \over {a}},$ since $D$ sits on the boundary of the 
effective cone. \\
If $\lambda = 0,$ then $b = 0$, $a = 1$ and $c_2 = 1.$ Since $\sF$ has a section without zeroes in codimenion 1, we obtain a sequence
$$ 0 \to \sO \to \sF \to \sI_p(-1) \to 0  \eqno (*)$$
with $p$ some point in $\bP_2.$ 
Conversely, we construct $\sF$ by the Serre construction as an extension (*). Then $\sF(2) $ is spanned and, setting $X = \bP(\sF)$,
$-K_X$ is big and nef but not ample. The morphism associated to $-K_X$ contracts a divisor, namely the union of all negative sections sitting
over the lines through $p.$ \\
If $\lambda = 1,$ then we cannot have $b = 0$ (since then $a = 1$, contradicting the minimality of $\lambda$), hence $(a,b) = (1,1),(2,2).$ 
The second case is impossible since then $D$ would be a multiple of some effective divisor. Hence $c_2(\sF) = 4$ and we end up with a sequence
$$ 0 \to \sO \to \sF(1) \to  \sI_Z(1) \to 0 $$
where $Z$ has length $ c_2(\sF(1)) = c_2(\sF) = 4. $ For the existence we argue as before. 
\\ If $\lambda = 2,$ then we have a sequence
$$ 0 \to \sO \to \sF(2) \to \sI_Z(3) \to 0 $$
where $Z$ has length $ c_2(\sF(2)) = c_2 + 2.$ Since $h^0(\sF(1)) = 0,$ the bundle $\sF$ is stable. 
 We have to exclude the cases $(a,b,c_2) = (3,5,6)$, $(3,4,5)$, $(6,10,6)$ from above. Then the following Example~\ref{case22} completes the proof. 
Let $H = \eta + \phi^*\sO_{\PN_2}(2)$, i.e. $-K_X = 2H$ and $H = \psi^*H'$ for some $H' \in \Pic(X')$ with $-K_{X'} = 2H'$. We have ${H'}^3 = 7-c_2$. 

\vspace{0.2cm} 

1.) Assume $(a,b,c_2) = (3,5,6)$ or $(6,10,6)$. In these cases $H^3 = 1$ and hence by \cite{Shin}, 
 \[X' \lra W \subset \PN_6\] 
is a 2:1-covering of the Veronese cone $W$, moreover $|H'|$ has a single base point $p$ and the general member $S' \in |H'|$ is smooth in a neighborhood of $p$. In particular, $p \in X'_{reg}$, meaning $B = \psi(D)$ does not contain $p$. 
The blowup $X'_p = Bl_p(X')$ sits in the diagram
 \[\xymatrix{X'_p \ar[dr]\ar[rr]^{2:1} & & \PN(\sO_{\PN_2} \oplus \sO_{\PN_2}(2)) \ar[dl] & \hspace{-1cm} = Bl_p(W)\\            & \PN_2 &&}\]
The map $X'_p \to \PN_2$ is an elliptic fibration and the 2:1-covering $X'_p \to \PN(\sO_{\PN_2} \oplus \sO_{\PN_2}(2))$ is ramified along the disjoint union of the minimal section of $\PN(\sO_{\PN_2} \oplus \sO_{\PN_2}(2))$ and a 3:1-covering $K$ of $\PN_2$. Since $p \not\in B$ we may view $B$ also as a curve in $X'_p$. It is the singular locus of the surface $K$.

The restriction of $K$ to any fiber of $\PN(\sO_{\PN_2} \oplus \sO_{\PN_2}(2)) \lra \PN_2$ that meets $B$ gives a singular cubic polynomial, the points of intersection of the fiber with $B$ corresponding to multiple roots. A cubic polynomial can have at most one multiple root. Hence $B$ is mapped isomorphically onto its image in $\PN_2$, i.e., $B$ is a plane curve.

For $(a,b,c_2) = (3,5,6)$ we have
 \[g(B) = h^1(D, \sO_D) = h^2(X, \sO_X(-D)) = h^1(\PN_2, \sF(1)) = 2\]
by Riemann Roch. This is impossible. For $(a,b,c_2) = (6,10,6)$, we
have $b = 2a-2$ and $\phi^*\sO_{\PN_2}(1)\cdot l_{\psi} = 1$. Then $D \to B$ is a $\PN_1$-bundle, hence $K_D^2 = 8(1-g(B))$. This gives us $g(B) = 30$. On the other hand, $\deg(B) = -\frac{H.D^2}{2} = 16$, which is again not possible for a smooth plane curve.  
      
\vspace{0.2cm} 

2.) Assume $(a,b,c_2) = (3,4,5)$. Then $H^3 = 2$ and by \cite{Shin}, $|H'|$ is base point free, defining a 2:1-covering $X' \to \PN_3$, ramified along a quartic in $\PN_3$ with singularities along the image of $B$. The strict transform of the ramification divisor is a section in $2H-D = \eta$, which is absurd.
\end{proof}

\begin{example} \label{case22}
Let $x_0, \dots, x_4$ be homogeneous coordinates of $\PN_4$. Let $B = B_4$ be the rational normal curve of degree $4$ in $\PN_4$, cut out by the six $2\times 2$ minors of
  \[\left(\begin{array}{ccccc}
      x_0 & x_1 & x_2 & x_3 \\
      x_1 & x_2 & x_3 & x_4
    \end{array}\right),\]
i.e., by $q_0 = x_0x_2-x_1^2$, $q_1 = x_0x_3-x_1x_2$, $q_2 = x_0x_4-x_1x_3$, $q_3 = x_1x_3-x_2^2$, $q_4 = x_1x_4-x_2x_3$ and $q_5 = x_2x_4-x_3^2$. For $X'$, take the secant variety $S_1(B)$ to $B$, defined by the cubic determinant
 \[\left|\begin{array}{ccc}
      x_0 & x_1 & x_2 \\
      x_1 & x_2 & x_3 \\
      x_2 & x_3 & x_4
    \end{array}\right| \in H^0(\PN_4, \I_{B/\PN_4}^2(3)).\]
Then $X'$ is a canonical Gorenstein Fano threefold with singularities of type $cA_1$ along $B$. The blowup $X = Bl_B(X')$ is almost Fano with $-K_X = \psi^*\sO_{\PN_4}(2)|_X$. 

Denote the blowup map $Bl_B(\PN_4) \lra \PN_4$ by $\psi$ as well, and by $\hat{D}$ the exceptional divisor. The quadrics $q_0, \dots, q_5$ define a morphism $\phi: Bl_B(\PN_4) \lra \PN_5$ onto the smooth quadric 
  \[Q = q_1q_4 - q_2q_3 - q_0q_5.\]
The map $\phi$ is nothing but the blowup of $Q$ in a $\PN_2$, embedded into $\PN_5$ via Veronese
  \[[s:t:u] \mapsto [s^2:-su:u^2-st:st:-tu:t^2],\] 
(\cite{Harris}, p.90). Our $X$ is the exceptional divisor of $\phi$, i.e., $X = \bP(N_{\bP_2/Q}^*)$. The Chern classes of the rank two bundle $N^*_{\bP_2/Q}$ on $\bP_2$ are easily computed:
  \[c_1(N^*_{\PN_2/Q}) = -5H, \quad c_2(N^*_{\PN_2/Q})= 10.\]
Then $\sF =N^*_{\PN_2/Q}(2)$ has Chern classes $c_1(\sF) = -1$ and $c_2(\sF) = 4$. Clearly $\sF$ and $\sF(1)$ do not have sections.

The bundle $\sF$ gives a single point in the moduli space $M(-1,4)$. For $\sF$ general, $-K_X$ remains big and nef, but $\psi$ is small, contracting the (finitely many) jumping lines. 
\end{example}

\begin{example} 
Consider $\bP_4 = {\rm Proj}\,\KC[v, w, x, y, z]$. Let $l$ be the line defined by the homogeneous ideal $I_l = (x, y, z)$. Let $X'$ be the zero set of a general cubic in $I_l^2$ such that in particular $X'_{sing} = l$. Then $X = Bl_l(X')$ is an almost Fano threefold. If $\psi: X \to X'$ denotes the blowdown as usual, then $-K_X = \psi^*\sO_{\bP_4}(2)$ and $(-K_X)^3 = 24$. 

On $X$, the net $|\psi^*\sO_{\bP_4}(1)-D|$ is spanned, defining a $\bP_1$-bundle structure $\phi: X \to \bP_2$. Let $\sF$ be the normalized rank two vector bundle on $\bP_2$ such that $X = \bP(\sF)$. Since $-K_X$ is divisible, $c_1(\sF) = -1$. A general $S \in |{-}\frac{1}{2}K_X|$ is an almost del Pezzo surface, obtained by blowing up a cubic in $\bP_3$ in a single du Val point. Hence $K_S^2 = 3$ and $c_2(\sF(2)) = 6$. Then $[\sF] \in M(-1, 4)$. Since $\sF(1)$ has a section given by $D$, $\sF$ is a Hulsbergen bundle.

If the cubic defining $X'$ is general enough, its singularities along $l$ will be generically $cA_1$, with three dissident points at the intersection with the discriminant cubic $\Delta \subset \bP_4$. For special $X'$, for example, the zero set of 
  \[x^3 + vy^2 + wz^2,\]
we get $cA_2$ singularities, up to two dissident points of type $cD_4$. Here $\Delta = xvw$.
\end{example}

\

\noindent From now on we shall assume $\Delta \ne \emptyset.$ 
%4.6
\begin{proposition}
\begin{enumerate} 
\item $\alpha (-K_X)^3 + \beta (12-d) = 0. $
\item $$ D \cdot l_{\psi} = \beta \deg (\phi \vert l_{\psi}) \deg (\phi(l_{\psi})), $$
in particular $\beta = -1,-2.$ 
\item If $D$ is not a $\bP_1-$bundle over $B$, then $\beta = -1.$ 
\item If $\tilde \psi$ is birational, then $(-K_X)^3 \geq 4.$
\end{enumerate}
\end{proposition}

\begin{proof} (1) The equation $K_X^2 \cdot D = 0$ reads 
$$ \alpha (-K_X)^3 + \beta K_X^2 \cdot \phi^*(\sO(1)) = 0.$$
Then (1) follows via 
$$ K_X^2 \cdot \phi^*(\sO(1)) = 12-d.$$ \\ 
The first formula in (2) is clear and the second follows from $D \cdot l_{\psi} = -2.$ \\
(3) Suppose that $D$ is not a $\bP_1-$bundle. Recall that $D \cdot l_{\psi} = -2$. Moreover
the family $l_{\psi} $ splits, hence the image family in $\bP_2$ splits, i.e.  $\phi(l_{\psi})$ is a conic and by (2) we obtain
$\beta =  -1.$ \\
(4) Finally, if $(-K_X)^3 = 2,$ then $h^0(\sE) = 4,$ so that $\vert \zeta \vert $ maps to $\bP_3$ and $\tilde \psi$ cannot be birational.
\end{proof}

%4.7 
\begin{proposition} 
\begin{enumerate}
\item If $\tilde \psi $ is divisorial, then it contracts a unique divisor $\tilde D$ and nothing else;
and $\tilde D \cap X = D.$ 
\item If $\tilde \psi $ is divisorial, then $\lambda \geq 0, c_1 \leq 3$ and 
$$\alpha(c_1^2-c_2) + \beta c_1 = 0.$$ 
\item If $D$ extends to $\tilde D \in \vert \alpha \zeta + \pi^*(\sO(\beta)) \vert,$ then 
$$2\alpha (c_1^2 - c_2) + (\alpha (3-c_1) + 2 \beta)c_1 + (3-c_1) \beta = 0.$$
\item If $\alpha \leq 2$, then $D$ extends. 
\item If $\dim Z = 3,$ i.e. $\tilde \psi$ is a fibration, then $c_1^2 = c_2.$ Moreover $\lambda > 0,$ i.e. $c_1 \leq 2.$ 
\end{enumerate}
\end{proposition}

\begin{proof} (1) Only the second assertion needs a proof. Let $\tilde D$ be the unique exceptional divisor of $\tilde \psi.$ 
Using $c_1 \leq 3$ (see (2)), either $c_1 \leq 2$ and $\sE$ is a Fano bundle or $c_1 = 3$ and $-K_{\bP(\sE)}$ is big and nef but not
ample. In this second case, $X \cdot l_{\psi} = (2 \zeta + \pi^*(\sO(3-c_1))) \cdot l_{\psi} = 0.$ Since every curve which is contracted by
$\tilde \psi$, is proportional to $l_{\psi},$ it is either contained in $X$ or is disjoint from $X$. Hence $\tilde D \cap X = D.$ 
In the first case one can just use the classification or show directly that $\tilde D \cdot l_{\phi} = 2$ which easily implies the claim.

\vskip .2cm \noindent 
(2) The first assertion follows from $X \cdot l_{\psi} \geq 0$ and the second from $\zeta^3 \cdot \tilde D = 0.$  

\vskip .2cm \noindent 
(3) This comes from $0 = K_X^2 \cdot D = \zeta^2 \cdot X \cdot \tilde D = 0.$ 

\vskip .2cm \noindent
(4) The obstruction for extending $D$ is in
$$H^1(\bP(\sE), -2 \zeta - \pi^*(\sO(\lambda)) + \alpha \zeta + \pi^*(\sO(\beta))) = 0. $$
This group vanishes if $\alpha \leq 2.$  

\vskip .2cm \noindent
(5) is again clear from $X \cdot l_{\psi} > 0.$ 

\end{proof}

%4.8
\noindent A direct numerical consequence is 

\begin{corollary}  
\begin{enumerate} 
\item If $\tilde \psi$ is divisorial, then $c_1 = 3$ or $c_2 = 0.$ If $\beta = -1$, we have $(\alpha,c_1,c_2) = (1,3,6), (1,1,0).$
If $\beta = -2,$ we have $(\alpha,c_1,c_2) = (1,3,3)$, $(2,3,6)$, $(3,3,7)$, $(6,3,8)$, $(1,2,0). $
\item If $\tilde \psi$ is not birational, then $(c_1,c_2) = (1,1), (2,4).$
\end{enumerate}
\end{corollary}

\begin{proof} (1) Using (3.8)(2) and (3), we obtain 
$$ (3-c_1)(1 - {{c_1^2} \over {c_1^2-c_2}}) = 0, $$
hence $c_1 = 3 $ or $c_2 = 0.$ 
Now we just use 3.8(2) to obtain the listed cases, having in mind $c_1^2 > c_2$ since $\tilde \psi$ is birational. \\
(2) is obvious from (3.8)(5). 
\end{proof}

%4.9
\begin{corollary} 
Suppose $\tilde \psi$ is not birational. Then $\beta = -1$ and  $X$ is a complete intersection in 
$\bP_2 \times \bP_3$ of degree $(1,1),(2,2)$ resp. $(1,2),(2,1).$ However in these cases $\psi$ is small.
\end{corollary}

\begin{proof} Applying (3.9)(2) we are reduced to two cases. In both cases $\bP(\sE)$ is Fano and we can apply the classification [SW90]. \\
If $(c_1,c_2) = (1,1)$, then $\sE = T_{\bP_2}(-1) \oplus \sO$ and $\lambda = 2.$ Thus $\bP(\sE) \subset \bP_2 \times \bP_3$ is a 
divisor of degree $(1,1)$. To see that $\psi$ is small, suppose the contrary. Then $h^0(\sO_X(D)) = 1;$ therefore
$D = 2(-K_X) + \phi^*(\sO(-1))$ yields 
$$ h^0(S^2(\sE)(-1)) = 1, $$
which is clearly not true. \\
If $(c_1,c_2) = (2,4),$ then $\sE$ is given by an exact sequence
$$ 0 \to \sO(-2) \to \sO^4 \to \sE \to 0 $$
which realises $\bP(\sE)$ as a divisor of degree $(2,1).$ Again we see easily that $\psi$ is small. 
\end{proof}

\noindent So from now we may assume (and do) that $\tilde \psi$ is birational.

%4.10
\begin{proposition}
Suppose $\beta = -2$.  Then $X$ is one of the following.
\begin{enumerate} 
\item $\sE = \sO(2) \oplus \sO^2$ and $X \in \vert 2 \zeta + \pi^*(\sO(1)) \vert,$  (A.3,no.5);
\item $\sE = \sO \oplus \sF$ or a non-split extension 
$$ 0 \to \sO \to \sE \to \sF \to 0 $$ where $\sF$ is given by a non-split extension
$$ 0 \to \sO(2) \to \sF \to \sI_p(1) \to 0 $$
for some $p \in \bP_2,$ (A.3,no.6);
\item $\sE = \sO(2) \oplus T_{\bP_2}(-1),$ (A.3,no.6).

\end{enumerate}
In all cases but the first, $X \in \vert 2 \zeta \vert.$ 

\end{proposition}

\begin{proof} Since $\beta = -2,$ then by (3.7)(3), $D$ is a $\bP_1-$bundle over $B.$ Let $e$ be its invariant [Ha77,V.2] and $C_0$ a section with $C_0^2 = -e.$ 
\vskip .2cm \noindent
{\bf(A)} First suppose that $\tilde \psi$ is {\it divisorial}. By (3.9) we have six possible triples for $(\alpha,c_1,c_2).$ If $c_2 = 0,$ 
then $c_1 = 1,2$ so that $\sE$ is a Fano bundle. Hence by [SW90], $\sE = \sO(a) \oplus \sO^2$ with $a = 1,2.$ 
The case $a = 2$ certainly occurs,
while in case $a = 1$ we obtain $\beta = -1.$ \\
Now assume $c_1 = 3.$ In that case $-K_{\bP(\sE)}$ is spanned and nef, but not ample. Equivalently $\sE $ is nef with $\zeta$ big,
but $\sE$ is not ample (otherwise $-K_X = \zeta \vert X$ would be
ample or apply [SW90]). Notice that we may assume that $\sE$ is spanned as already observed in (3.2).
\vskip .2cm \noindent
(1) $c_1 = 3, c_2 = 3.$ Choose a general section of $\sE $ without zeroes to obtain an exact sequence of vector bundles
$$ 0 \to \sO \to \sE \to \sF \to 0. \eqno (S)$$ 
Then $\sF$ is a spanned rank 2-bundle with $c_1 = 3, c_2 = 3$. First notice that $\sF$ is not ample. In fact, if $\sF$ is ample,
then it is a Fano bundle, hence $\sF = T_{\bP_2}.$ If $(S)$ splits, then $\tilde  \psi$ would be small, contracting simply $\bP(\sO).$
If $(S)$ does not split, then $\sE$ would be ample (as jet bundle), contradiction. \\
So $\sF$ is not ample. Hence $-K_{\bP(\sF)}$ is big and spanned, so that (3.4) applies and $\sF$ is given by an extension 
$$ 0 \to \sO(2) \to \sF \to \sI_p(1) \to 0. \eqno (S') $$
If $(S)$ splits, then $\tilde \psi $ contracts exactly $\tilde D = \bP(\sO \oplus \I_p(1)),$ providing an example. \\
If $(S)$ does not split, we argue as follows.

The divisor $\tilde D$ is given by $\tilde D = \bP(\sG) $, where the torsion free sheaf $\sG$ is given by 
$$ 0 \to \sO(2) \to \sE \to \sG \to 0.$$ 
In fact, $ \tilde D = \zeta + \pi^*(\sO(-2))$ since $\beta = -2$ and $\alpha = 1.$  
Notice $c_1(\sG) = 1$ and $c_2(\sG) = 1.$ 
Consider the exact sequence 
$$ 0 \to \sG \to \sG^{**} \to Q \to 0,$$
where $\sG^{**}$ is locally free and $Q$ supported on a finite set. Then 
$$  c_2(\sG^{**}) = c_2(\sG) - l(Q) = 1 - l(Q). $$
Now $\sG^{**}$ has sections vanishing in codimension 2 or nowhere, hence $c_2(\sG^{**}) \geq 0$ so that $l(Q) \leq 1.$ 
If $l(Q) = 0,$ then $\sG$ is already locally free and $\sG_l = \sO(1) \oplus \sO$ for all lines, so that $\sG = \sO(1) \oplus \sO$ or
$T_{\bP_2}(-1)$ by the classification of uniform bundles. The first case is ruled out by $c_2(\sG) = 1,$ but the second of course exists. \\
If $l(Q) = 1,$ then $c_2(\sG^{**}) = 0$ and therefore $\sG^{**}$ is spanned, actually $\sG^{**} = \sO(1) \oplus \sO.$ 
This case exists, too; here $\sF$ is defined by $(S')$ and $\sE$ by $(S)$. The sheaf $\sG$ is given by an extension
$$ 0 \to \sO \to \sG \to \sI_p(1) \to 0.$$

\vskip .2cm \noindent
(2) $c_1 = 3, c_2 = 6. $ 
Again the divisor $\tilde D$ is given by $\tilde D = \bP(\sG) $, with
$$ 0 \to \sO(2) \to \sE \to \sG \to 0.$$ 
Here $c_1(\sG) = 1$ and $c_2(\sG) = 4.$ Since $\sG^{**}$ is generated outside a finite set, it is nef, hence $\sG^{**}_l = \sO(1) \oplus \sO$
for all lines $l \subset  \bP_2$. Using the classification of uniform bundles as above and taking into account $c_2(\sG) = 4,$ we obtain in the
same notations as above that $l(Q) = 4$ and $\sG^{**} = \sO(1) \oplus \sO$, or $l(Q) = 3$ and $\sG^{**} = T_{\bP_2}(-1)$.
In both cases we immediately see that $\sG$ cannot be spanned.

\vskip .2cm \noindent 
(3) $c_1 = 3, c_2 = 7.$ 
This case is ruled out as in (2). 
\vskip .2cm \noindent
(4) $c_1 = 3, c_2 = 8.$ Here $h^0(\sE) = 3 = h^0(-K_X)$ by Riemann-Roch, so that $\psi$ cannot be birational. 

\vskip .2cm \noindent
{\bf (B)} Now suppose that $\tilde \psi $ is small. 
\vskip .2cm \noindent
(1) If $\alpha \geq 3$ we argue as follows.
By (3.7)(1) we have $$ (-K_X)^3 = {{24 - 2d} \over {\alpha }}. $$
On the other hand, we must have $h^0(\sE) \geq 5,$ otherwise $\tilde \psi$ cannot be birational. \\
First suppose $h^0(\sE) \geq 6.$ 
This gives $\alpha = 3 $ since we assume $d > 0.$  Moreover $d = 3, (-K_X)^3 = 6$ and $h^0(\sE) = 6.$ Also notice that $D$ is a
$\bP_1-$bundle over $B$ since $\beta = -2.$ 
Via the adjunction formula and again using the notation $g = g(B)$ we obtain
$$ 8(1-g) = K_D^2 =  (K_X+D)^2 \cdot D = - 160,$$
hence $g = 21.$ This shows that $\phi_D$ must be finite because otherwise the exceptional section
$C_0$ is contracted by $\phi_D$ and therefore $C_0 $ is rational contradicting $g = 21.$ If $\phi_D$ is finite, we must argue
differently. Notice first
$$ g = h^1(\sO_D) = h^2(\sO_X(-D)) = h^1(K_X+D) = $$
$$ =h^1(-2K_X+\phi^*(\sO(-2))) = h^1(S^2(\sE(-1))).$$ 
Via the exact sequence
$$ 0 \to \sI_X \otimes (2 \zeta \otimes \pi^*(\sO(-2))) \to 2 \zeta \otimes \pi^*(\sO(-2)) \to -2K_X \otimes \phi^*(\sO(-2)) \to 0$$
and 
$$ H^0(-2K_X \otimes \phi^*(\sO(-2))) = 0$$
(since $\alpha = 3$) and
$$ H^0(\sI_X \otimes (2 \zeta \otimes \pi^*(\sO(-2)))) = H^0(\sO_{\bP_2}(c_1-5)),$$
we obtain
$$ g = h^1(S^2(\sE(-1))) = - \chi(S^2(\sE(-1))) + h^0(\sO(c_1-5)). \eqno (E)$$
Now Riemann-Roch gives
$$ \chi (S^2(\sE(-1))) = 3c_1^2 - 2c_1 - 5c_2. \eqno (RR)$$
By $ 6 = h^0(\sE) = \chi(\sE) $ and with Riemann-Roch we conclude
$$ c_2 = {{c_1^2} \over {2}} + {{3} \over {2}} c_1 - 3. $$
Putting this into (RR) yields 
$$ \chi (S^2(\sE(-1))) = {{c_1^2} \over {2}} - {{19} \over {2}} c_1 + 15$$
and finally with (E):
$$ {{c_1^2} \over {2}} - {{19} \over {2}} c_1 + 15 = {{(c_1-3)(c_1-4)} \over {2}} - 21.$$
This leads to $c_1 = 5$ and $c_2 = 17.$ 
Now consider the exact sequence 
$$ H^0(3 \zeta \otimes \pi^*(\sO(-2))) \to H^0(\sO_X(D)) \to H^1(3 \zeta \otimes \pi^*(\sO(-2)) \otimes \sO(-X)) = $$
$$ = H^1(\sE) =  H^1(-K_X) = 0.$$
The sequence shows that $D$ extends, hence we may apply 3.8(3) and obtain a contradiction for the specific values $c_1 = 5$ and $c_2 = 17.$ 
 \\ 
Finally if $h^0(\sE) = 5$, then  $\alpha = 3,4,5$. The case $\alpha = 3$ is excluded in the same way as before. 
If $\alpha = 4,5$, probably the same can be done, but it is more convenient to argue as follows. Since $(-K_X)^3 = 4,$ in our standard
notation either $W \subset \bP_4$ is a quartic and $\mu$ is an isomorphism or $W$ is a quadric and $\mu $ has degree 2. 
In the first case, $B$ is the singular locus of $X' \subset \bP_4$ and cut out by cubics so that $\sI_B(3)$ is spanned. Hence
$-3K_X - D$ is spanned which implies $\alpha \leq 3. $ \\
In the second case we observe that $h^0(\sI_B(2)) \ne 0,$ since $\mu$ is ramified along a quartic. Hence $H^0(-2K_X-D) \ne 0.$
This contradicts $\alpha = 4,5.$ 

\vskip .2cm \noindent

\vskip .2cm \noindent 
(2) Suppose now that $\alpha \leq 2.$  Hence $D$ extends 
to $\tilde D$, but $\tilde D$ is not contracted
by $\tilde \psi.$ Using $N^*_{X/\bP(\sE)} \vert D = N^*_{D/\tilde D}$ and  the conormal bundle sequence
$$ 0 \to N^*_{D/\tilde D} \vert l_{\psi} \to N^*_{l_{\psi }/\tilde D} \to N^*_{l_{\psi}/D} = \sO \to 0, $$
we conclude that $X\cdot l_{\psi} < 0$, because otherwise $N_{l_{\psi} /\tilde D}$ is nef and therefore $l_{\psi} $ moves in
$\tilde D$ in a covering family. This is only possible if $\tilde \psi$ contracts $\tilde D.$ \\
So $X \cdot l_{\psi} < 0$ and from $X = 2 \zeta + \pi^*(\sO(\lambda))$, it follows $\lambda < 0,$ i.e. 
$$ c_1 = c_1(\sE) \geq 4.$$
\vskip .2cm \noindent 
(2.a) If $\alpha = 2,$ then $(-K_X)^3 = 12 - d$, so that $(-K_X)^3 = 4,6,8,10$ so that $h^0(\sE) = 5,6,7,8.$
Now (3.8(3)) gives
$$ c_1^2 - 2c_2 + 2c_1 - 3 = 0.$$
By Riemann-Roch for $\chi(\sE)$ we obtain   
$$ c_1^2 - 2c_2 + 3c_1 = 4,6,8,10. $$ 
Both quadratic equations yield $c_1 = 1,3,5,7.$ Since $c_1 \geq 4$ we
end up either with
$$ c_1 = 5, c_2 = 16, h^0(\sE) = h^0(-K_X) = 7, (-K_X)^3 = 8 $$
or with
$$ c_1 = 7, c_2 = 30, h^0(\sE) = h^0(-K_X) = 8, (-K_X)^3 = 10.$$
Suppose first $c_1 = 5.$ 
With the same computations as in (1) we obtain $8(1-g) = (K_X+D)^2 \cdot D = -32,$ so $g = 5.$ 
On the other hand, as in (1),
$$ g = h^1(\sO_D) = h^2(\sO_X(-D)) = h^1(K_X+D) = h^1(-K_X \otimes \phi^*(\sO(-2))) = h^1(\sE(-2)).$$
Now Riemann-Roch shows $\chi(\sE(-2)) = -6$. Since $h^0(\sE(-2)) = h^2(\sE(-2)) = 0,$ we obtain $h^1(\sE(-2)) = 6,$
contradiction. \\
The case $c_1 = 7$ is ruled out in the same way. 
\vskip .2cm \noindent
(2.b) $\alpha  = 1.$ \\
Since $D = -K_X + \phi^*(\sO(-2)),$ we have $h^0(\sE(-2)) = 1,$ and, more generally
$$ h^0(S^m(\sE(-2))) = 1 $$
for all positive integers $m.$ This already shows $\lambda \geq -3.$ 
As in (1), we compute $8(1-g) = 8,$ hence $g = 0.$ Therefore $ h^1(\sE(-2)) = 0$ so that 
$$ \chi(\sE(-2)) = 1.$$ 
By Riemann-Roch we obtain
$$ c_1^2 - 2c_2 - c_1 = 2.$$ 
On the other hand, (3.8)(3) yields
$$ c_1^2 - 2c_2 + c_1 = 6.$$
Hence $c_1 = 2, c_2 = 0.$ This implies $\sE = \sO(2) \oplus \sO^2$ by [SW90]. But then $\tilde \psi$ is divisorial; namely it
contracts $\tilde D = \bP(\sO^2).$
\end{proof}

%4.11
\begin{proposition} Suppose that $\beta = -1$ and $\tilde \psi$ birational.
Then $X$ is one the following.
\begin{enumerate}
\item $X$ is a divisor in $\bP(\sO(1) \oplus \sO^2) $ over  $\bP_2$ in the linear system $\vert 2 \zeta + \pi^*(\sO(2))\vert,$ (A.3,no.7);
\item Consider a vector bundle $\sE$ over $\bP_2$ given by the exact sequence
$$ 0 \to \sO(-2) \to \sO(1) \oplus \sO^3 \to \sE \to 0;$$ 
then 
$$ X \in \vert 2 \zeta \vert,$$ (A.3,no.8). 

\end{enumerate}
\end{proposition}

\begin{proof}  By (3.7)(1)  we have
$$ \alpha (-K_X)^3 = 12 - d. $$  
This reduces already by $d > 0$ to $(-K_X)^3 \leq 10$. We also know $(-K_X)^3 \geq 4$ (3.7(4)) so that $\alpha \leq 2$. \\
(1) Suppose that $\tilde \psi$ is divisorial. \\
By (3.9) $\alpha = 1,$ and $c_1 = 1,3.$ In the first case $c_2 = 0$ and 
by the classification of the Fano bundles, $\sE = \sO(1) \oplus \sO^2,$ which is one of the cases listed. \\
If $c_1 = 3,$ then $c_2 = 6,$ then $h^0(\sE) = 6$ and $\zeta^4 = 3$ so that $\bP(\sE)$ maps to a cubic in $\bP_5,$ 
and $\deg \tau  = 1.$ 
The exact sequence
$$ 0 \to H^0(\zeta - \tilde D) \to H^0(\zeta) \to H^0(\zeta \vert \tilde D) \to H^1(\zeta - \tilde D) = 0$$
shows together with $\zeta - \tilde D = \pi^*(\sO(1))$ that $h^0(\zeta \vert \tilde D) = 3$, so that $\tau \tilde \psi(\tilde D)$ is
a plane or a curve in a plane. It cannot be a curve, since $\zeta^2 \cdot \tilde D^2 \ne 0.$ 
So $\tau \tilde \psi(\tilde D)$ is a plane and therefore $\bP(\sE) $ is the blowup of its image in $\bP_5$ along a plane. Thus
$$ \bP(\sE) \subset \bP(\sO_{\bP_2}(1) \oplus \sO_{\bP_2}^3),$$
coming from an epimorphism
$$ \sO(1) \oplus \sO^3 \to \sE \to 0.$$
Conversely, if $\sE$ is given by a sequence 
$$ 0 \to \sO(-2) \to \sO(1) \oplus \sO^3 \to \sE \to 0,$$
then $\sE$ is spanned. Let $\tilde D = \bP(\sE) \cap \bP(\sO^3)$ and take 
$$X \in \vert 2 \zeta_{\sE} \vert  $$
general. Then $-K_X = \zeta_{\sE} \vert X $ is big and nef, and $K_X^2 \cdot D = 0$ so that $-K_X$ is not ample.

\vskip .2cm \noindent
(2) $ \tilde \psi$ is small. \\
\noindent
By 3.8(4) the extension $\tilde D$ exists, but $\tilde \psi$ does not contract $\tilde D.$ 
\vskip .2cm \noindent
(2.a) The case $\alpha = 2$ is ruled out as follows.
First apply 3.8(3) to obtain
$$ 2c_1^2 - 4c_2 + 5c_1 - 3 = 0.$$ 
Since $\alpha = 2,$ we have $(-K_X)^3 = 4,$ hence $h^0(\sE) = 5.$ Thus Riemann-Roch gives
$$ c_1^2 - 2c_2 + 3c_1 = 4.$$ Comparing both formulas gives $c_1 = 5$ and $c_2 = 18.$ 
Thus $\zeta^4 = c_1^2 - c_2 = 7$  and so $\zeta^3 \cdot X = \zeta^3 \cdot (2 \zeta + \pi^*(\sO(-2))) = 4.$ 
Now $\vert \zeta \vert$ maps $\bP(\sE) $ to $\bP_4$ such that the image $W$ of $X$ has degree $2$ or $4$.  
This case is ruled out by the method of Step 1 and 2 of Theorem 4.6 below. To be a little more precise, let us consider
only the non-hyperelliptic case. We first compute 
$$32 = (2 (-K_X))^3 = D^3 + 3 D^2 \cdot \phi^*(\sO(1)) + 3 D \cdot \phi^*(\sO(1))^2. \eqno (*) $$
Then we consider the blowup $\hat \bP_4 \to \bP_4$ along $B$ with exceptional divisor $\hat D$ whose restriction to $X \subset \hat \bP_4$
is of course just $D.$ We compute easily
$$ g(B) = h^1(\sO_X(D)) = h^1(\sE(-1)) = 2. $$
Next we compute 
$$D^3 = \hat D^3 \cdot X = -6 \deg B - 4(g-1) $$  as in (4.6). 
By $g = 2,$ we obtain $\deg B = 8.$ 
Now $X \in \vert \hat \psi^*(\sO(4)) - 2 \hat D \vert.$ Since $ \hat \psi^*(\sO(4)) - 2 \hat D \vert X = \phi^*(\sO(2))$ and since
$H^1(\sO_X) = 0,$ we conclude that $\hat \psi^*(\sO(4)) - 2 \hat D$ is spanned. On the other hand, we compute easily that 
$(\hat \psi^*(\sO(4)) - 2 \hat D)^4 = -96,$ a contradiction. \\
The hyperelliptic case works similarly.
\vskip .2cm \noindent
(2.b) If $\alpha = 1,$ then $(-K_X)^3 $ can take the values $4,6,8,10 $ so that $h^0(\sE) = 5,6,7,8.$ 
Hence Riemann-Roch yields 
$$ c_1^2 - 2c_2 + 3c_1 = 4,6,8,10.$$ 
On the other hand, 3.8(3) gives
$$ c_1^2 - 2c_2 + 2c_1 = 3.$$ 
Putting things together, we obtain $(c_1,c_2) = (1,0), (3,6), (5,16), (7,30).$
Since $\zeta \cdot \tilde D > 0,$ we have $c_1^2 - c_2 - c_1 > 0,$ which deletes the first two cases.  
\vskip .2cm \noindent We next exclude the case $(7,30).$ Suppose the contrary and suppose also that the general 
$M \in \vert 2 \zeta + \pi^*(\sO(-2)) \vert$ is irreducible and reduced. Then we conclude from $X \in \vert 2 \zeta + \pi^*(\sO(-4)) \vert $ that
$$ h^0( 2 \zeta + \pi^*(\sO(-2))) > h^0(\pi^*(\sO(2))) + 1,$$  
since $M$ is not in the space generated by $X + \pi^*(H^0(\sO(2))) $ and $2 \tilde D.$ 
Since $2 \zeta + \pi^*(\sO(-2)) \vert X = 2D,$ we have $h^0(2 \zeta + \pi^*(\sO(-2) \vert X)) = 1$, which yields a contradiction by
taking cohomology of
$$ 0 \to 2 \zeta + \pi^*(\sO(-2)) - X \to 2 \zeta + \pi^*(\sO(-2)) \to 2 \zeta + \pi^*(\sO(-2)) \vert X \to 0.$$ 
It remains to show that the general $M$ is irreducible and reduced. Suppose the contrary and write
$$ M = \sum a_i M_i.$$ 
We see immediately that after possibly renumbering, we have $M_1 \in \vert 2 \zeta + \pi^*(\sO(a)) \vert  $ with $a \leq -3, a_1 = 1$ and
$M_j \in \vert \pi^*(\sO(b_j)) \vert $ for $j \geq 2$ with $b_j \geq 1.$ Hence the general element of $\vert 2 \zeta + \pi^*(\sO(a)) \vert $
is irreducible and reduced. This contradicts the exact sequence
$$ 0 \to H^0(2 \zeta + \pi^*(\sO(a)) - X) \to H^0(2 \zeta + \pi^*(\sO(a))) \to H^0(2 \zeta + \pi^*(\sO(a)) \vert X) = 0$$
(notice $  2 \zeta + \pi^*(\sO(a)) \vert X = 2D - \phi^*(\sO(a+2)))$. 

\vskip .2cm \noindent
In the last case $(5,16),$ observe that $\zeta^4 = 9.$ Hence $\deg \tilde W \subset \bP_6$ has degree 3 or 9. In the second case
$\deg \tau = 1,$ and $h^0(2 \zeta - X) \subset  H^0(\sI_W(2)).$ 
Since $X'$ not hyperelliptic, $X' \subset \bP_6$ is cut out by three quadrics, hence 
$h^0(2 \zeta - X) = 3.$ But $2 \zeta - X = \pi^*(\sO(2)),$ contradiction. \\
If $\deg \tau = 3,$ then $\tilde W $ is ``minimal'', hence by classification, $\tilde W$ is the image of a scroll $\bP(\sF) $ 
with $\sF $ a spanned bundle over $\bP_1$ under the morphism defined by the tautological bundle, see e.g. [IP99,2.2.11]. 
Moreover the singular locus of $\tilde W$, which certainly contains $B$, is linear. This contradicts 
$\deg B = {{1} \over {2}} K_X \cdot D^2 = 3.$ 
\end{proof}

\noindent In conclusion we obtain the following 

%4.12
\begin{theorem} Let $X$ be a smooth projective threefold with $-K_X$ big and nef but not ample. Assume $\rho(X) = 2$ and that
$X$ is a proper conic bundle $\phi: X \to \bP_2$, i.e. the
discriminant locus is non-empty. Assume that the anticanonical
morphism is divisorial and $\sE = \phi_*(-K_X)$ is spanned. Then $X$ is one of the
following: 
\begin{enumerate}
\item $X$ is a divisor in $\bP(\sO(1) \oplus \sO^2) $ over  $\bP_2$ in the linear system $\vert 2 \zeta + \pi^*(\sO(2))\vert,$ (A.3,no.7);
\item Consider a vector bundle $\sE$ given by the exact sequence
$$ 0 \to \sO(-2) \to \sO(1) \oplus \sO^3 \to \sE \to 0;$$ 
then 
$$ X \in \vert 2 \zeta \vert,$$ (A.3,no.8);
\item $\sE = \sO(2) \oplus \sO^2$ and $X \in \vert 2 \zeta + \pi^*(\sO(1)) \vert,$  (A.3,no.5);
\item $\sE = \sO \oplus \sF$ or a non-split extension 
$$ 0 \to \sO \to \sE \to \sF \to 0 $$ where $\sF$ is given by a non-split extension
$$ 0 \to \sO(2) \to \sF \to \sI_p(1) \to 0 $$
for some $p \in \bP_2$; here $X \in \vert 2 \zeta \vert,$ (A.3,no.6);
\item $\sE = \sO(2) \oplus T_{\bP_2}(-1)$ and $X \in \vert 2 \zeta \vert,$ (A.3,no.6). 
\end{enumerate}
\end{theorem}

\

%%%%%%%%%%%%%%%%%%%%%

\section{Blowdown to a curve} \label{blowdowncurve}
\setcounter{lemma}{0}

\begin{abs} 
{\bf Setup.} In this section $\phi: X \to Y$ is the blowup of a smooth Fano threefold $Y$ with $\rho(Y) = 1$ in a smooth curve $C$. As usual, $X$ is a smooth almost Fano threefold, such that the anticanonical map $\psi: X \to X'$ is divisorial. Let $H \in \Pic(Y)$ be the fundamental divisor and $r$ the index of $Y$, i.e. 
 \[-K_Y = rH\] 
in $\Pic(Y)$. Let $E \subset X$ be the exceptional divisor of $\phi$,
so that $E \simeq \PN(N^*_{C/Y})$ is a ruled surface. Then $\Pic(X)$ is generated by $\phi^*H$ and $E$. Let
 \[g_C = g(C) \quad \mbox{ and } \quad d = \deg(C) = H\cdot C\]
the genus and degree of $C$, respectively. We have 
 \[-K_X = \phi^*(rH) -E\] 
and $\phi^*H^2\cdot E = 0$, $\phi^*H\cdot E^2 = -d$, $E^3 = -\deg(N_{C/Y}) = rd + 2g_C-2$. Therefore
 \[({-}K_X)^3 = r^3H^3 -2rd + 2g_C-2.\]
The restriction $-K_X|_E = \sO_E(1) \otimes \phi^*(rH)|_E$ is still nef and also big, since $E \not= D$. Hence $\deg(N^*_{C/Y}\otimes (rH)) > 0$, implying
 \begin{equation} \label{absch}
   2g_C-2 < rd, \quad d < r^2H^3,
 \end{equation}
where the second inequality follows from $2rd = r^3H^3 + K_X^3 +2g_C-2 < r^3H^3 +rd$. 

\vspace{0.2cm}

We note moreover 
 \begin{equation} \label{pic}
   \Pic(X') = \KZ \cdot ({-}K_{X'}).
 \end{equation}
Indeed: by assumption, the Picard number of $X'$ is one. Let $\Pic(X') = \KZ \cdot L$ for some line bundle $L$ on $X'$. Then $-K_{X'} = kL$ for some integer $k$, hence $-K_X = k\psi^*L$. The intersection with an exceptional fiber $l_{\phi}$ of $\phi$ gives $1 = -K_X\cdot l_{\phi} = k\psi^*L\cdot l_{\phi}$, hence $k = 1$. This proves the claim.
\end{abs}

The following Lemma shows how we find examples:

\begin{lemma}
 Let $Y$ be a smooth Fano threefold of index $r$. Assume that $H$ is
 generated. Let $X = Bl_C(Y)$ be the blowup of $Y$ in a smooth curve $C$. Then $-K_X$ is generated if and only if the curve $C$ is cut out by hypersurfaces of degree $\le r$.  
\end{lemma}

\begin{proof}
If the curve $C$ is cut out by hypersurfaces of degree $\le r$, then ${\mathcal I}_C \otimes \sO_Y(rH)$ is generated. Then also $\sO_X(-E) \otimes \phi^*\sO_Y(r) \simeq \sO_X({-}K_X)$ is generated. Conversely, if $-K_X$ is generated, then we have a surjection
  \[\sO_X^{\oplus n} \lra \sO_X({-}K_X) \simeq \sO_X(-E) \otimes \phi^*\sO_Y(rH) \lra 0.\]
The kernel is a vector bundle ${\mathcal E}$ on $X$. If we can show $R^1\phi_*{\mathcal E} = 0$, then $C$ is cut out by hypersurfaces of degree $\le r$.

To this end let $l_{\phi} \simeq \PN_1$ be a fiber of $E \to C$. Then $N^*_{l_{\phi}/X} \simeq \sO_{\PN_1} \oplus \sO_{\PN_1}(1)$ and ${\mathcal E}|_{l_{\phi}} \simeq \sO_{\PN_1}^{\oplus n-2} \oplus \sO_{\PN_1}(-1)$. We conclude
 \[H^1(S^k N^*_{l_{\phi}/X} \otimes {\mathcal E}) = 0\]
for any $k > 0$. The formal function theorem gives $R^1\phi_*{\mathcal E} = 0$.
\end{proof}

%%%%%%%%%%%%%%%%%%%%%%%%%%%%%%%%

\

\noindent {\bf (A) Assume $\psi$ contracts $D$ to a curve.} Then by Proposition~\ref{Wilson}, $X$ is the blowup of $X'$ along the smooth curve $B = \psi(D)$ and $X'$ is a Gorenstein Fano threefold with cDV singularities along $B$. We have $\rho(X') = 1$. 

\begin{abs}
\noindent {\bf Some numerical data.} Since $-K_X = \psi^*({-}K_{X'})$ and $D$ is contracted to a curve, we have
 \[K_X^2\cdot D = 0 \quad \mbox{ and } \quad K_X\cdot D^2 > 0.\]
 Since $\Pic(X)$ is generated over $\KZ$ by $E$ and $\phi^*H$, we may write
 \[D = \alpha\phi^*H -\beta E\]
in $\Pic(X)$ for some $\alpha,\beta \in \KZ$. Since $\phi$ is an extremal contraction, we have $D \not= E$, hence $D\cdot l_{\phi} \ge 0$ for a fiber $l_{\phi}$ of $E \to C$. This shows $\beta \ge 0$. On the other hand, $D$ is effective, but exceptional, hence 
 \[\alpha, \beta > 0.\] 

Let $l_{\psi}$ be the general exceptional fiber of $\psi$ as usual. Then $-K_X\cdot l_{\psi} = 0$ gives $r\phi^*H\cdot l_{\psi} = E\cdot l_{\psi}$, hence $-2 = D\cdot l_{\psi} = (\alpha-\beta r)\phi^*H\cdot l_{\psi}$ by Proposition~\ref{Wilson}. We define
 \[\epsilon := \beta r -\alpha; \quad \mbox{ then } \epsilon = 1, 2.\]
From $K_X^2\cdot D = 0$ and the formula for $({-}K_X)^3$ from above we obtain
 \begin{eqnarray} \label{ab}
   \beta ({-}K_X)^3 & = & \epsilon (r^2H^3-d),\\\nonumber
    K_X\cdot D^2       & = & \epsilon (\alpha rH^3 -\beta d),\\\nonumber
   \beta (2g_C-2) & = & (\beta r +\alpha)d -\alpha r^2H^3.
 \end{eqnarray}
\end{abs}

\begin{lemma} \label{finite}
Assume $\psi|_E: E \to \psi(E)$ is not finite. Then $C \simeq \PN_1$ is either a line or a conic in $Y$ with $N_{C/Y} = \sO_{\PN_1}(rd) \oplus \sO_{\PN_1}(-2)$. In particular, $r \le 2$.
\end{lemma}

\begin{proof}
We have already seen that $-K_X$ is still free and big on $E$. Hence, if $\psi|_E$ is not finite, then the minimal section $C_0$ of the ruled surface $E$ is contracted by $\psi$. This means $C_0$ is a fiber or a component of a fiber of $D \to B$, hence $C_0 \simeq \PN_1$ with $D\cdot C_0 = -1, -2$ by Proposition~\ref{Wilson}. With the above notation we have 
 \[D\cdot C_0 = (\beta({-}K_X) - \epsilon \phi^*H)\cdot C_0 = - \epsilon \phi^*H\cdot C_0,\] 
since ${-}K_X$ is trivial on $C_0$ by assumption. We conclude $d \le 2$ and $N^*_{C/Y}(rd)$ is nef, but not ample. The adjunction formula yields $\deg N_{C/Y} = dr -2$, completing the proof.
\end{proof}

\begin{abs} \label{lines}
\noindent {\bf Lines and Conics.} 1.)~Assume $C \subset Y$ is a line. For $r \ge 3$, the blowup $X = Bl_C(Y)$ is a Fano threefold. For $r = 1$, we obtain an almost Fano threefold, where $\psi$ is small (\cite{AG5}, Corollary~4.3.2). Assume $r = 2$. For $H^3 \ge 3$ again $X$ is Fano (\cite{AG5}, Proposition~3.4.1). Hence assume $H^3 \le 2$. Then the third line in (\ref{ab}) reads $-2\beta = 2\beta + \alpha -4\alpha H^3$, hence $4\mid \alpha$. We conclude $\epsilon = 2$. Using $(-K_X)^3 = 8H^3 -6$, the first line of (\ref{ab}) reads $\beta (4H^3 -3) = 4H^3 -1$. For $H^3 = 2$ we obtain $5\beta = 7$, which is impossible; for $H^3 = 1$ we get $\beta = 3$ and hence $\alpha = 4$. Since $\epsilon = 2$, we have $\phi^*H\cdot l_{\psi} = 1$, hence $D$ is a smooth ruled surface over $B$. By the following Lemma~\ref{B}, $g_B = 2$, hence $K_D^2 = 8(1-g_B) = -8$, contradicting $K_D^2 = (K_X+D)^2\cdot D = -24$ by adjunction formula.

2.)~Assume $C \subset Y$ is a conic. Again, $X = Bl_C(Y)$ is Fano for $r \ge 3$. Assume $r = 1$. First note that $-K_X$ is not big for $H^3 \le 6$. Then \cite{AG5}, Corollary~4.4.3 says: $-K_X$ is big and nef for $H^3 \ge 8$; $\psi$ is small for $H^3 \ge 12$ and $C$ general; $\psi$ will always be small for $H^3 \ge 16$. We obtain no.~1 and 3 in table~\ref{tabdivcurve} for $H^3 = 8$ and $10$ (see \cite{AG5}, p.86 for $H^3 = 10$). The case $H^3 = 12$ is excluded by the first line of (\ref{ab}); for $H^3 = 14$ we obtain $\epsilon = 2$, $\beta = 3$ and $\alpha = 1$. Then $D$ is a smooth ruled surface over $B$, but the adjunction formula yields $K_D^2 = -18$, which is not divisible by $8$. Assume now $r = 2$. The third line in (\ref{ab}) gives $\alpha(2H^3 -1) = 3\beta$, which is together with $\epsilon = 2\beta - \alpha$ only possible for $H^3 = 2$. Here we have $\alpha = \beta = \epsilon$. The case $\epsilon = 2$ is again impossible by Lemma~\ref{B}. We obtain no.~5 in table~\ref{tabdivcurve}. 
\end{abs}

From now on we will assume $C$ is neither a line nor a conic. Then $\psi|_E$ is finite by Lemma~\ref{finite}, i.e., the image of $\psi$ still contains a one-dimensional family of disjoint lines. If $X'$ is hyperelliptic, this means $W$ cannot be the Veronese cone, hence Proposition~\ref{hyp} implies here

\begin{corollary} \label{hypfinite}
Assume $X'$ is hyperelliptic and $C$ is neither a line nor a conic. Then $({-}K_X)^3 \le 4$.
\end{corollary}

\begin{abs}
\noindent {\bf Mukai's classification.} \label{Mukaiclass}
In \cite{Mukai}, Mukai describes the embedding of all Gorenstein Fano threefolds with canonical singularities, such that the anticanonical divisor does not admit a {\em moving decomposition} (see table~\ref{Mukaitab}): by definition, a linear system $|L|$ on a normal projective variety admits a moving decomposition, if $L \sim A + B$, where $A$ and $B$ are Weil divisors, such that $|A|$ and $|B|$ are both of positive dimension.

In our case, $-K_{X'}$ generates $\Pic(X')$ as seen in (\ref{pic}). This means, $|{-}K_{X'}|$ cannot admit a moving decomposition, whenever $X'$ is factorial. With notations from \cite{KoMo}, the map $\psi: X \to X'$ is a (divisorial) log contraction of the klt pair $(X, \nu D)$ for any $0 < \nu < 1$ rational. By \cite{KoMo}, Corollary~3.18, $X'$ is $\KQ$-factorial. A closer look at the proof shows in fact $X'$ factorial whenever $\phi^*H\cdot l_{\psi} = 2$. This means, $\epsilon = 1$ implies $X'$ factorial and that Mukai's classification applies. 
\end{abs}

\begin{lemma} \label{B} Let $H'$ be a general hyperplane section of $\PN_{g+1}$. Then
 \begin{enumerate}
   \item $B$ is a smooth curve of degree $\mu^*H'\cdot B = \frac{1}{2}K_X\cdot D^2$ and genus 
 \[g_B = 1 - \frac{2\alpha}{r} + \frac{\epsilon d}{12} - \frac{1}{4}K_X\cdot D^2 - \frac{1}{6}D^3,\]
where $D^3 = \alpha^3H^3 + \beta^2(\beta r -3\alpha)d + 2\beta^3(g_C -1)$.
  \item If $\epsilon = 2$, then $D^3 = 4d - \frac{48\alpha}{r}$.
 \end{enumerate}
\end{lemma}

\begin{proof}
Consider for $m \gg 0$ the twisted ideal sequence of $D$ in $X$
 \[0 \lra \sO_X(-mK_X-D) \lra \sO_X(-mK_X) \lra \sO_X(-mK_X)|_D \lra 0.\]
By Kawamata-Viehweg vanishing, $h^i(X, -mK_X) = 0$ for $i > 0$. By Proposition~\ref{Wilson} $D\cdot l_{\psi} < 0$ for any irreducible curve $l_{\psi}$ contracted by $\psi$, i.e. $-mK_X -D$ is big and nef for $m \gg 0$. Therefore $h^i(X, -mK_X-D) = 0$ for $i > 0$. We obtain $h^i(X, -mK_X|_D) = 0$ for $i > 0$, hence
 \[h^0(B, m\mu^*H'|_B) = h^0(D, -mK_X|_D) = \chi(X, -mK_X) - \chi(X, -mK_X -D).\]
By the Riemann-Roch formula for threefolds, then
 \[h^0(B, m\mu^*H'|_B) = \frac{1}{12}c_2(X)\cdot D + \frac{1}{4}K_X\cdot D^2 + \frac{1}{6}D^3 + m \cdot \frac{1}{2}K_X\cdot D^2.\]
Since $h^1(B, m\mu^*H') = 0$ for $m \gg 0$, the linear term of the right hand side is $1-g_B$:
 \[g_B = 1-\frac{1}{12}c_2(X)\cdot D - \frac{1}{4}K_X\cdot D^2 - \frac{1}{6}D^3.\] 
To compute $c_2(X)\cdot \phi^*H$ and $c_2(X)\cdot E$ note 
 \[\chi(Y, H) = h^0(Y, H) = h^0(X, \phi^*H) = \chi(X, \phi^*H)\]
and $-K_X\cdot c_2(X) = c_1(X)\cdot c_2(X) = 24$. We get $c_2(X)\cdot \phi^*H = \frac{24}{r} + d$ and $c_2(X)\cdot E = rd$. The first claim follows. Pulling back the cycle $\mu^*H'|_B$ to $X$, we obtain $-K_X\cdot D^2 = D\cdot  \sum l_i = (\mu^*H'\cdot B)\cdot (D\cdot l_1)$ for general exceptional fibers $l_i$ of $\psi$, hence $\mu^*H'\cdot B = \frac{1}{2}K_X\cdot D^2$. 

Finally assume $\epsilon = 2$. Then $D$ is a $\PN_1$-bundle over $B$, hence smooth with $K_D^2 = 8(1-g_B)$. By adjunction, $K_D^2 = 2K_X\cdot D^2 + D^3$.
\end{proof}

\begin{theorem} \label{divcurve}
Let $X$ be an almost Fano threefold with $\rho(X) = 2$. Assume $X = Bl_C(Y)$ with $Y$ a smooth Fano threefold of index $r$ and $C \subset Y$ a smooth curve, and assume $|{-}K_X|$ induces a divisorial map $\psi: X \to X'$, contracting $D$ to a curve $B \subset X'$. Then we are in one of cases no.~1-24 in table~\ref{tabdivcurve} and all of them really exist.
\end{theorem}

\begin{lemma} \label{beta}
Assume $C$ is neither a line nor a conic. Then in the situation of the theorem
 \[D = \alpha \phi^*H - \beta E \quad \Longrightarrow \quad 1 \le \beta \le 4\]
and either
 \begin{enumerate}
  \item $\alpha = \beta r -1$, $\phi^*H\cdot l_{\psi} = 2$ or
  \item $\alpha = \beta r -2$, $\phi^*H\cdot l_{\psi} = 1$ and $D$ is a smooth ruled surface over $B$.
 \end{enumerate} 
More precisely: if $X'$ is not hyperelliptic, then $\beta \le 3$. If $X'$ is hyperelliptic, then one of the following holds
 \begin{enumerate}
  \item $({-}K_X)^3 = 2$ and either $\beta \le 3$ or $\beta = 4$ and $\epsilon = r$;
  \item $({-}K_X)^3 = 4$ and either $\beta \le 2$ or $\beta = 3$ and $\epsilon = r$.
 \end{enumerate} 
\end{lemma}
\noindent No.~22 in table~\ref{tabdivcurve} shows that the bound for $\beta$ in the Lemma is sharp in the non-hyperelliptic case.
 
\begin{proof}
1.) The case $X'$ not hyperelliptic and $g \le 5$. By Proposition~\ref{ci}, $X'$ is a complete intersection in $\PN_{g+1}$, namely, a quartic hypersurface in $\PN_4$, the intersection of a quadric and a cubic in $\PN_5$ or the intersection of three quadrics in $\PN_6$. In all of these cases, the Jacobian ideal is generated by cubics, defining some subscheme $\tilde{B}$ of $\PN_{g+1}$, such that $\tilde{B}_{red} = B$. Denote by $\I = \I_{\tilde{B}/\PN_{g+1}}$ the ideal sheaf of $\tilde{B}$. Then $\I(3)$ is globally generated, and if we define
 \[\J = \psi^{-1}\I \cdot \sO_X,\]
then also $\J \otimes \sO_X(-3K_X)$ is generated. Outside of some
codimension two subset of $X$, the sheaf $\J$ coincides with $\sO_X(-\lambda D)$ for some $\lambda \ge 1$. If we restrict $\sO_X(-\lambda D) \otimes \sO_X(-3K_X)$ to the general exceptional fiber $l_{\phi}$ of $\phi$, we still have sections, implying
 \[(-3K_X -\lambda D)\cdot l_{\phi} \ge 0.\]
This gives $\beta \le \frac{3}{\lambda} \le 3$.

2.) The case $X'$ not hyperelliptic and $g \ge 6$. By (\ref{ab}),
 \begin{equation} \label{kxeven}
   ({-}K_X)^3 = \frac{\epsilon (r^2H^3-d)}{\beta} \in 2\KZ.
 \end{equation}
From $g \ge 6$ we infer $({-}K_X)^3 \ge 10$. In order to see $\beta \le 3$ we have to prove $r^2H^3-d \le 19$. From Iskovskikh's classification we get $r^2H^3\le 22$, hence $r^2H^3 -d \le 21$. If $r^2H^3 -d \ge 20$, then $r = 1$, $H^3 = 22$ and $d \le 2$, implying $C$ is a line or a conic, since $H$ is very ample. 

3.) The case $X'$ hyperelliptic. Since $C$ is neither a line nor a conic, we have $({-}K_X)^3 \le 4$ by Corollary~\ref{hypfinite}. If $W = \PN_3$, then $\mu$ is ramified along a sextic $S \in |\sO_{\PN_3}(6)|$, which is singular along the image $\mu(B) \simeq B$. The reduced strict transform $\hat{S}$ of $S$ in $X$ gives a smooth section in $|{-}3K_X -D| = |(3r-\alpha)\phi^*H -(3-\beta)E|$. If $E$ is not a connected component of $\hat{S}$, then the intersection with $l_{\phi}$ yields $\beta \le 3$. If $E$ is a component, we claim $\hat{S} = E$. Assume $\hat{S} = E + \hat{S}'$ with some effective $\hat{S}' \in |\phi^*(r-\epsilon)H|$. This is impossible, since $\hat{S}$ is smooth, but any section of $|H|$ meets the curve $C$ on $Y$. We conclude $\beta = 4$ and $\epsilon = r$ in this case. If $W \subset \PN_4$ is a quadric, then $\mu$ is branched along a quartic. The same argument completes the proof of the lemma.
\end{proof}

\begin{proof}[Proof of \ref{divcurve}.]
By Corollary~\ref{gencor}, $-K_X$ is globally generated. We have already considered lines and conics in \ref{lines} and obtained no.~1,3, and 5. Assume therefore $C$ is neither a line nor a conic. The structure of the proof is as follows: First we show some numerical conditions in Step~1,2 and 3. By the lemma, $\beta \le 4$. We consider the cases $\beta = 1,2,3$ and $4$ in Step~4-7 separately. Finally, we will construct almost Fano threefolds corresponding to the data listed in table~\ref{tabdivcurve} in Step~8.

\vspace{0.2cm}

\noindent {\bf Step~1.} {\em Numerical formulas in the case $({-}K_X)^3 \le 8$ and $X'$ not hyperelliptic.} Then $X'$ is a complete intersection in $\PN_{g+1}$ by Proposition~\ref{ci} and hence $X$ is a complete intersection in the blowup 
 \[\widehat{\psi}: \widehat{\PN}_{g+1} \lra \PN_{g+1}\] 
of $\PN_{g+1}$ in the smooth curve $B$. Call the exceptional divisor $\widehat{D}$ and let again be $H'$ be a general hyperplane section of $\PN_{g+1}$. Then $\widehat{\psi}|_X = \psi$ and $\widehat{D} \cap X = D$. We have 
 \[\epsilon \phi^*H = \beta ({-}K_X) -D = \beta\psi^*H' -D.\] 
The following formula's on intersection numbers are well known: $\widehat{\psi}^*{H'}^{g+1} = 1$,  $\widehat{\psi}^*H'\cdot \widehat{D}^g = (-1)^{g+1}\deg(B)$, $\widehat{D}^{g+1} = (-1)^{g+1}\big(2g_B-2 + (g+2)\deg(B)\big)$ and $\widehat{\psi}^*{H'}^i\cdot \widehat{D}^{g+1-i} = 0$ for $i \not= 0, 1, g+1$. 

If $({-}K_X)^3 = 4$, then $X \in |4\widehat{\psi}^*H'-2\widehat{D}|$; if $({-}K_X)^3 = 6$, then $X$ is the complete intersection either of type $(3\widehat{\psi}^*H'-2\widehat{D}, 2\widehat{\psi}^*H'-\widehat{D})$, or $(3\widehat{\psi}^*H'-\widehat{D}, 2\widehat{\psi}^*H'-2\widehat{D})$; if $({-}K_X)^3 = 8$, then $X$ is of type $(2\widehat{\psi}^*H'-\widehat{D}, 2\widehat{\psi}^*H'-\widehat{D}, 2\widehat{\psi}^*H'-2\widehat{D})$. We obtain
\begin{eqnarray}
 ({-}K_X)^3 = 4: \quad \epsilon^3H^3 & = & 4\beta^3 -(6\beta -6)\deg(B) + 4(g_B -1); \label{n4}\\
 ({-}K_X)^3 = 6: \quad \epsilon^3H^3 & = & 6\beta^3 -(6\beta -5)\deg(B) + 4(g_B -1); \label{n5}\\
     \quad \mbox{ or } \quad \epsilon^3H^3 & = & 6\beta^3 -(6\beta -4)\deg(B) + 4(g_B -1); \nonumber\\
 ({-}K_X)^3 = 8: \quad \epsilon^3H^3 & = & 8\beta^3 -(6\beta -4)\deg(B) + 4(g_B -1). \label{n6}
\end{eqnarray}

\vspace{0.2cm} 

\noindent {\bf Step~2.} {\em Numerical formulas in the case $X'$ hyperelliptic.} Here $\mu: X' \to W$ is ramified along some Cartier divisor in $W$, and $W$ is either $\PN_3$ or a quadric in $\PN_4$. We obtain an induced double cover $\sigma: X \to Bl_{\mu(B)}(W)$. Since $D \subset X$ is exceptional, it is stabilized by the automorphism inducing $\sigma$, i.e. $\sigma^*D' = D$ with $D'$ the exceptional divisor for $\sigma$. 

As in the last step, write $\epsilon \phi^*H = \beta ({-}K_X)
-D$. Then $\epsilon^3H^3 = 2(\beta\sigma^*H' -D')^3$ yields again (\ref{n4}) for $({-}K_X)^3 = 4$ and
 \begin{equation} \label{h2}
   ({-}K_X)^3 = 2: \quad \epsilon^3H^3 = 2\beta^3 - (6\beta -8)\deg(B) + 4(g_B-1).
 \end{equation}

\vspace{0.2cm} 

\noindent {\bf Step~3.} {\em If $r = 1$, then $H^3 \ge 10$.} Indeed: by (\ref{absch}), $2g_C-2 < d$, hence $h^0(C, H|_C) = 1-g_C+d$ by the Riemann-Roch theorem. Again, by (\ref{absch}), we obtain
 \[h^0(C, H|_C) \le \frac{H^3}{2} -1.\]
On the other hand, $H$ is ample and globally generated by \cite{AG5}, Corollary~2.4.6, hence $h^0(C, H|_C) \ge 2$. This shows $H^3 \ge 6$. If $H^3 = 6$, then $h^0(C, H|_C) = 2$, i.e. $C$ is a line. If $H^3 = 8$, then $C$ is a plane curve, meaning $2g_C = (d-1)(d-2)$. Then
 \[2g_C-2 = d(d-3) > (2g_C-2)(2g_C-5)\]
by (\ref{absch}), hence $g_C = 0$, i.e. $C$ is a line or a conic. 

\vspace{0.2cm} 

\noindent {\bf Step~4.} {\em The case $\beta=1$.} Using (\ref{ab}) we obtain no.~6, 9, 10, 11, 15, 16, 17, 18, 20, 21, 22, 23 and 24.

\vspace{0.2cm}

\noindent {\bf Step~5.} {\em The case $\beta=2$.} {\bf (i)} $r = 1$, $\epsilon = 1$. Then $\alpha = 1$, $2({-}K_X)^3 = H^3 -d$, $K_X\cdot D^2 = H^3 -2d$, $4(g_C-1) = 3d-H^3$ and $4(g_B -1) = -3d+H^3$, hence $g_C = g_B = 0$. We obtain no.~4 for $H^3 = 16$. The case $H^3 = 22$ is impossible, since $X'$ cannot be hyperelliptic by Lemma~\ref{beta}, but (\ref{n6}) does not hold.

{\bf (ii)} $r = 2$, $\epsilon = 1$. Then $\alpha = 3$, $2({-}K_X)^3 = 4H^3 -d$, $K_X\cdot D^2 = 6H^3 -2d$, $4(g_C-1) = 7d -12H^3$ and $4g_B = 8H^3 -8 -3d$, hence $4\mid d$ and $d \le 8$. For $d = 8$ we obtain no.~8. For $d = 4$ we have $H^3 = 3$, $({-}K_X)^3 = 4$, $g_C = 0$, $g_B = 1$, $K_X\cdot D^2 = 10$, but (\ref{n4}) does not hold.

{\bf (iii)} $r = 2$, $\epsilon = 2$. Then $\alpha = 2$, $({-}K_X)^3 = 4H^3 -d$, $K_X\cdot D^2 = 8H^3 -4d$, $4(g_C-1) = 6d -8H^3$, $D^3 = -24H^3 + 16d$ and $2g_B = -2 + 4H^3 - 3d$, hence $d \le 6$ is even. By Lemma~\ref{B} (2), here $D^3 = 4d -24a$, which gives $2H^3 = 4+d$. This is no.~12.

{\bf (iv)} $r = 3$, $\epsilon = 1$. Then $\alpha = 5$, $H^3 = 2$, $2({-}K_X)^3 = 18 -d$, $K_X\cdot D^2 = 30 -2d$, $4(g_C-1) = 11d -90$ and $4g_B = 34 -3d$. Hence $d$ is even, $8 \le d \le 11$ and $4 \not|\; d$. We obtain no.~14.

{\bf (iv)} $r = 3$, $\epsilon = 2$. Then $\alpha = 4$, $H^3 = 2$ and hence $D^3 = 16(d-10)$. On the other hand, $D^3 = 4d -16\alpha$ by Lemma~\ref{B} (2), hence $d = 8$ and $g_C = 3$. Since $\epsilon = 2$, $D$ is a smooth $\PN_1$-bundle over $B$. The image $D_1 = \phi(D)$ is a nonnormal, Gorenstein surface of degree $4$ in $Q_3$ with normalization $\phi_D:D \to D_1$. Let $\tau: Q_3 \to \PN_3$ be a double cover, ramified along some quadric $Q$ and define $D_2 = \tau(D_1)$, which is again Gorenstein. Assume $\tau|_{D_1}$ is generically 2:1. Then $D_2$ is a quadric in $\PN_3$ containing the image $\tau(C)$. The pullback to $X$ defines a section in $|2\phi^*H -E|$, which is impossible. Hence $\tau|_{D_1}$ is 1:1 and $D_2$ is a quartic in $\PN_3$. The same argument shows that $C$ cannot be contained in the ramification divisor $Q$, i.e., $\tau|_{D_1}$ is an isomorphism at least outside some curve $N$, meeting $C$ at most in points. Then by subadjunction, $K_{D_1} = \tau^*K_{D_2} - \lambda N$ for some $\lambda \ge 0$ outside the finite set $C \cap N$. On the other hand, $K_{D_1} = H|_{D_1}$ and $K_{D_2} = \sO_{D_2}$ by adjunction formula. Since $N \backslash (N \cap C)$ is in the smooth locus of $D_1 \backslash (N \cap C)$, the pullback to $D$ is well defined and we find $\phi^*H|_D = -\lambda\phi_D^*N$ outside a finite set. This is impossible.  

{\bf (v)} $r = 4$, $\epsilon = 1$. Then $\alpha = 7$, $H^3 = 1$, $2({-}K_X)^3 = 16 -d$, $K_X\cdot D^2 = 28 -2d$, $4(g_C-1) = 15d -112$ and $4g_B = 32-3d$. We conclude $d = 8$. Then $({-}K_X)^3 = 4$, $K_X\cdot D^2 = 12$, $g_C = 3$ and $g_B = 2$, but (\ref{n4}) does not hold.

{\bf (vi)} $r = 4$, $\epsilon = 2$. Then $\alpha = 6$ and $H^3 = 1$ and hence $D^3 = 4(4d-42)$. On the other hand, $D^3 = 4d -12a$ by Lemma~\ref{B} (2), hence $d = 8$. We obtain $({-}K_X)^3 = 8$, $K_X\cdot D^2 = 16$, $g_C = 5$ and $g_B = 2$, but (\ref{n6}) does not hold.   

\vspace{0.2cm}

\noindent {\bf Step~6.} {\em The case $\beta=3$.} {\bf (i)} $r = 1$, $\epsilon = 1$. Then $\alpha = 2$, $3({-}K_X)^3 = H^3 -d$, $K_X\cdot D^2 = 2H^3 -3d$, $6(g_C-1) = 5d - 2H^3$ and $6g_B = -18 + 7H^3 -13d$. We get either $H^3 = 18$ or $H^3 = 10$. In the first case, we have $d = 6$, $({-}K_X)^3 = 4$, $K_X\cdot D^2 = 18$, $g_C = 0$, but (\ref{n4}) does not hold. For $H^3 = 10$, we obtain no.~2.

{\bf (ii)} $r = 1$, $\epsilon = 2$. Then $\alpha = 1$, $3({-}K_X)^3 = 2(H^3 -d)$, $K_X\cdot D^2 = 2(H^3 -3d)$, $6(g_C-1) = 4d - H^3$ and $2g_B = 14 - H^3 + 2d$. From Lemma~\ref{B} (2) we infer $4d = H^3 -6$, hence $4 \mid (H^3 -6)$ and $2\mid d$. Then $H^3 = 14$, $18$ or $22$. In the first case we get $({-}K_X)^3 = 8$, but (\ref{n6}) does not hold. Assume $H^3 \ge 18$. Then $\alpha = 1$ and $\epsilon = 2$ imply that the image $\phi(D) \subset Y$ is a member in $|H|$, cut out by lines and ${\rm Sing}(\phi(D)) = C$. We claim that the divisor $R$ cut out by all the lines is an irreducible, generically reduced member of $|dH|$ for $d \ge 2$, except $Y = V_{22}^{MU}$, the almost homogeneous Mukai-Umemura threefold. Indeed, if $Y \not= V_{22}^{MU}$, but $H^3 = 18$ or $22$, every component of $R$ is generically reduced by \cite{ProkhEx} and by \cite{IS}, $Y$ contains a $1$-dimensional family of pairs of intersecting lines, implying $R \in |\sO_Y(d)|$ with $d \ge 2$ as claimed (see \cite{Tyurin}). But for $Y = V_{22}^{MU}$, the divisor cut out by lines is singular along the closed orbit, which is a rational curve of degree $12$ (see \cite{MU} or \cite{almc}). 

{\bf (iii)} $r = 2$, $\epsilon = 1$. Then $\alpha = 5$, $3({-}K_X)^3 = 4H^3 -d$, $K_X\cdot D^2 = 10H^3 -3d$, $6(g_C-1) = 11d - 20H^3$ and $6g_B = -24 +40H^3 -13d$. We get $H^3 = 3$ or $H^3 = 4$. In the first case we have $({-}K_X)^3 = 2$, $d = 6$, $g_C = 2$, $K_X\cdot D^2 = 12$ and $g_B = 3$, but (\ref{h2}) does not hold. For $H^3 = 4$ we get no.~7.

{\bf (iv)} $r = 2$, $\epsilon = 2$. Then $\alpha = 4$, $3({-}K_X)^3 = 2(4H^3 -d)$, $K_X\cdot D^2 = 2(8H^3 -3d)$, $6(g_C-1) = 10d - 16H^3$ and $g_B = 13 -4H^3 +d$. From Lemma~\ref{B} (2) we infer $2d = 5H^3 -6$ hence $H^3$ is even. We get $H^3 = 4$, $({-}K_X)^3 = 6$, $d = 7$, $g_C = 2$, $K_X\cdot D^2 = 22$ and $g_B = 4$, but (\ref{n5}) does not hold.

{\bf (v)} $r = 3$, $\epsilon = 1$. Then $\alpha = 8$, $3({-}K_X)^3 = 18 -d$, $K_X\cdot D^2 = 48-3d$ and $6(g_C -1) = 17d-144$, hence $6\mid d$. We get no.~13.

{\bf (vi)} $r = 3$, $\epsilon = 2$. Then $\alpha = 7$, $3({-}K_X)^3 = 2(18 -d)$, $K_X\cdot D^2 = 2(42-3d)$, $6(g_C-1) = 16d- 126$. From Lemma~\ref{B} (2) we infer $2d = 21$, which is impossible.

{\bf (vii)} $r = 4$, $\epsilon = 1$. Then $\alpha = 11$, $3({-}K_X)^3 = 16 -d$ and $6(g_C-1) = 23d-176$. We get $({-}K_X)^3 = 2$, $d = 10$ and $g_C = 10$, but (\ref{h2}) does not hold.

{\bf (viii)} $r = 4$, $\epsilon = 2$. Then $\alpha = 10$ and from Lemma~\ref{B} (2) we infer $d = 10$, which is case no.~19.

\vspace{0.2cm}

\noindent {\bf Step~7.} {\em The case $\beta=4$.} Here $X'$ is hyperelliptic with $({-}K_X)^3 = 2$ and $\epsilon = r$.

{\bf (i)} $r = 1$, $\epsilon = 1$. Then $\alpha = 3$, $d = H^3 -8$, hence $K_X\cdot D^2 = 32 -H^3$, $2g_C = H^3 -12$, $2g_B = 58 -3H^3$. We obtain $12 \le H^3 \le 18$. From (\ref{h2}) we infer $H^3 = 16$, $d = 8$, $g_C = 2$, $K_X\cdot D^2 = 16$ and $g_B = 5$. The double cover $\mu: X' \to \PN_3$ is ramified along the image of $E$. Denote $E_1 = \psi(\mu(E))$. Then $E_1$ is a sextic in $\PN_3$, containing the smooth curve $\mu(B) \simeq B$ with multiplicity $2$. Define $\tau: \tilde{X} = Bl_{\mu(B)}(\PN_3) \to \PN_3$ and call the exceptional divisor $\tilde{E}$. We get an induced double cover $\sigma: X \to \tilde{X}$, ramified along the strict transform $E_2$ of $E_1$, and such that $\sigma^*\tilde{E} = D$. The pullback $\sigma^*({-}K_{\tilde{X}}) = \sigma^*(\tau^*\sO_{\PN_3}(1) - \tilde{E}) = -4K_X - D = \phi^*H$ is big and nef on $X$, hence $\tilde{X}$ is again an almost Fano threefold with $\rho = 2$. Since the anticanonical map of $\tilde{X}$ contracts the divisor $E_2$ to a curve isomorphic to $C$ by construction, $E_2$ plays the role of $D$ on $\tilde{X}$ and we obtain the invariants $\beta(\tilde{X}) = 2$ and $\alpha(\tilde{X}) = 6$. Such an almost threefold does not exist by Step~5.  

{\bf (ii)} $r = 2$, $\epsilon = 2$. Then $\alpha = 6$, $d = 4H^3 -4$, $K_X\cdot D^2 = 2(16-4H^3)$, $g_C = 4H^3 -6$, hence $2 \le H^3 \le 3$. From Lemma~\ref{B} (2) we infer $9H^3 = 12$, which is impossible.

%%%%%%%%%%%%%%%%%%%%%%

\vspace{0.2cm}

\noindent {\bf Step~8.} {\em Constructions.}

\vspace{0.2cm}

\noindent {\bf No.~1.} 
Let $\tilde{\phi}: \tilde{X} \to V_{2,4}$ be the blowup of the Fano threefold $V_{2,4}$ along a line $\tilde{C}$, with exceptional divisor $\tilde{E}$. Then $\tilde{X}$ is a Fano threefold with $\rho = 2$, such that the second contraction is birational, contracting a divisor $\tilde{D} \in |\tilde{\phi}^*\sO_{Q_3}(2) -3\tilde{E}|$ onto a curve of genus $2$ and degree $5$ in $\PN_3$ (see \cite{AG5}, \S~12.3., no.~19). Consider the double cover $Y = V_{1,8} \to V_{2,4}$, ramified along a general quadric $S$. Then $S$ meets $\tilde{C}$ in $2$ points, i.e., the pullback of $\tilde{C}$ yields a smooth rational curve $C \subset Y$ of degree $2$. The blowup $X = Bl_C(Y)$ is a double cover $\sigma: X \to \tilde{X}$, ramified along $\tilde{\phi}^*S$, i.e., $-K_X = \sigma^*(-K_{\tilde{X}} - \tilde{\phi}^*\sO_{V_{2,4}}(1))$ is trivial on the pullback of a general fiber of the divisor $\tilde{D}$. Then $X$ is an almost Fano threefold with $\rho(X) = 2$, such that $\psi$ contracts the pullback $D$ of $\tilde{D}$ to a curve of genus $2$ and degree $5$ in $\PN_3$.

Note that every Fano threefold $\tilde{X}$ with $\rho(\tilde{X}) = 2$ admitting two birational contractions, both contracting a divisor to a curve, provides an example for our case. By classification, there are $6$ such $\tilde{X}$. We obtain no.~1, 2, 3, 7, 8, and 13 this way.
 
\vspace{0.2cm}

\noindent {\bf No.~2.} 
Analogously to no.~1 consider a double cover $V_{1,10} \to V_{2,5}$ ramified along a general quadric $S$ in $V_{2,5}$ and let $C$ be a double cover of a conic $\tilde{C}$ in $V_{2,5}$ (see also \cite{Gushel}). The blowup $\tilde{X} = Bl_{\tilde{C}}(V_{2,5})$ is a Fano threefold, such that the second contraction is birational (see \cite{AG5}, \S~12.3., no.~22). 

\vspace{0.2cm}

\noindent {\bf No.~3.} 
For a general conic in $V_{1,10}$, the anticanonical map $\psi$ is small, but \cite{AG5}, p.86, Example is a special one with $\psi$ divisorial. We may also construct $X$ as double cover of a Fano threefold $\tilde{X}$ as in no.~1: take $\tilde{X} = Bl_{\tilde{C}}(V_{2,5})$, where $\tilde{C}$ is a line in $V_{2,5}$. Then $\tilde{X}$ is a Fano threefold, such that the second contraction is birational (see \cite{AG5}, \S~12.3., no.~26).

\vspace{0.2cm}

\noindent {\bf No.~4.} 
Let $x_0, \dots, x_5$ be homogeneous coordinates of $\PN_5$ and consider the rational normal curve $B = B_4$ of degree $4$ in $\PN_4 = \{x_5 = 0\}$ from Example~\ref{case22}, now embedded into $\PN_5$. Let $f \in \KC[x_0, \dots, x_4]$ be the cubic determinant, defining the secant variety $S_1(B)$ in $\PN_4$. Define a new cubic
 \[K = f + x_5q \in H^0(\PN_5, \I_B^2(3))\]
with $q \in \KC[x_0, \dots, x_5]$ a general quadric. Let 
 \[Q_4 = \sum\lambda_iq_i + x_5l\]
be a general quadric containing $B$, with $l \in \KC[x_0, \dots, x_5]$ linear. Then the strict transforms $\hat{K}$ and $\hat{Q}_4$ of $K$ and $Q_4$ in $Bl_B(\PN_5)$ are smooth and their intersection is a smooth almost Fano threefold $X$ with anticanonical model $X' = K \cap Q_4$. 

The quadrics vanishing on $B$ define a birational, small contraction $\phi$ of $Bl_B(\PN_5)$ contracting the strict transform $\hat{S}_1(B)$ of $S_1(B) = K \cap \{x_5 = 0\}$ onto a $\PN_2$, embedded via Veronese into $\PN_{11}$. Define $Y = \phi(X)$ in $\PN_{11}$. On the intersection $X = \hat{K} \cap \hat{Q}_4$, the map $\phi$ becomes divisorial, contracting $\hat{S}_1(B) \cap \hat{Q}_4$ onto a conic in $\PN_2$. This shows $X = Bl_C(Y)$ with $C$ a rational curve of degree $4$ in $Y$.

\vspace{0.2cm}

\noindent {\bf No.~5.} 
Let $S \in |H|$ be a general section on $Y = V_{2,2}$. Since $H$ is globally generated (\cite{AG5}, Theorem~2.3.1), $S$ is a smooth. By adjunction, $-K_S = H|_S$, hence $S$ is a del Pezzo surface of degree $2$ and we have a map $\pi: S \to \PN_2$, which is a blowup in $7$ general points. Let $l_1, \dots, l_7$ be the $(-1)$-curves. Then $-K_S = \pi^*\sO(3) - \sum l_i$. Take
 \[C \in |\pi^*\sO(5) - 2l_1 - \dots -2l_6 -l_7|\]
general. We claim that $C$ is a smooth curve: $C$ is the strict transform of a curve in $\PN_2$ with exactly $6$ double points containing the $7$th point in the smooth locus. Let $\pi': S' \to \PN_2$ be the blowup in the first $6$ points with exceptional curves $l_1', \dots, l_6'$. Then $S'$ is the cubic in $\PN_3$ and a linear system $|W|$ contains a smooth irreducible member, if $W^2 > 0$ and $W\cdot l \ge 0$ for each of the $27$ lines $l$ on $S'$. This is true for $|{\pi'}^*\sO(5) - 2l'_1 - \dots -2l'_6|$. By construction, $C$ is a conic.

Define now $X = Bl_C(Y)$ and let $D$ be the strict transform of $S$. We claim that $X$ is almost Fano, such that the anticanonical map contracts $D$ to a rational curve. From $C \subset S$ we get a short exact sequence 
 \[0 \to \I_{S/Y} \to \I_{C/Y} \to \I_{C/S} \to 0.\] 
The sequence remains exact on $H^0$-level after twisting by $\sO_Y(2H)$. By construction, $\I_{C/S}(2) = 2H|_S -C = \pi^*\sO(1) - l_7$. This shows: $-K_X$ is generated with $({-}K_X|_D)^2 = 0$, but $-K_X|_D \not\equiv 0$, hence the morphism defined by $|{-}K_X|$ contracts $D$ to a curve. Proposition~\ref{ci} implies $X'$ is a complete intersection of a quadric and a cubic in $\PN_5$.

\vspace{0.2cm}

\noindent {\bf No.~6.} 
Analogously to no.~5, let $S \in |H|$ be general on $Y = V_{2,3}$. Then $S$ is a del Pezzo surface of degree $3$, i.e. $\pi: S \to \PN_2$ is the blowup in $6$ general points. Here take
 \[C \in |\pi^*\sO(5) - 2l_1 - \dots -2l_5 -l_6|\]
general. Consider again the blowup $\pi': S' \to \PN_2$ in the first $5$ points and note that ${\pi'}^*\sO(3) - l_1' - \dots - l'_4 -2l_5'$ is globally generated by \cite{Ha}, V, Proposition~4.3 and that ${\pi'}^*\sO(2) - l_1' - \dots - l'_4$ is globally generated as well. Moreover, $\sO(1)$ has enough sections, i.e. Bertini shows the existence of a smooth $C$. 

By construction, $C$ is a smooth elliptic curve of degree $4$ and $X = Bl_C(Y)$ is almost Fano with $({-}K_X)^3 = 8$, such that $\psi$ contracts the strict transform $D$ of $S$ to a rational curve. By Proposition~\ref{ci}, $X' \subset \PN_6$ is a complete intersection of $3$ quadrics.

\vspace{0.2cm}

\noindent {\bf No.~7.} 
Analogously to no.~1, let $\tilde{\phi}: \tilde{X} \to Q_3$ be the blowup of a quadric $Q_3$ along an elliptic curve $\tilde{C}$ of degree $5$, with exceptional divisor $\tilde{E}$. Then $\tilde{X}$ is a Fano threefold with $\rho = 2$, such that the second contraction is birational, contracting a divisor $\tilde{D} \in |\tilde{\phi}^*\sO_{Q_3}(5) -3\tilde{E}|$ onto an elliptic curve of degree $5$ in $\PN_3$ (see \cite{AG5}, \S~12.3., no.~17). Consider the double cover $Y = V_{2,4} \to Q_3$, ramified along a general quadric $S$. Then $S$ meets $\tilde{C}$ in $10$ points, i.e., the pullback of $\tilde{C}$ yields a smooth curve $C \subset Y$ of genus $6$ and degree $10$. As seen in no.~1, the blowup $X = Bl_C(Y)$ is an almost Fano threefold with $\rho(X) = 2$, such that $\psi$ contracts the pullback $D$ of $\tilde{D}$ to an elliptic curve of degree $5$ in $\PN_3$.

\vspace{0.2cm}

\noindent {\bf No.~8.} 
Let $x_0, \dots, x_4$ be homogeneous coordinates of $\PN_4$. Let $B = B_4$ be the rational normal curve of degree $4$ in $\PN_4$ and let $f$ be the cubic determinant, defining the secant variety $S_1(B)$ to $B$ as in Example~\ref{case22}. Then $S_1(B)$ is a canonical Gorenstein Fano threefold with singularities of type $cA_1$ along $B$ and the blowup $Bl_B(S_1(B))$ is smooth. Let on the other hand 
 \[Q_3 = \sum\lambda_iq_i, \quad Z = \sum\lambda_{ij}q_iq_j\]
be a general quadric containing $B$ and a general member in $H^0(\PN_4, \I_{B/\PN_4}^2(4))$. Denote the blowup map $Bl_B(\PN_4) \to \PN_4$ by $\psi$ with exceptional divisor $D$. It is easily checked that the strict transforms $\hat{Q}_3$ and $\hat{Z}$ in $Bl_B(\PN_4)$ are smooth. 

The quadrics $q_0, \dots, q_5$ cutting out $B$ define a morphism
 \[\phi: Bl_B(\PN_4) \lra \PN_5\] 
onto the smooth quadric $Q = q_1q_4 - q_2q_3 - q_0q_5$, which is the blowup of $Q$ in a $\PN_2$, embedded
into $\PN_5$ via Veronese, with exceptional divisor $Bl_B(S_1(B)) =
\PN(N^*_{\PN_2/Q})$ (see Example~\ref{case22}). By construction, $\hat{Q}_3 \in |\phi^*\sO_{\PN_5}(1)|$ is the pullback of a hyperplane section, hence $\hat{Q}_3 \cap Bl_B(S_1(B))$ is a $\PN_1$-bundle over a conic in $\PN_2$. The same holds true for $\hat{Z} \in |\phi^*\sO_{\PN_5}(2)|$, meaning $\hat{Z} \cap Bl_B(S_1(B))$ is a $\PN_1$-bundle over a smooth curve of degree $4$ in $\PN_2$.

1.) Define $X' = Z$ a quartic in $\PN_4$. Then the strict transform $\hat{Z} = X$ is a smooth almost Fano threefold, mapped onto the complete intersection $Y$ of two quadrics in $\PN_5$ by the birational map $\phi|_X$. The exceptional divisor of $\phi$ on $X$ is $\hat{Z} \cap Bl_B(S_1(B))$, i.e. $X$ is the blowup of $Y$ in a curve of degree $8$ and genus $3$.     

2.) Define a double cover $\sigma: V \to Bl_B(\PN_4)$, ramified along $\hat{Z}$. Let $X$ be defined by the restriction to $\hat{Q}_3$. Then $X$ is a smooth almost Fano threefold and the Stein factorization of $\psi \circ \sigma$ gives 
 \[\xymatrix{X \ar[r]^{\sigma} \ar[d] & \hat{Q}_3 \ar[d]^{\psi}\\
             X' \ar[r]^{2:1} & Q_3,}\]
where the induced map $X \to X'$ is the anticanonical map of $X$, i.e. $X'$ is hyperelliptic. The Stein factorization of $\phi \circ \sigma$ gives a similar diagram, defining a double cover
 \[Y \stackrel{2:1}{\lra}  Q \cap H,\]
with $H = \phi(\hat{Q_3})$ a hyperplane in $\PN_5$. Then $Y$ is a smooth Fano threefold of index $2$ and degree $4$. On the exceptional divisor $S_1(B) \to \PN_2$, the intersection with $\hat{Q}_3$ gives a $\PN_1$-bundle, $\sigma$ is ramified over $8$ fibers, since the ramification divisor $\hat{Z} \in |\phi^*\sO_{\PN_5}(2)|$. This shows $X \to Y$ contracts a divisor onto a smooth curve $C$, which is a double cover of $\PN_1$, ramified over $8$ points. Hence $g_C = 3$ and the degree of $C$ in $Y$ is $8$.  

We may also construct $X$ as double cover as in no.~1: take $\tilde{X} = Bl_{\tilde{C}}(Q_3)$, where $\tilde{C}$ is a twisted quartic. Then $\tilde{X}$ is a Fano threefold admitting a second birational contraction onto a quadric $Q_3$ (see \cite{AG5}, \S~12.3., no.~21).

\vspace{0.2cm}

\noindent {\bf No.~9.}
Analogously to no.~5, let $S \in |H|$ be general on $Y = V_{2,4}$. Then $S$ is a del Pezzo surface of degree $4$, i.e. $\pi: S \to \PN_2$ is the blowup in $5$ general points. Here take
 \[C \in |\pi^*\sO(5) - 2l_1 - \dots - 2l_4 -l_5|\]
general. The same argument shows $C$ is a smooth curve of genus $2$ and degree $6$, and $X = Bl_C(Y)$ is almost Fano with $({-}K_X)^3 = 10$, such that $\psi$ contracts the strict transform $D$ of $S$ to a rational curve. Mukai's classification applies, since $\epsilon = 1$.

\vspace{0.2cm}

\noindent {\bf No.~10.}
Analogously to no.~5, let $S \in |H|$ be general on $Y = V_{2,5}$. Then $S$ is a del Pezzo surface of degree $5$, i.e. $\pi: S \to \PN_2$ is the blowup in $4$ general points. Here take
 \[C \in |\pi^*\sO(5) - 2l_1 - \dots - 2l_3 -l_4|\]
general. The same argument shows $C$ is a smooth curve of genus $3$ and degree $8$, and $X = Bl_C(Y)$ is almost Fano with $({-}K_X)^3 = 12$, such that $\psi$ contracts the strict transform $D$ of $S$ to a rational curve. Mukai's classification applies, since $\epsilon = 1$.

\vspace{0.2cm}

\noindent {\bf No.~11.}
Analogously to no.~5, let $S \in |H|$ be general on $Y = V_{2,5}$. Then $S$ is a del Pezzo surface of degree $5$, i.e. $\pi: S \to \PN_2$ is the blowup in $4$ general points. Here take
 \[C \in |\pi^*\sO(4) - 2l_1 - \dots - 2l_3|\]
general. The same argument shows $C$ is a smooth rational curve of degree $6$, and $X = Bl_C(Y)$ is almost Fano with $({-}K_X)^3 = 14$, such that $\psi$ contracts the strict transform $D$ of $S$ to a rational curve. Mukai's classification applies, since $\epsilon = 1$.

\vspace{0.2cm}

\noindent {\bf No.~12.}
The following description of $V_{2,5}$ as an $\Sl_2(\KC)$-almost homogeneous threefold is due to Mukai and Umemura (\cite{MU}): think of $\PN_6$ as $\PN(H^0(\PN_1, \sO_{\PN_1}(6)))$ and denote by $t_0, t_1$ homogeneous coordinates of $\PN_1$. Then
  \[V_{2,5} = \overline{\Sl_2(\KC)[t_0t_1(t_0^4+t_1^4)]},\]
i.e., $V_{2,5}$ is the closure of $[t_0t_1(t_0^4+t_1^4)] \in \PN(H^0(\PN_1, \sO_{\PN_1}(6)))$ under the natural $\Sl_2(\KC)$ action. There are two more orbits besides $\Sl_2(\KC)[t_0t_1(t_0^4+t_1^4)]$ in $V_{2,5}$:

i) the two-dimensional orbit $\Sl_2(\KC)[t_0t_1^5]$,

ii) the one-dimensional closed orbit $C = \Sl_2(\KC)[t_1^6]$.

The union of these two orbits gives an irreducible surface $S \subset V_{2,5}$, contained in $|2H|$, and singular along $C$. From the above description one sees that the normalization of $S$ is isomorphic to $\PN_1 \times \PN_1$, the normalization map
  \[\nu: \PN_1 \times \PN_1 \lra S\]
is given by an incomplete subsystem of bidegree $(1, 5)$, and $\nu(\Delta) = C$. In particular $C$ is a rational normal curve in $\PN_6$.

Define $X = Bl_C(V_{2,5})$. The strict transform $D \simeq \PN_1 \times \PN_1$ of $S$ is in $|2\phi^*H - 2E|$. The restriction of $-K_X$ to $D$ is of bidegree $(0,8)$. We claim that $-K_{X'}$ does not admit a moving decomposition: assume $-K_{X'} \sim W_1 + W_2$ for two Weil divisors. If $\widehat{W}_1$ and $\widehat{W}_2$ are their strict transforms in $X$, then $-K_X = \widehat{W}_1 + \widehat{W}_2 + \lambda D$ for some integer $\lambda \ge 0$. This means $(1-2\lambda)(\phi^*H -E)$ is effective, implying $\lambda = 0$. Then $\widehat{W}_1 = \phi^*H - kE$ and $\widehat{W}_2 = (k-1)E$ for some $k \ge 1$ (or vice versa). But this is impossible, since $D \in |2\phi^*H - 2E|$ is exceptional. 

\vspace{0.2cm}

\noindent {\bf No.~13.} 
Analogously to no.~1, here $X$ is a double cover of the Fano threefold $\tilde{X}$, obtained as blowup of $\PN_3$ along a curve $\tilde{C}$ of degree $6$ and genus $3$, which is cut out by cubics (\cite{AG5}, \S~12.3., no.~12). The second contraction of $\tilde{X}$ is birational, contracting a divisor $\tilde{D} \in |\tilde{\phi}^*\sO_{\PN_3}(8) - 3\tilde{E}|$ again to a curve of degree $6$ and genus $3$ in $\PN_3$. Let $Y = Q_3 \to \PN_3$ be a double cover, ramified along a general quadric and let $C$ be the pullback of $\tilde{C}$. Then $C$ is a smooth curve of genus $11$ and degree $12$ and the blowup $X = Bl_C(Y)$ is almost Fano, such that $\psi$ contracts the pullback $D$ of $\tilde{D}$ to a curve $B$ of degree $6$ and genus $3$ in $\PN_3$. 

\vspace{0.2cm}

\noindent {\bf No.~14.} 
Let $x_0, \dots, x_4$ be homogeneous coordinates of $\PN_4$. Let $B$ be the elliptic normal curve of degree $5$ in $\PN_4$, cut out by the $5$ quadrics $q_0, \dots, q_4$ vanishing on $B$. Let $\psi: Bl_B(\PN_4) \to \PN_4$ the blowup of $B$ and define a quartic $K$ in $\PN_4$ by $\sum\lambda_{ij}q_iq_j$ with $\lambda_{ij} \in \KC$  general.

1.) Define $X' = K$. Then the strict transform $X$ of $X'$ in $Bl_B(\PN_4)$ is a smooth almost Fano threefold. The map $Bl_B(\PN_4) \to \PN_4$ defined by $q_0, \dots, q_4$ realizes $X$ as blowup of the quadric $Q_3$ defined by $\sum\lambda_{ij}y_iy_j$.

2.) Define a smooth quadric $W \subset \PN_4$ containing $B$ by $\sum\lambda_iq_i$. Consider the double cover $\sigma: V \to Bl_B(\PN_4)$, ramified along the (smooth) strict transform of the quartic $K$. The induced covering of the strict transform of $W$ is a smooth almost Fano threefold $X$, and the Stein factorization of $\sigma \circ \psi$ is the anticanonical map $X \to X'$, followed by a double cover $X' \to W$, ramified along the intersection of the quadric and the quartic, which still contains $B$ as singular locus. 

\vspace{0.2cm}

\noindent {\bf No.~15.}
Analogously to no.~5, let $S \in |2H|$ be general on $Y = Q_3$. Then $S$ is a del Pezzo surface of degree $4$, i.e. $\pi: S \to \PN_2$ is the blowup in $5$ general points. Here take
 \[C \in |\pi^*\sO(8) - 3l_1 - \dots - 3l_4 - 2l_5|\]
general. Write $\pi^*\sO(8) - 3l_1  - \dots - 3l_4 - 2l_5  =  -2K_S + (\pi^*\sO(2) - l_1 - \dots - l_4)$. The existence of a smooth $C$ follows, since $4$ general points in $\PN_2$ are cut out by $2$ quadrics and $-K_S$ is very ample. Hence $C$ is a smooth curve of genus $8$ and degree $10$. By Proposition~\ref{ci}, $X' \subset \PN_6$ is a complete intersection of $3$ quadrics.

\vspace{0.2cm}

\noindent {\bf No.~16.}
Analogously to no.~15, let $S \in |2H|$ be general on $Y = Q_3$. Here take
 \[C \in |\pi^*\sO(7) - 3l_i - \dots - 3l_4 - l_5|\]
general. Write $\pi^*\sO(7) - 3l_1 - \dots -3l_4  - l_5  =  -K_S + 2(\pi^*\sO(2) - l_1 - \dots -l_4)$. The existence of a smooth $C$ follows as above. Hence $C$ is a smooth curve of genus $3$ and degree $8$. Since $\epsilon = 1$, Mukai's classification for $X'$ applies.

\vspace{0.2cm}

\noindent {\bf No.~17.}
Analogously to no.~15, let $S \in |H|$ be general on $Y = Q_3$. Then $S \simeq \PN_1 \times \PN_1$. Here take $C \in |\sO(2,3)|$ general. Then $C$ is a smooth curve of genus $2$ and degree $5$ on $Q_3$. Define as usual $X = Bl_C(Y)$ and let $D$ be the strict transform of $S$. From $3H|_S -C \in |\sO(1,0)|$ we infer that $-K_X$ is generated with $({-}K_X|_D)^2 = 0$, but $-K_X|_D \not\equiv 0$.

Note that $X'$ cannot be hyperelliptic by Proposition~\ref{hyp}. Analogously to no.~22, $X'$ is contained in the double cone over $Q_3 \hookrightarrow \PN_{13}$ embedded with $|\sO_{Q_3}(2)|$. But here $X \subset \PN(\sO_{Q_3}^{\oplus 2} \oplus \sO_{Q_3}(2))$ is not a complete intersection.

\vspace{0.2cm}

\noindent {\bf No.~18.}
Analogously to no.~17, let $S \in |H|$ be general on $Y = Q_3$. Then $S \simeq \PN_1 \times \PN_1$. Here take $C \in |\sO(1,3)|$ general. Then $C$ is a smooth rational curve of degree $4$ on $Q_3$. Analogously to the last case, $X'$ is contained in the triple cone over $Q_3 \hookrightarrow \PN_{13}$. Again $X \subset \PN(\sO_{Q_3}^{\oplus 3} \oplus \sO_{Q_3}(2))$ is not a complete intersection.

\vspace{0.2cm}

\noindent {\bf No.~19.}
Let $x_0, \dots, x_4$ be homogeneous coordinates of $\PN_4$ and let $M$ be a symmetric $(4\times 4)$-matrix with linear entries in $\KC[x_0, \dots, x_4]$. Let $X'_4 \subset \PN_4$ be defined by the determinant of $M$. By \cite{Catanese}, the $(3\times 3)$-minors of $M$ cut out a smooth curve $B$ of genus $6$ and degree $10$ in $\PN_4$. This is the singular locus of $X' = X'_4$. 

Let on the other hand $y_0, \dots, y_3$ be homogeneous coordinates of $\PN_3$. Multiplying $M$ with the vector $(y_0, \dots, y_3)^t$, we obtain $4$ sections of $|\sO(1,1)|$ on $\PN_3 \times \PN_4$. These cut out a smooth almost Fano threefold $X$, the desingularisation of $X'$. The projection onto $\PN_4$ maps $X$ birationally onto $X'$, the exceptional locus is a divisor $D$ over the singular locus $B$ of $X'$. 

The projection onto $\PN_3$ is birational as well. Write $M\cdot (y_0, \dots, y_3)^t = N \cdot (x_0, \dots, x_4)^t$, where $N$ is a $(4\times 5)$-matrix with linear entries in $\KC[y_0, \dots, y_3]$. By the Hilbert-Burch theorem the $(4\times 4)$-minors of $N$ cut out a smooth curve $C$ of genus $11$ and degree $10$ in $\PN_3$ and $X = Bl_C(\PN_3)$.     
By Lemma~\ref{beta}, $X'$ cannot by hyperelliptic. 

Note that if $M$ is not symmetric, the singular locus of $X'$ consists of $20$ points, i.e. the contraction $X \to X'$ is small.

\vspace{0.2cm}

\noindent {\bf No.~20.}
Analogously to no.~5, let $S \in |3H|$ be general on $Y = \PN_3$. Then
$S$ is a del Pezzo surface of degree $3$. Take
 \[C \in |\pi^*\sO(11) - 4l_1 - \dots - 4l_5 -3l_6|\]
general. Then a linear system $|W|$ contains a smooth irreducible member if $W^2 > 0$ and $W\cdot l \ge 0$ for each of the $27$ lines $l$ on $S$ (\cite{Ha}, V, Theorem~4.11.). This is true in our case. By construction, $C$ is a curve of genus $12$ and degree $10$. By Proposition~\ref{ci}, $X' \subset \PN_5$ is a complete intersection of a quadric and a cubic.

\vspace{0.2cm}

\noindent {\bf No.~21.}
Analogously to no.~20, let $S \in |3H|$ be general on $Y = \PN_3$. Here take
 \[C \in |\pi^*\sO(10) - 4l_1 - \dots -4l_5 -2l_6|\]
general. This gives a smooth curve of genus $5$ and degree $8$ on $\PN_3$. By Proposition~\ref{ci}, $X' \subset \PN_6$ is a complete intersection of $3$ quadrics.

\vspace{0.2cm}

\noindent {\bf No.~22.}
Analogously to no.~17, let $S \in |2H|$ be general on $Y = \PN_3$. Then $S \simeq \PN_1 \times \PN_1$. Take $C \in |\sO(3,4)|$ general. Then $C$ is a smooth curve of genus $6$ and degree $7$ on $\PN_3$. Define as usual $X = Bl_C(Y)$ and let $D$ be the strict transform of $S$. From $4H|_S -C \in |\sO(1,0)|$ we infer that $-K_X$ is generated with $({-}K_X|_D)^2 = 0$, but $-K_X|_D \not\equiv 0$.

Note that $X'$ cannot be hyperelliptic by Proposition~\ref{hyp}. The anticanonical map $\psi: X \to X'$ is the resolution of the rational map $\PN_3 \rightharpoonup \PN_{11}$ defined by all quartics of $\PN_3$ vanishing on $C$. These are $qq_0, \dots, qq_9$ and two further quartics, where $q$ defines $S$ and $q_0, \dots, q_9$ are all quadrics of $\PN_3$. This shows $X'$ is contained in the double cone over the $2$-uple embedding $\PN_3 \hookrightarrow \PN_9$. The surface $S$ is mapped to the vertex of the cone, which is the line $B$. Blowing up $B$, we find $X$ as complete intersection of two general elements from $|\sO(1) \otimes \pi^*\sO_{\PN_3}(1)|$ in $\PN(\sO_{\PN_3}^{\oplus 2} \oplus \sO_{\PN_3}(2))$. 

The projection onto $\PN_3$ is our birational map $\phi$. A section from $|\sO(1) \otimes \pi^*\sO_{\PN_3}(1)|$ corresponds to a section of the vector bundle $\sO_{\PN_3}(1)^{\oplus 2} \oplus \sO_{\PN_3}(3)$, is hence a vector $(l, l', f)$ with two linear forms $l, l'$ and a cubic $f$ on $\PN_3$. The exceptional locus of $\phi$ then lies over the curve in $\PN_3$ cut out by the minors of the matrix
 \[\left(\begin{array}{ccc} l_1 & l_1' & f_1\\
                            l_2 & l_2' & f_2\end{array}\right),\]
the rows corresponding two the two sections cutting out $X$.

\vspace{0.2cm}

\noindent {\bf No.~23.}
Analogously to no.~22, let $S \in |2H|$ be general on $Y = \PN_3$ and take $C \in |\sO(2,4)|$ general. Then $C$ is a smooth curve of genus $3$ and degree $6$ on $\PN_3$ (compared to the curve $\tilde{C} \subset \PN_3$ in no.~13, here $C$ is not cut out by cubics). 

Note that $X'$ cannot be hyperelliptic by Proposition~\ref{hyp}. Analogously to the last case, here $X'$ is contained in the tripel cone over $\PN_3 \hookrightarrow \PN_9$, since $\PN_3 \rightharpoonup \PN_{12}$ is defined by $qq_0, \dots, qq_9$ and three further quartics. The surface $S$ is mapped to a conic in the vertex, i.e. $X'$ is moreover contained in some quadric, not containing the cone. In this case, $X \subset \PN(\sO_{\PN_3}^{\oplus 3} \oplus \sO_{\PN_3}(2))$ is not a complete intersection.

\vspace{0.2cm}

\noindent {\bf No.~24.}
Analogously to no.~22, let $S \in |2H|$ be general on $Y = \PN_3$ and take $C \in |\sO(1,4)|$ general. Then $C$ is a smooth rational curve of degree $5$ on $\PN_3$. 

Note that $X'$ cannot be hyperelliptic by Proposition~\ref{hyp}. Analogously to the last case, $X'$ is contained in the $4$-cone over $\PN_3 \hookrightarrow \PN_9$, since $\PN_3 \rightharpoonup \PN_{12}$ is defined by $qq_0, \dots, qq_9$ and $4$ further quartics. The surface $S$ is mapped to a rational curve of degree $3$ in the vertex. Again, $X \subset \PN(\sO_{\PN_3}^{\oplus 4} \oplus \sO_{\PN_3}(2))$ is not a complete intersection.
\end{proof}

\

%%%%%%%%%%%%%%%%%%%%

\noindent {\bf (B) Assume $\psi$ contracts $D$ to a point.} Let again be $D = \alpha\phi^*H - \beta E$ in $\Pic(X)$. Since $E \not= D$ and $D$ is irreducible, we have $\alpha, \beta > 0$ for the same reason as above. Concerning the anticanonical model $X'$ note that $|{-}K_{X'}|$ never admits a moving decomposition in this case (see \ref{Mukaiclass}): assume $-K_{X'} \sim W_1+W_2$ with Weil divisors $W_1$ and $W_2$. Since $|{-}K_{X'}|$ is base point free, we may assume a general member of $|W_1|$ does not meet the point $X'_{sing} = \psi(D)$. Then $W_1$ is Cartier, implying $W_2$ is Cartier as well. But $\Pic(X') \simeq \KZ \cdot ({-}K_{X'})$ by \ref{pic}. We conclude that either $X'$ is hyperelliptic or can be found in Mukai's table~\ref{Mukaitab}.  

\begin{theorem} \label{divpoint}
Let $X$ be an almost Fano threefold with $\rho(X) = 2$. Assume $X = Bl_C(Y)$ with $Y$ a smooth Fano threefold of index $r$ and $C \subset Y$ a smooth curve, and assume $|{-}K_X|$ induces a divisorial map $\psi: X \to X'$, contracting $D$ to a point. Then $\beta = 1$ and $C$ is the complete intersection of $\phi(D) \in |r'H|$ for some $r' < r$ and an element in $|{-}K_Y|$. We obtain the data of no.~25 in table~\ref{divcurve} and all of these cases really exist. 
\end{theorem}

\begin{proof}
 By assumption, $X'$ is smooth outside a single point $p = \psi(D)$, and $-K_{X'} = \mu^*H'$ is globally generated and ample. Assume $X'$ is not hyperelliptic. Then $H'$ is very ample on $X'$, i.e. the point $p$ may be cut out by the sections in $|{-}K_{X'}| = |H'|_{X'}|$ vanishing on $p$. Then the same argument as in the proof of Lemma~\ref{beta} shows $\beta = 1$ and $r-\alpha > 0$ by restricting to some exceptional fiber of $\phi$, or to some general irreducible curve in $X$ not meeting $E$. In particular, $r \ge 2$ and $\alpha \le 3$ in this case. If $X'$ is hyperelliptic, then $-2K_{X'}$ is very ample, and the same argument as above shows $\beta \le 2$ and $\alpha < 2r$. 

Assume first $\beta = 2$. Then $X'$ is hyperelliptic, hence $({-}K_X)^3 \le 8$ by Proposition~\ref{hyp}. Considering the twisted ideal sequence of $p$ in $X'$ shows
 \[h^0(X, -K_X-D) \ge h^0(X, -K_X) -1.\]
We have $-K_X - D = (r-\alpha)\phi^*H + E$, implying $\alpha < r$. In particular, $r \ge 2$. As above we have  
 \begin{eqnarray*}
   2,4,8 & = & r^3H^3 -2rd + 2g_C-2,\\
   0 & = & \alpha r^2H^3-(\alpha +2r)d + 4g_C-4,\\
   0 & = & \alpha^2rH^3 -4\alpha d + 8g_C-8.
 \end{eqnarray*}
If $r = 2$, then $\alpha = 1$ and the second line implies $4\mid d$, hence $H^3 = 4$ by the third line. We obtain $d = 4$ and $g_C = 2$, contradicting the first line. Assume $r = 3$. If $\alpha = 1$, the last two lines come up with $10d = 32$, which is impossible. For $\alpha = 2$ we obtain $d = 5$ and $g_C = 3$, contradicting the first line. If $r = 4$, the last two lines give $\alpha = 2$, $g_C = 3$ and $d = 4$, contradicting the first line.  

For $\beta = 1$, let $S \in |{-}K_X|$ be general. Since $X'$ has only an isolated singularity and $|{-}K_{X'}|$ is base point free, $S$ will not meet $D$. This means the images $\phi(S)$ and $\phi(D)$ meet transversally in $C$, i.e., $C$ is a complete intersection. The following construction completes the proof:

\vspace{0.2cm}

\noindent {\bf No.~25.}
Analogously to no.~5 above, let $S \in |r'H|$ be general on $Y$ for $r' = \alpha < r$. Then $S$ is a del Pezzo surface. We choose $C \in |rH|_S|$ general. Since $rH = -K_Y$ is globally generated, the existence of a smooth $C$ is clear. Then $d = r'H|_S\cdot (rH)|_S = rr' H^3$ and $2g_C-2 = rr'^2H^3$ by adjunction. We have $rH|_S -C \sim \sO_S$ on $S$, hence $-K_X$ is generated and $\psi$ contracts the strict transform of $S$ to a point. If $X'$ is not hyperelliptic, then Mukai's classification applies as seen above.
\end{proof}

\

%%%%%%%%%%%%%%%%%%%

\section{Blowdown to a point} \label{blowdownpoint}
\setcounter{lemma}{0}

\begin{abs} 
{\bf Setup.} In this section $\phi: X \to Y$ is the blowdown to a Fano threefold $Y$ with $\rho(Y) = 1$, such that the exceptional divisor $E$ is mapped to a point $p$. Then $Y$ is either smooth or has a terminal singularity at $p$. As usual, $X$ is a smooth almost Fano threefold with $\rho(X) = 2$, such that the anticanonical map $\psi$ is divisorial with exceptional divisor $D$. We have the following possibilities for $(Y, E)$ (see \cite{Mori}):
 \begin{enumerate}
  \item $(E,\sO_E(-E)) = (\PN_2, \sO_{\PN_2}(1))$ and $Y$ is smooth;
  \item $(E,\sO_E(-E)) = (\PN_2, \sO_{\PN_2}(2))$ and $Y$ is singular, $2$-Gorenstein;
  \item $(E,\sO_E(-E)) = (Q, \sO(1))$ and $Y$ is singular, Gorenstein.
 \end{enumerate}
Here $(Q, \sO(1))$ denotes an irreducible quadric in $\PN_3$ with restriction of $\sO_{\PN_3}(1)$, i.e., $Q$ is either $\PN_1 \times \PN_1$ or the quadric cone. We show the following classification result
\end{abs}

\begin{theorem} 
Let $X$ be an almost Fano threefold with $\rho(X) = 2$. Assume $X = Bl_p(Y)$ with $Y$ a smooth or terminal Fano threefold of index $r$ and $p \in Y$ a point. Assume $|{-}K_X|$ induces a divisorial map $\psi: X \to X'$. Then $\psi$ contracts an irreducible divisor $D$ to a smooth curve $B \subset X'$ and we are in one of the cases in table~\ref{tabpointpoint} and all of these cases really exist.
\end{theorem}

\begin{proof}
We first show that $\psi(D) = B$ is a smooth curve in $X'$. Let $l_{\psi}$ be any irreducible curve in the exceptional divisor $D$ of $\psi$. Since $-K_X$ is trivial on $l_{\psi}$, whereas $\phi$ is an extremal contraction and $Y$ is Fano, $l_{\psi}$ must meet $E$ in points. This shows $E \cap D \not= \emptyset$ and $\psi(D)$ is a curve $B$ in $X'$. By Proposition~\ref{Wilson}, $B$ is smooth and $D\cdot l_{\psi} = -2$ for the general exceptional fiber of $D$. As in \ref{B} we see that $B$ in $X'$ is a curve of degree 
 \[\deg(B) = \frac{K_X\cdot D^2}{2}.\]

\noindent {\bf Step~1.} {\em The case $(E,\sO_E(-E)) = (\PN_2, \sO_{\PN_2}(1))$.} Here $Y$ is a smooth Fano threefold with $\rho(Y) = 1$ and we may use Iskovskikh's classification for $Y$. Let $-K_Y = rH$ for some $r \in \KN$ and $H$ a generator of $\Pic(Y)$. Then 
 \[-K_X = r\phi^*H - 2E.\]
Since the blowup of $\PN_3$ and $Q_3$ in a point is a Fano variety, we have $r = 1$ or $2$. The Picard group of $X$ is generated by $\phi^*H$ and $E$ over $\KZ$, i.e. we may write
 \[D = \alpha\phi^*H - \beta E, \quad \mbox{ for some } \; \alpha, \beta \in \KN.\]
We have $\phi^*H^2\cdot E = \phi^*H\cdot E^2 = 0$, $E^3 = 1$. Therefore
 \[({-}K_X)^3 = r^3H^3 - 8.\]
From $K_X^2\cdot D = 0$ follows
 \[\beta = \frac{\alpha r^2}{4}H^3.\]
We have $r\phi^*H\cdot l_{\psi} = 2E\cdot l_{\psi}$ and $\alpha\phi^*H\cdot l_{\psi} - \beta E\cdot l_{\psi} = -2$ for the general exceptional fiber $l_{\psi}$, hence 
 \[4 = (\beta r -2\alpha)\phi^*H\cdot l_{\psi}.\]
Defining $\epsilon = \beta r -2\alpha$ as usual, we obtain  
 \[\epsilon = \frac{\alpha}{4}({-}K_X)^3 = 1,2,4 \quad \mbox{ and } \quad K_X\cdot D^2 = \frac{\alpha^2 rH^3}{8}({-}K_X)^3.\]
Since $Y$ is smooth, we may compute $\chi(Y, H)$ by Riemann-Roch and obtain by the same method as in Lemma~\ref{B}
 \[g_B = 1- \frac{2\alpha}{r} + \frac{\alpha^2 H^3(3r-2\alpha)}{12} + \frac{\beta^2(\beta -3)}{6}.\]
We continue case by case.

1.) $r = 1$. Then $-K_X$ is not divisible in $\Pic(X)$. If $({-}K_X)^3 \le 8$, then $X'$ is either hyperelliptic or a complete intersection by Proposition~\ref{ci}. Applying the same argument as in the proof of Lemma~\ref{beta}, we obtain in this case: if $X'$ is a complete intersection, then $\beta \le 6$. If $X'$ is hyperelliptic and $W = \PN_3$, then either $\beta \le 6$ or $\beta = 7$ and $r \le \epsilon$; if $W \subset \PN_4$ is a quadric, then either $\beta \le 4$ or $\beta = 5$ and $r \le \epsilon$; if $W$ is the Veronese cone, then either $\beta \le 3$ or $\beta = 4$, $r \le \epsilon$ and $2\mid (\epsilon -r)$. Moreover, we have the following formulas, analogously to (\ref{n4}) etc. 
 \begin{eqnarray} \label{p2}
  \hspace{0.8cm} ({-}K_X)^3 = 2: & \hspace{-0.3cm}& \epsilon^3H^3 = 2\beta^3 - 8(3\beta -8)\deg(B) + 32(g_B-1);\\\label{p4}
  \hspace{0.8cm} ({-}K_X)^3 = 4: & \hspace{-0.3cm}& \epsilon^3H^3 = 4\beta^3 - 24(\beta -2)\deg(B) + 32(g_B-1);\\\label{p6}
  \hspace{0.8cm} X' \mbox{ c.i.}, ({-}K_X)^3 = 8: & \hspace{-0.3cm}& \epsilon^3H^3 = 8\beta^3 - 8(3\beta -4)\deg(B) + 32(g_B-1).
 \end{eqnarray}

1.1.) $\epsilon = 1$. Then $\alpha(H^3 -8) = 4$ and $H^3 -8 = ({-}K_X)^3$ is even: (i) $\alpha = 1$ and $({-}K_X)^3 = 4$. Then $H^3 = 12$, $\beta = 3$, $K_X\cdot D^2 = 6$ and $g_B = 0$, contradicting (\ref{p4}); (ii) $\alpha = 2$ and $({-}K_X)^3 = 2$. Then $H^3= 10$, $\beta = 5$, $K_X\cdot D^2 = 10$ and $g_B = 2$, contradicting (\ref{p2}).

1.2.) $\epsilon = 2$. Then $\alpha(H^3 -8) = 8$: (i) $\alpha = 1$ and $({-}K_X)^3 = 8$. Then $H^3 = 16$, $\beta = 4$, $K_X\cdot D^2 = 16$ and $g_B = 3$. Since $\epsilon -r$ is not even, $X'$ is not hyperelliptic and we obtain a contradiction to (\ref{p6}); (ii) $\alpha = 2$ and $({-}K_X)^3 = 4$. Then $H^3 = 12$, $\beta = 6$, $K_X\cdot D^2 = 24$ and $g_B = 11$, contradicting (\ref{p4}); (iii) $\alpha = 4$ and $({-}K_X)^3 = 2$. We have $H^3 = 10$ and $\beta = 10$. This is impossible.

1.3.) $\epsilon = 4$. Then $\alpha(H^3 -8) = 16$: (i) $\alpha = 1$ and $({-}K_X)^3 = 16$. Then $H^3 = 24$, which is impossible; (ii) $\alpha = 2$ and $({-}K_X)^3 = 8$. Then $H^3 = 16$ and $\beta = 8$, a contradiction; (iii) $\alpha = 4$ and $({-}K_X)^3 = 4$. We have $H^3 = 12$ and $\beta = 12$, which is again impossible; (iv) $\alpha = 8$ and $({-}K_X)^3 = 2$. We have $H^3 = 10$ and $\beta = 20$, a contradiction.

2.) $r = 2$. Then $\epsilon = 2\alpha (H^3 -1)$, hence $\epsilon$ is even and $H^3 \ge 2$. 

2.1.) $\epsilon = 2$. Then $\alpha = 1$, $H^3 = 2$, $({-}K_X)^3 = 8$, $\beta = 2$, $K_X\cdot D^2 = 4$ and $g_B = 0$. This is no.~5.

2.2.) $\epsilon = 4$. (i) $\alpha = 1$ and $H^3 = 3$. Then $({-}K_X)^3 = 16$, $\beta = 3$, $K_X\cdot D^2 = 12$ and $g_B = 1$. This no.~7. (ii) $H^3 = 2$ and $\alpha = 2$. Then $({-}K_X)^3 = 8$, $\beta = 4$, $K_X\cdot D^2 = 16$ and $g_B = 3$. This is no.~6.

\vspace{0.2cm}

\noindent {\bf Step~2.} {\em The case $(E,\sO_E(-E)) = (\PN_2, \sO_{\PN_2}(2))$.} Now $Y$ is a terminal $2$-Gorenstein, $2$-factorial Fano threefold. Let $-2K_Y = rH$ for some $r \in \KN$ and $H$ a generator of $\Pic(Y)$. Then 
 \[-2K_X = r\phi^*H - E.\]
The Picard group of $X$ is generated by $\phi^*H$ and $E$ over $\frac{1}{2}\KZ$, i.e. we may write
 \[D = \alpha\phi^*H - \beta E, \quad \mbox{ for some } \; \alpha, \beta \in \frac{1}{2}\KN.\]
We have $\phi^*H^2\cdot E = \phi^*H\cdot E^2 = 0$, $E^3 = 4$. Therefore 
 \[({-}K_X)^3 = \frac{r^3H^3 - 4}{8}.\] 
From $K_X^2\cdot D = 0$ we deduce $4\beta = \alpha r^2 H^3$. For the general exceptional fiber $l_{\psi}$ we have $r\phi^*H\cdot l_{\psi} = E\cdot l_{\psi}$ and $\alpha\phi^*H\cdot l_{\psi} - \beta E\cdot l_{\psi} = -2$, hence $\epsilon := 2(\beta r - \alpha) = 1,2,4$. This gives
 \[4\alpha ({-}K_X)^3 = 1,2,4.\]
Since $({-}K_X)^3$ is even, we have $({-}K_X)^3 = 2$, $\alpha = \frac{1}{2}$, $\beta = \frac{5}{2}$, $r = 1$, $H^3 = 20$ and $K_X\cdot D^2 = 10$. Since $\phi^*H\cdot l_{\psi} = 1$, we have $K_D^2 = 8(1-g_B)$ and obtain $g_B = 6$; $X'$ is hyperelliptic with $W = \PN_3$ and the strict transform $\hat{S} \in |\phi^*H + E|$ is disconnected. This is no.~4.

\vspace{0.2cm}

\noindent {\bf Step~3.} {\em The case $(E,\sO_E(-E)) = (Q, \sO(1))$, with $Q$ a quadric in $\PN_3$.} Here $Q$ may be either $\PN_1 \times \PN_1$, or the quadric cone, and $Y$ is a terminal Gorenstein, but factorial Fano threefold. Let as usual be $-K_Y = rH$ for some $r \in \KN$ and $H$ a generator of $\Pic(Y)$. Then 
 \[-K_X = r\phi^*H - E.\]
The Picard group of $X$ is generated by $\phi^*H$ and $E$ over $\KZ$, i.e. we may write
 \[D = \alpha\phi^*H - \beta E, \quad \mbox{ for some } \; \alpha, \beta \in \KN.\]
We have $\phi^*H^2\cdot E = \phi^*H\cdot E^2 = 0$, $E^3 = 2$. Therefore
 \[({-}K_X)^3 = r^3H^3 - 2.\]
From $K_X^2\cdot D = 0$ follows
 \[\beta = \frac{\alpha r^2}{2}H^3,\]
and $r\phi^*H\cdot l_{\psi} = E\cdot l_{\psi}$ and $\alpha\phi^*H\cdot l_{\psi} - \beta E\cdot l_{\psi} = -2$ for the general exceptional fiber $l_{\psi}$ implies
 \[2 = (\beta r - \alpha) \phi^*H\cdot l_{\psi}.\]
Defining $\epsilon := \beta r - \alpha$ as usual, we get
 \[\epsilon = \frac{\alpha}{2}({-}K_X)^3 = 1,2 \quad \mbox{ and } \quad K_X\cdot D^2 = \frac{\alpha^2 rH^3}{2}({-}K_X)^3.\]

1.) $\epsilon = 1$. Then $2 = \alpha ({-}K_X)^3$ implies $\alpha = 1$
and $({-}K_X)^3 = 2$, hence $r = 1$, $H^3 = 4$, $\beta = 2$ and $K_X\cdot D^2 = 4$. Then $X' \to W = \PN_3$ is hyperelliptic and $\deg(\mu(B)) = 2$ implies $g_B = 0$. Here $Y$ is a Gorenstein Fano threefold with canonical (even terminal) singularities, hence either $Y \subset \PN_4$ is a quartic or $Y$ is a double cover of a quadric in $\PN_4$. In the first case, let $\pi: Bl_p(\PN_4) \to \PN_4$ be the blowup with exceptional divisor $\hat{E}$. Then $X \in |\pi^*\sO_{\PN_4}(4) -2\hat{E}|$ and $D \in |(\pi^*\sO_{\PN_4}(1) -2\hat{E})|_X|$. The twisted ideal sequence of $X$ in $Bl_p(\PN_4)$ shows that this system is empty. Hence $Y$ is hyperelliptic. This is no.~1.    

2.) $\epsilon = 2$. Then $4 = \alpha ({-}K_X)^3$: (i) $\alpha = ({-}K_X)^3 = 2$. Then $r = 1$, $H^3 = 4$, $K_X\cdot D^2 = 16$ and $\beta = 4$. We have $\phi^*H\cdot l_{\psi} = 1$, hence $D$ is a $\PN_1$-bundle over $B$ and $K_D^2 = 8(1-g_B)$. This gives $g_B = 9$. Here $X'$ is hyperelliptic with $W = \PN_3$ and the strict transform of the branch divisor $\hat{S} \in |\phi^*H + E|$ is disconnected. As in the last case we find that $Y$ is hyperelliptic. This no.~2. (ii) $\alpha = 1$ and $({-}K_X)^3 = 4$. Then $r = 1$, $H^3 = 6$, $\beta = 3$ and $K_X\cdot D^2 = 12$. We have $\phi^*H\cdot l_{\psi} = 1$, hence $D$ is a $\PN_1$-bundle over $B$ and $K_D^2 = 8(1-g_B)$. We obtain $g_B = 4$. This is no.~3. 

\vspace{0.2cm}

\noindent {\bf Step~4.} {\em Constructions.}

\vspace{0.2cm}

\noindent {\bf No.~1.}
Let $x_0, \dots, x_3$ be homogeneous coordinates of $\PN_3$ and let $q \in \KC[x_0, x_1, x_2]$ be homogeneous of degree $2$, defining a smooth conic $B'$ in $H' = \{x_3 = 0\} \simeq \PN_2$. Let $\sigma: \tilde{X} \to \PN_3$ be the blowup of $\PN_3$ along $B'$ with exceptional divisor $\tilde{D}$. Take a sextic defined by 
 \[S = \lambda_1q^3 + \lambda_2q^2x_3^2 + \lambda_3qx_3^4 + \lambda_4x_3^6,\] 
the $\lambda_i \in \KC$ general, and define $X$ as double cover $\tau: X \to \tilde{X}$, ramified along the strict transform $\tilde{S}$ of $S$. By construction, $-K_X = \tau^*\sigma^*\sO_{\PN_3}(1)$ is big and nef, hence $X$ is a smooth almost Fano threefold with $({-}K_X)^3 = 2$. Let $D \subset X$ be the induced cover of $\tilde{D}$. Then the Stein factorization of $X \to \PN_3$ factorizes over the anticanonical model $X'$ and $\psi: X \to X'$ contracts the divisor $D$ to a smooth curve $B \simeq B'$.  

Since $B'$ is cut out by quadrics, the system $|\sigma^*\sO_{\PN_3}(2) - \tilde{D}|$ is base point free and big, the pullback to $X$ defines a birational map $\phi: X \to Y$ onto some Fano threefold $Y$. The strict transform $\tilde{H'}$ of $H'$ in $\tilde{X}$ is isomorphic to $H' \simeq \PN_2$, and the restriction of $\tilde{S}$ to $\tilde{H'}$ is the conic $B'$. The induced cover of $\tilde{H}$ hence yields a smooth quadric $E \in |{-}K_X - D|$. We find that $\tau^*(\sigma^*\sO_{\PN_3}(2) - \tilde{D})$ is trivial on $E$, i.e., $\phi$ contracts $E$ to a point.  

\vspace{0.2cm}

\noindent {\bf No.~2.}
Let $x_0, \dots, x_3$ be homogeneous coordinates of $\PN_3$ and let $B'$ be the complete intersection of a smooth quadric $Q_2$ and a smooth quartic $K$. Then $B'$ is a curve of degree $8$ and genus $9$. Let $\sigma: \tilde{X} \to \PN_3$ be the blowup of $\PN_3$ along $B'$ with exceptional divisor $\tilde{D}$ and define $X$ as double cover $\tau: X \to \tilde{X}$, ramified along the disjoint union of the strict transforms $\tilde{Q_2}$ and $\tilde{K}$ of $Q_2$ and $K$. As in no.~1, $X$ is almost Fano and $\psi$ maps the pullback $D$ of $\tilde{D}$ to a curve $B \simeq B'$. 

The pullback of  $|\sigma^*\sO_{\PN_3}(4) - \tilde{D}|$ is base point
free and big, defining a birational map $\phi: X \to Y$. Since
$\tilde{K}$ is part of the ramification divisor, the sytem becomes
divisible on $X$. Let $E \subset X$ be the reduced pullback of
$\tilde{Q_2}$. Since $Q_2$ and $K$ intersect in $B'$, the restriction
of $\tilde{K} - \tilde{D}$ is trivial on $\tilde{Q_2}$, i.e., $\phi$ contracts $E$ to a point.    

\vspace{0.2cm}

\noindent {\bf No.~3.}
Let $B \subset \PN_4$ be the complete intersection of a linear subspace $L$, a smooth quadric $Q$ and a general cubic $K$. Then $B$ is a smooth curve of degree $6$ and genus $4$. As a curve in $E' = L \cap Q$, the curve $B$ is a divisor of degree $3$. Note that $E'$ is either $\PN_1 \times \PN_1$ or a quadric cone, depending on $L$ and $Q$. 

1.) Define $X'$ by the equation
 \[f = \lambda_1lk + \lambda_2q^2,\]
where $l,q,k$ are defining equations of $L, Q, K$, respectively, and $\lambda_1, \lambda_2 \in \KC$. Define $X = Bl_B(X')$. We find that $X$ is a smooth almost Fano threefold with $({-}K_X)^3 = 4$. 

Consider the rational map $\PN_4 \rightharpoonup \PN_{19}$ defined by the cubics vanishing on $B$. The resolution of this map yields a morphism $\phi': X \to \PN_{19}$, which is birational. By construction, $\phi'$ contracts the strict transform $E$ of $E'$ to a point. We have constructed the morphism associated to $|2\phi^*H| = |-3K_X -D|$. 

2.) Consider $B \subset Q$ and define $X'$ as double cover of $Q$, ramified along $R = L + K$. Note that the strict transforms of $L$ and $K$ become disjoint after blowing up $B$. 

\vspace{0.2cm}

\noindent {\bf No.~4.}
Let $K \subset \PN_3$ be a general smooth quintic and $H'$ a general hyperplane. Then $B' = K \cap H'$ is a smooth complete intersection curve of genus $6$, which is the singular locus of the sextic $S = K + H'$. In the blowup $Bl_{B'}(\PN_3)$, the strict transform $\hat{S}$ is the disjoint union of the smooth strict tramsforms of $K$ and $H'$. Define $X$ as double cover of $Bl_{B'}(\PN_3)$, ramified along $\hat{S}$.

\vspace{0.2cm}

\noindent {\bf No.~5.}
Consider the double cover $\nu: Y = V_{2,2} \to \PN_3$, ramified along a smooth quartic $K$ in $\PN_3$. Let $R \subset Y$ be the reduced pullback of $K$ and $p \in R$ a point. Note that $K \simeq R$ is smooth. Define $\phi: X = Bl_p(Y) \to Y$ the blowup in $p$ with exceptional divisor $E$.

Define as usual $H = \nu^*\sO_{\PN_3}(1)$, i.e., $-K_Y = 2H$. In order to prove that $X$ is almost Fano, it is sufficient to show that $|\phi^*H -E|$ is base point free. We have
 \[0 \lra \I_R(H) \lra \I_p(H) \lra \I_{p/R}(H) \lra 0,\]
hence $H^0(Y, \I_p(H)) \simeq H^0(Y, \I_{p/R}(H))$. This means the base locus of $|\I_p(H)|$ is contained in the smooth surface $R$. But on $R \simeq K$, the system $|H|$ is very ample.

It remains to ensure that $\psi$ is divisorial. Let $S \subset \PN_3$ be a hyperplane, such that the intersection with $K$ is a plane quartic curve with a double point. Let $p$ be the point over the singularity. Then the strict transform of $S$ is $D$.

\vspace{0.2cm}

\noindent {\bf No.~6.}
By \cite{AG5}, Lemma~3.3.4, there exists a $Y = V_{2,2}$ with a line $Z$ whose splitting type of the normal bundle is $(2,-2)$. Using the notation of no.~5, take the point $p$ on $Z \cap R$. Then $-K_X$ is big and nef, since $p \in R$, and there is a one-dimensional family of lines through $p$, the strict transform is $D$.

\vspace{0.2cm}

\noindent {\bf No.~7.}
Let $x_0, \dots, x_4$ be homogeneous coordinates of $\PN_4$. Let $f_3 \in \KC[x_0, x_1, x_2]$ be homogeneous of degree $3$, defining a smooth elliptic curve $B$ in $\PN_2$, and let $q \in \KC[x_0, \dots, x_4]$ be a general quadric. Define 
 \[Y = \{f_3 + x_4q = 0\} \subset \PN_4.\] 
Then $Y$ is a smooth cubic in $\PN_4$, hence of type $V_{2,3}$. The intersection with the hyperplane $H = \{x_4 = 0\}$ is by construction the cone over the smooth elliptic curve $B$, which we denote by $D'$. The vertex of $D'$ is the point 
 \[p = [0:0:0:1:0] \in Y \cap H = D'.\] 
Define $\phi: Bl_p(\PN_4) \to \PN_4$ with exceptional divisor $E$ and $X: = Bl_p(Y)$ the strict transform of $Y$. Then $-K_X = 2(\phi^*\sO_{\PN_4}(1) - E)|_X$ is big and nef. The strict transform $D$ of $D'$ is a smooth ruled surface in $|\phi^*H -3E|$ on $X$, such that $-K_X$ is trivial on the ruling, i.e., $\psi$ contracts $D$ to the elliptic curve $B$. 
\end{proof}

%%%%%%%%%%%%%%%%%%%

\vfill

\appendix

\section{Tables}

\subsection{Mukai's list} \label{Mukaitab}
The following table can be found in \cite{Mukai}, it contains all anticanonically embedded Gorenstein Fano threefolds with canonical singularities, such that the anticanonical system does not admit a moving decomposition. Note that all of these threefolds are degenerations of smooth ones.

\vspace{0.5cm}

\noindent \begin{tabular}{l|p{11.5cm}|}
 $g$ & Anticanonical model \\\hline
 $3$ & $X'_4 = (4) \subset \PN_4$ a quartic.\\
 $4$ & $X'_6 = (2) \cap (3) \subset \PN_5$ a complete intersection.\\
 $5$ & $X'_8 = (2) \cap (2) \cap (2) \subset \PN_6$ a complete intersection.\\
 $6$ & $X'_{10} \subset \PN_7$ is a quadric hypersurface section of a quintic del Pezzo $4$-fold $V_5 \subset \PN_7$. $V_5 \subset \PN_7$ may be chosen as the cone over a quintic del Pezzo $3$-fold.\\
 $7$ & $X'_{12} \subset \PN_8$ is a linear section of a $10$-dimensional orthogonal Grassmannian variety. $[\Sigma_{12}^{10} \subset \PN_{15}] \cap H_1 \cap \dots \cap H_7$.\\
 $8$ & $X'_{14} \subset \PN_9$ is a linear section of an $8$-dimensional Grassmannian variety. $[G(2,6) \subset \PN_{14}] \cap H_1 \cap \dots \cap H_5$.\\
 $9$ & $X'_{16} \subset \PN_{10}$ is a linear section of a $6$-dimensional symplectic Grassmannian variety. $[\Sigma_{16}^6 \subset \PN_{13}] \cap H_1 \cap H_2 \cap H_3$.\\
 $10$ & $X'_{18} \subset \PN_{11}$ is a linear section of a $G_2$-variety. $[\Sigma_{18}^5 \subset \PN_{13}] \cap H_1 \cap H_2$.\\
 $12$ & $X'_{22} \subset \PN_{13}$ is isomorphic to a threefold $G(3,7,N) \subset \PN_{13}$ obtained from a non-degenerate $3$-dimensional subspace $N \subset \wedge^2\KC^7$.\\\hline
\end{tabular}

\newpage

\subsection{Del Pezzo fibrations} \label{tabdP}
We use the notation of section~\ref{delPezzo}. In the following table $d_B/g_B$ denotes genus and degree of the curve $B$, where the degree is by definition $H'\cdot B$ with $H'$ a general hyperplane section of $X'$. A $(k)$ means that here $X'$ is the image of $|{-}\frac{1}{k}K_X|$. The abbreviation {\em a.m.} stands for {\em anticanonical model}. Finally ``?'' indicates that the existence is still open. 

\vspace{0.5cm}

\noindent \begin{tabular}{c||c||c||p{6cm}|c||c|}
  No. & $({-}K_X)^3$ & $K_F^2$ & $X'$ & $d_B/g_B$ & $(\alpha, \beta)$\\\hline
  $1$ & $54$ & $9$ & $\PN(1^2,2^2)$ \hfill  $(3)$ & $1/0$ & $(\frac{1}{3}, 2)$\\\hline
  $2$ & $4$ & $4$ & $X'_4 \subset \PN_4$ or \newline 2:1 over $Q_3$, ram. quartic & $8/5$ & $(2,2)$\\
  $3$ & $2$ & $3$ & 2:1 over $\PN_3$, ram. sextic & $9/10$ & $(3,2)$\\\hline
  $4$ & $4$ & $2$ & 2:1 over $\PN(1^2,2^2)$ & $2/1$ & $(1,2)$\\
  $5$ & $6$ & $3$ & $X'_{2,3} \subset \PN_5$ & $3/1$ & $(1,2)$\\
  $6$ & $8$ & $4$ & $X'_{2,2,2} \subset \PN_6$ & $4/1$ & $(1,2)$\\
  $7$ & $10$ & $5$ & a.m. of $X \subset \PN(\sO_{\PN_1}(2) \oplus \sO_{\PN_1}^{\oplus 5})$ \hfill $(?)$ & $5/1$ & $(1,2)$\\
  $8$ & $12$ & $6$ & a.m. of $X \subset \PN(\sO_{\PN_1}(2) \oplus \sO_{\PN_1}^{\oplus 6})$ \hfill $(?)$ & $6/1$ & $(1,2)$\\
  $9$ & $16$ & $8$ & 2:1 over $\PN_3$, ram. quartic \hfill $(2)$ & $4/1$ & $(1,2)$\\\hline
  $10$ & $2$ & $2$ & 2:1 over $\PN_3$, ram. sextic & $1/0$ & $(1,1)$\\
  $11$ & $4$ & $4$ & $X'_4 \subset \PN_4$ or \newline 2:1 over $Q_3$, ram. quartic & $2/0$ & $(1,1)$\\ 
  $12$ & $8$ & $8$ & 2:1 over Ver.c. & $4/0$ & $(1,1)$\\ 
  $13$ & $2$ & $4$ & 2:1 over $\PN_3$, ram. sextic & $4/1$ & $(2,1)$\\\hline
  $14$ & $32$ & $8$ & $X'_2 \subset \PN(1^2, 2^3)$ \hfill $(2)$ & $2/0$ & $(\frac{1}{2},2)$\\ 
  $15$ & $16$ & $8$ & 2:1 over $\PN_3$, ram. quartic \hfill $(2)$ & $2/0$ & $(\frac{1}{2}, 1)$\\\hline
\end{tabular}

\vspace{0.8cm}

\subsection{Conic bundles} \label{tabcb}

We use the notation of section~\ref{conic}. As above, $d_B/g_B$ denotes genus and degree of the curve $B$; a $(2)$ means that here $X'$ is the image of $|{-}\frac{1}{2}K_X|$; the degree of the discriminant $\Delta$ is $d$. If $\Delta = 0$, then $(c_1,c_2)$ are the Chern classes of $\sF$, such that $X = \PN(\sF)$, else we denote with $(c_1, c_2)$ the Chern classes of $\E = \phi_*({-}K_X)$. By $S_1(B_4)$ we denote the secant variety to the rational normal curve of degree $4$. The abbreviation {\em a.m.} stands for {\em anticanonical model}.

\vspace{0.5cm}

\noindent \begin{tabular}{c||c||c|c||p{5cm}|c||c|}
  No. & $({-}K_X)^3$ & $d$ & $(c_1,c_2)$ & $X'$ & $d_B/g_B$ & $(\alpha, \beta)$\\\hline
  $1$ & $72$ & $0$ & $(-1,-2)$ & $\PN(1^3, 3)$ \hfill $(2)$ & -- & $(\frac{1}{2}, -3)$\\\hline
  $2$ & $48$ & $0$ & $(-1,1)$ & $X'_4 \subset \PN_4$ \hfill $(2)$ & $1/0$ & $(\frac{1}{2}, -2)$\\
  $3$ & $24$ & $0$ & $(-1,4)$ & $X'_3 \subset \PN_4$ \hfill $(2)$& $1/0$ & $(\frac{1}{2}, -1)$\\
  $4$ & $24$ & $0$ & $(-1,4)$ & $S_1(B_4) \subset \PN_4$ \hfill $(2)$ & $4/0$ & $(1, -2)$\\\hline
  $5$ & $10$ & $7$ & $(2,0)$ & a.m. of $X \subset \PN(\sO_{\PN_2}(2) \oplus \sO_{\PN_2}^{\oplus 2})$ & $1/0$ & $(1,-2)$\\
  $6$ & $12$ & $6$ & $(3,3)$ & a.m. of $X \subset \PN(\sE)$ as in 3.10. or \newline a.m. of $X \subset \PN(\sO_{\PN_2}(2) \oplus T_{\PN_2}(-1))$ & $2/0$ & $(1,-2)$\\\hline
  $7$ & $4$  & $8$ & $(1,0)$ & $X'_4 \subset \PN_4$ & $1/0$ & $(1,-1)$\\
  $8$ & $6$  & $6$ & $(3,6)$ & $X'_{2,3} \subset \PN_5$ & $2/0$ &
  $(1,-1)$\\\hline
  $9$ & $2$ & $10$ & $(-1,-2)$ & 2:1 over $\PN_3$, ram. sextic & -- &
  $(1,-1)$\\
  $10$ & $4$ & $8$ & $(1,0)$ & 2:1 over $Q_3$, ram. quartic & $1/0$ &
  $(1,-1)$\\
  $11$ & $2$ & $8$ & $(1,1)$ & 2:1 over $\PN_3$, ram. sextic & $3/0$ &
  $(2,-1)$\\
  $12$ & $8$ & $8$ & $(1,-2)$ & 2:1 over Ver. cone, ram. cubic & -- & $(1,-2)$\\\hline
\end{tabular}  

\newpage

\subsection{Blowdown to a curve}\label{tabdivcurve}
We use the notation of section~\ref{blowdowncurve}; $d/g_C$ and $d/g_B$ denote the degree and genus of $C$ and $B$, a ``--'' indicates that $D$ maps to a point. By $V_{r,d}$ we denote a Fano threefold of index $r$ and degree $d$; {\em a.m.} stands for {\em anticanonical model} and an ``(M)'' means that $X'$ can be found in table~\ref{Mukaitab}.

\vspace{0.1cm}

\noindent \begin{tabular}{c||c||c|p{1.1cm}||p{4.6cm}|c||c|}
  No. & $({-}K_X)^3$ & $Y$ & $d/g_C$ & $X'$ & $d/g_B$ & $(\alpha, \beta)$\\\hline
  $1$ & $2$ & $V_{1,8}$ & $2/0$ & 2:1 over $\PN_3$, ram. sextic & $5/2$ & $(2,3)$\\
  $2$ & $2$ & $V_{1,10}$ & $4/1$ & 2:1 over $\PN_3$, ram. sextic & $4/0$ & $(2,3)$\\
  $3$ & $4$ & $V_{1,10}$ & $2/0$ & $X'_4 \subset \PN_4$ or \newline 2:1 over $Q_3$, ram. quartic & $3/0$ & $(1,2)$\\ 
  $4$ & $6$ & $V_{1,16}$ & $4/0$ & $X'_{2,3} \subset \PN_5$ & $4/0$ & $(1,2)$\\
  $5$ & $6$  & $V_{2,2}$ & $2/0$ & $X'_{2,3} \subset \PN_5$ & $1/0$ & $(1,1)$\\
  $6$ & $8$  & $V_{2,3}$ & $4/1$ & $X'_{2,2,2} \subset \PN_6$ & $1/0$ & $(1,1)$\\
  $7$ & $2$ & $V_{2,4}$ & $10/6$ & 2:1 over $\PN_3$, ram. sextic & $5/1$ & $(5,3)$\\
  $8$ & $4$ & $V_{2,4}$ & $8/3$ & $X'_4 \subset \PN_4$ or \newline 2:1 over $Q_3$, ram. quartic & $4/0$ & $(3,2)$\\
  $9$ & $10$ & $V_{2,4}$ & $6/2$ & (M) & $1/0$ & $(1,1)$\\
  $10$ & $12$ & $V_{2,5}$ & $8/3$ & (M) & $1/0$ & $(1,1)$\\
  $11$ & $14$ & $V_{2,5}$ & $6/0$ & (M) & $2/0$ & $(1,1)$\\
  $12$ & $14$ & $V_{2,5}$ & $6/0$ & (M) & $8/0$ & $(2,2)$\\
  $13$ & $2$  & $Q_3$ & $12/11$ & 2:1 over $\PN_3$, ram. sextic & $6/3$ & $(8,3)$\\
  $14$ & $4$  & $Q_3$ & $10/6$ & $X'_4 \subset \PN_4$ or \newline 2:1 over $Q_3$, ram. quartic & $5/1$ & $(5,2)$\\
  $15$ & $8$  & $Q_3$ & $10/8$ & $X'_{2,2,2} \subset \PN_6$ & $1/0$ & $(2,1)$\\
  $16$  & $10$ & $Q_3$ & $8/3$ & (M) & $2/0$ & $(2,1)$\\ 
  $17$  & $26$ & $Q_3$ & $5/2$ & a.m. of $X \subset \PN(\sO_{Q_3}^{\oplus 2} \oplus \sO_{Q_3}(2))$ & $1/0$ & $(1,1)$\\
  $18$  & $28$ & $Q_3$ & $4/0$ & a.m. of $X \subset \PN(\sO_{Q_3}^{\oplus 3} \oplus \sO_{Q_3}(2))$ & $2/0$ & $(1,1)$\\
  $19$  & $4$  & $\PN_3$ & $10/11$ & $X'_4 \subset \PN_4$ & $10/6$ & $(10,3)$\\
  $20$  & $6$  & $\PN_3$ & $10/12$ & $X'_{2,3} \subset \PN_5$ & $1/0$ & $(3,1)$\\
  $21$  & $8$  & $\PN_3$ & $8/5$ & $X'_{2,2,2} \subset \PN_6$ & $2/0$ & $(3,1)$\\
  $22$  & $18$ & $\PN_3$ & $7/6$ & a.m. of $X \subset \PN(\sO_{\PN_3}^{\oplus 2} \oplus \sO_{\PN_3}(2))$, \newline $X = H_1 \cap H_2$ is a c.i. & $1/0$ & $(2,1)$\\
  $23$  & $20$ & $\PN_3$ & $6/3$ & a.m. of $X \subset \PN(\sO_{\PN_3}^{\oplus 3} \oplus \sO_{\PN_3}(2))$ & $2/0$ & $(2,1)$\\
  $24$  & $22$ & $\PN_3$ & $5/0$ & a.m. of $X \subset \PN(\sO_{\PN_3}^{\oplus 4} \oplus \sO_{\PN_3}(2))$ & $3/0$ & $(2,1)$\\\hline
  $25$ & $r(r{-}r')^2H^3$ & $V_{r,d}$ & $rr'H^3/$\newline $\frac{2+rr'^2H^3}{2}$ & (M) or hyperelliptic as in \ref{hyp} & -- & $(r', 1)$\\
\hline
 \end{tabular}

\subsection{Blowdown to a point}\label{tabpointpoint}
We use the notation of section~\ref{blowdownpoint}; $d/g_B$ denotes the degree and genus of $B$; $V_{r,d}$ is a (terminal) Fano threefold of index $r$ and degree $d$; a $(2)$ means that $X'$ is the image of $|{-}\frac{1}{2}K_X|$; {\em a.m.} stands for {\em anticanonical model}.

\vspace{0.1cm}

\noindent \begin{tabular}{c||c||c|c||p{5cm}|c||c|}
  No. & $({-}K_X)^3$ & $Y$ & $(E, -E|_E)$ & $X'$ & $d/g_B$ & $(\alpha, \beta)$\\\hline
  $1$ & $2$ & $V_{1,4}$ & $(Q, \sO(1))$ & 2:1 over $\PN_3$, ram. sextic & $2/0$ & $(1,2)$\\
  $2$ & $2$ & $V_{1,4}$ & $(Q, \sO(1))$ & 2:1 over $\PN_3$, ram. sextic & $8/9$ & $(2,4)$\\
  $3$ & $4$ & $V_{1,6}$ & $(Q, \sO(1))$ & $X'_4 \subset \PN_4$ or \newline 2:1 over $Q_3$, ram. quartic & $6/4$ & $(1,3)$\\
  $4$ & $2$ & $V_{1,20}$ & $(\PN_2, \sO(2))$ & 2:1 over $\PN_3$, ram. sextic & $5/6$ & $(\frac{1}{2}, \frac{5}{2})$\\
  $5$ & $8$ & $V_{2,2}$ & $(\PN_2, \sO(1))$ & 2:1 over Ver.c., ram. cubic & $2/0$ &$(1,2)$\\
  $6$ & $8$ & $V_{2,2}$ & $(\PN_2, \sO(1))$ & 2:1 over Ver.c., ram. cubic & $8/3$ & $(2,4)$\\
  $7$ & $16$ & $V_{2,3}$ & $(\PN_2, \sO(1))$ & a.m. of $X \subset \PN(\sO_{\PN_3} \oplus \sO_{\PN_3}(1))$ \hfill $(2)$ & $6/1$ & $(1,3)$\\\hline
 \end{tabular}

%%%%%%%%%%%%%%%%%


\begin{thebibliography}{Mum69}
\bibitem[BS95]{BS95} M. Beltrametti, A.J. Sommese: The adjunction theory of complex projective varieties. de Gruyter 1995.
\bibitem[B07]{Bertini} E. Bertini: Introduzione alla geometria proiettiva degli iperspazi. E. Spoerri, Pisa, 1907.
\bibitem[CC97]{Catanese} G. Casnati, F. Catanese: Even sets of nodes are bundle symmetric. J. Diff. Geom. {\bf 47}, 237-256 (1997); erratum ibid. {\bf 50}, 415 (1998).
\bibitem[CSP04]{Cheltsov} I. Cheltsov, C. Shramov, V. Przyjalkowski: Hyperelliptic and trigonal Fano threefolds. math.AG/0406143.
\bibitem[Fu90]{Fu90} T. Fujita: Classification theories of polarized varieties. London Math. Soc. Lect. Notes Ser. {\bf 155} (1990).
\bibitem[Gu82]{Gushel} N.P. Gushel: Fano varieties of genus $6$. Math. USSR, Izv. {\bf 21}, 445-459 (1983).
\bibitem[Hr92]{Harris} J. Harris: Algebraic Geometry. Springer 1992.
\bibitem[Ha77]{Ha} R. Hartshorne: Algebraic Geometry. Springer 1977.
\bibitem[IF89]{IF} A.R. Iano-Fletcher: Working with weighted complete intersections. Bonn preprint MPI/89-35 (1989).
\bibitem[IS00]{IS} A. Iliev, C. Schuhmann: Tangent scrolls in prime Fano threefolds. Kodai Math. J. {\bf 23}, 411-431 (2000).
\bibitem[I78]{Isk} V.A. Iskovskikh: Fano $3$-folds I, II. Math. USSR, Izv. {\bf 11}, 485-527 (1977); {\bf 12}, 469-506 (1978).
\bibitem[I80]{Iskantican} V.A. Iskovskikh: Anticanonical models of three-dimensional algebraic varieties. J. Soviet Math. {\bf 13}, 745-814 (1980).
\bibitem[IP99]{AG5} V.A. Iskovskikh, Yu.G. Prokhorov: Algebraic Geometry V: Fano varieties. Springer 1999.
\bibitem[JR03]{almc} P. Jahnke, I. Radloff: Fano threefolds with sections in $\Omega_V^1(1)$. math.AG/0310390.
\bibitem[JR04]{genprep} P. Jahnke, I. Radloff: Gorenstein Fano threefolds with base points in the anticanonical system. math.AG/0404156, to appear in Comp. Math.
\bibitem[KM98]{KoMo} J. Koll\'ar, S. Mori: Birational Geometry of Algebraic Varieties. Camb. Univ. Press 1998.
\bibitem[Me99]{Mella} M. Mella: Existence of good divisors on Mukai varieties. J. Alg. Geom. {\bf 8}, 197-206 (1999).
\bibitem[Mi03]{Minagawa} T. Minagawa: Global smoothing of singular weak Fano $3$-folds. J. Math. Soc. Japan {\bf 55}, 695-711 (2003).
\bibitem[Mo82]{Mori} S. Mori: Threefolds whose canonical bundles are not numerically effective. Ann. Math. {\bf 116}, 133-176 (1982). 
\bibitem[Muk95]{Mukai} S. Mukai: New developments in the theory of Fano $3$-folds: Vector bundle method and moduli problem. Sugaku {\bf 47}, 125-144 (1995); Sugaku Expositions {\bf 15}, 125-150 (2002).
\bibitem[MU83]{MU} S. Mukai, H. Umemura: Minimal rational threefolds. Lect. Notes Math. {\bf 1016}, 490-518 (1983).
\bibitem[Mum69]{Mumford} D. Mumford: Varieties defined by quadric equations. C.I.M.E. (Varenna), 29-100 (1969).
\bibitem[Pa98]{Paoletti} R. Paoletti: The K\"ahler cone in families of quasi-Fano threefolds. Math. Z. {\bf 227}, 45-68 (1998).
\bibitem[Pr90]{ProkhEx} Yu.G. Prokhorov: On exotic Fano varieties. Moscow Univ. Math. Bull. {\bf 45}, No.4, 36-38 (1990).
\bibitem[Pr03]{Prokh} Yu.G. Prokhorov: A remark on Fano threefolds with canonical singularities. math.AG/0312012, to appear in Proc. Fano Conf.
\bibitem[R83]{Reid} M. Reid: Projective morphisms according to Kawamata. Preprint (1983).
\bibitem[Shi89]{Shin} K.-H. Shin: $3$-dimensional Fano Varieties with Canonical Singularities. Tokyo J. Math. {\bf 12}, 375-385 (1989).
\bibitem[Sho80]{Shokurov} V.V. Shokurov: Smoothness of a general anticanonical divisor on a Fano variety. Math. USSR, Izv {\bf 14}, 395-405 (1980). 
\bibitem[SW90]{SW90} M. Szurek, J. A. Wi\'sniewski: On Fano manifolds, which are $\PN^k$-bundles over $\PN^2$. Nagoya Math. J. {\bf 120}, 89-101 (1990).
\bibitem[T72]{Tyurin} A.N. Tyurin: Five lectures on three-dimensional varieties. Russ. Math. Surv. {\bf 27}, No.2, 1-53 (1972).
\bibitem[W92]{Wilson} P.M.H. Wilson: The K\"ahler cone on Calabi-Yau threefolds. Inv. Math. {\bf 107}, 561-583 (1992); {\em Erratum}: {\bf 114}, 231-233 (1993).
\bibitem[W97]{Wilson2} P.M.H. Wilson: Symplectic deformations of Calabi-Yau threefolds. J. Diff. Geom. {\bf 45}, 611-637 (1997). 
\end{thebibliography}
\end{document}